%% file: WEK.tex
\documentclass[12pt,reqno]{amsart}
\usepackage[utf8]{inputenc}
\usepackage{fancyhdr}
\usepackage{mathabx}
\input{Macros.tex}

\usepackage{amsrefs}
\usepackage[hidelinks]{hyperref}
\makeatletter
\@namedef{subjclassname@2020}{\textup{2020} Mathematics Subject Classification}
\makeatother

\begin{document}
\title{Weighted Erd\H{o}s--Kac Theorems via Computing Moments}
\author{Kai (Steve) Fan}
\address{Department of Mathematics, Dartmouth College, Hanover, NH 03755, USA}
\email{steve.fan.gr@dartmouth.edu}
\setlength{\footskip}{8mm}
\begin{abstract}
By adapting the moment method developed by Granville and Soundararajan \cite{GS}, Khan, Milinovich and Subedi \cite{KMS} recently obtained a weighted version of the Erd\H{o}s--Kac theorem for $\omega(n)$ with multiplicative weight $d_k(n)$, where $\omega(n)$ denotes the number of distinct prime divisors of a positive integer $n$, and $d_k(n)$ is the $k$-fold divisor function with $k\in\N$. In this paper, we generalize their method to study the distribution of additive functions $f(n)$ weighted by nonnegative multiplicative functions $\alpha(n)$ in a wide class. In particular, we establish uniform asymptotic formulas for the moments of $f(n)$ with suitable growth rates. We also prove a qualitative result on the moments which extends a theorem of Delange and Halberstam \cite{DH}. As a consequence, we obtain a weighted analogue of the Kubilius--Shapiro theorem.
\end{abstract}

\maketitle
\tableofcontents
\section{Introduction}\label{S:Intro}
The celebrated Erd\H{o}s--Kac theorem, first proved by Erd\H{o}s and Kac \cite{EK} in 1940, states that if $\omega(n)$ denotes the number of distinct prime divisors of a positive integer $n$, then 
\begin{equation}\label{S1Equ1}
\lim\limits_{x\rightarrow\infty}\frac{1}{x}\cdot\#\left\{n\leq x\colon \frac{\omega(n)-\log\log n}{\sqrt{\log\log n}}\leq V\right\}=\Phi(V)
\end{equation}
for any given $V\in\R$, where
\[\Phi(V)\colonequals\frac{1}{\sqrt{2\pi}}\int_{-\infty}^{V}e^{-v^2/2}\,dv\]
is the cumulative distribution function of the standard Gaussian distribution. This statistical result is a direct upgrade of an earlier theorem of Hardy and Ramanujan on the \emph{normal order} of $\omega$ (see \cite{HR} and \cite[Theorem 431]{HW}), which asserts that given any $\epsilon>0$, the inequality $|\omega(n)-\log\log n|<\epsilon\log\log n$ holds for all but $o(x)$ values of $n\leq x$. In fact, Erd\H{o}s and Kac proved in the same paper a more general result in which the function $\omega$ is replaced by any strongly additive function $f$ that is bounded on primes and admits an unbounded ``variance" $\sum_{p\leq x}f(p)^2/p$. Recall that an arithmetic function $f\colon\N\to\C$ is \emph{additive} if $f(mn)=f(m)+f(n)$ for all positive integers $m,n\in \N$ with $\gcd(m,n)=1$, and it is \emph{strongly additive} if it also satisfies $f(p^{\nu})=f(p)$ for all prime powers $p^{\nu}$. Thus, strongly additive functions are completely determined by their values at primes, which makes them a particularly nice subclass of additive functions. In fact, it can be shown that (\ref{S1Equ1}) also holds for $\Omega(n)$ in place of $\omega(n)$, where $\Omega(n)$ denotes the total number of prime factors of $n$, counting multiplicity, by exploiting the fact that $\omega$ and its cousin $\Omega$ do not differ very much on average:
\begin{equation}\label{S1Equ:Omea-omega}
\sum_{n\le x}(\Omega(n)-\omega(n))=O(x).
\end{equation}
\par The original proof of the Erd\H{o}s--Kac theorem \eqref{S1Equ1} combines the central limit theorem with Brun's sieve and is quite complicated. Later, LeVeque \cite[Theorem 1]{LeVeque} introduced some modifications and obtained a quantitative version of (\ref{S1Equ1}) with a rate of convergence $O(\log\log\log x/\sqrt{\log\log x})$. 
Using a deep analytic approach, R\'{e}nyi and Tur\'{a}n \cite{RT} improved upon LeVeque's result with a rate of convergence $O(1/\sqrt{\log\log x})$, which is best possible in the sense that it cannot be improved to $o(1/\sqrt{\log\log x})$ without loss of uniformity in $V$. 
\par A third approach to \eqref{S1Equ1}, first suggested by Kac \cite{Kac}, examines the moments of $\omega$. In probability theory, the moments of a random variable $X$ often provides valuable information about its distribution. For example, an application of Markov's inequality yields $\PP(|X|>c)\le\mathbb{E}|X|^k/c^k$. Given all the moments $\mathbb{E}|X|^k<\infty$, one may select $k$ that minimizes this tail estimate. If $X$ happens to obey a Gaussian law, then it is completely determined by its moments, a direct consequence of \cite[Theorem 30.1]{BillingsleyText} or \cite[Theorem 3.3.26]{Dur}. Consequently, by \cite[Theorem 30.2]{BillingsleyText}, one reduces the proof of (\ref{S1Equ1}) to that of the asymptotic formula 
\begin{equation}\label{S1Equ2}
\frac{1}{x}\sum_{n\leq x}\left(\omega(n)-\log\log x\right)^m=(\mu_m+o(1))(\log\log x)^{\frac{m}{2}}
\end{equation}
for every $m\in\N$. Here $\mu_m$ is the $m$th moment of a standard Gaussian distribution given by
\[\mu_m=\begin{cases}
	~m!/m!!,&\text{~~~~~if $2\mid m$},\\
	~0,&\text{~~~~~otherwise}.
\end{cases}\]
The case $m=1$ follows from Mertens' theorem \cite[Theorem 427]{HW}, and the case $m=2$ was settled by Tur\'an \cite{Turan}. Early proofs of (\ref{S1Equ2}) for all $k$ are due to Delange \cite{Delange1} in 1953 and Halberstam \cite{Hal1} in 1955, both of which are very complicated. 
Later, Delange \cite{Delange2} provided an elementary proof of (\ref{S1Equ2}) for strongly additive functions.
Halberstam's proof was also simplified and rendered more transparent by Billingsley \cite{Billingsley} in 1969.
In 2007, Granville and Soundararajan \cite{GS} derived asymptotic formulas for the moments which hold uniformly in the range $m\leq(\log\log x)^{1/3}$. Their method, which may be viewed to some extent as a clever and efficient repacking of Billingsley's, is so flexible that it can be modified to study the distribution of values of additive functions in a rather general sieve-theoretic framework.
\par More generally, one can study the distribution of  values of $\omega(n)$ weighted by certain nonnegative multiplicative functions $\alpha(n)$. For instance, Elliott \cite{Elliott2} showed, based on the Landau--Selberg--Delange method, that 
\[\lim\limits_{x\to\infty}\left(\sum_{n\leq x}d(n)^c\right)^{-1}\sum_{\substack{n\leq x\\\omega(n)\leq 2^c\log\log x+V\sqrt{2^c\log\log x}}}d(n)^c=\Phi(V)\]
for any given $c\in\R$ and $V\in\R$, where $d(n)$ denotes the number of positive divisors of $n$. 
Building on the method of Granville and Soundararajan, Khan, Milinovich and Subedi \cite{KMS} recently proved an analogue for the weight $d_k(n)$ with mean $k\log\log x$ and variance $k\log\log x$, where
$d_k(n)\colonequals\#\left\{(a_1,...,a_k)\in\N^k\colon a_1\cdots a_k=n\right\}$.
There is now a vast literature on weighted versions of the Erd\H{o}s--Kac theorem with general weights, including the early work of Alladi \cite{Alla} and the more recent works of Elboim and Gorodetsky \cite{EG} and Tenenbaum \cites{Tenen1,Tenen2}. Alladi made use of Halberstam's approach to prove weighted Erd\H{o}s--Kac type results for strongly additive functions with the weights being the characteristic functions of the subsets of $\N$ which satisfy certain sieve type conditions. On nonnegattive multiplicative weights $\alpha$, Elboim and Gorodetsky \cite[Theorem 1.1]{EG} generalized Billingsley's proof \cite{Billingsley} to handle the distribution of $\Omega(n)$ weighted by those $\alpha$ having constant mean values and satisfying certain growth conditions, while Tenenbaum's result \cite[Corollary 2.5]{Tenen1} was proved by means of charactersitic functions and allows for general additive functions and a large class of multiplicative weights $\alpha$ with the property that $\alpha(p)=O(p^{\sigma_0-1})$ for some constant $\sigma_0>0$.
\par The main purpose of this paper is to establish weighted versions of the Erd\H{o}s--Kac theorem by pushing the method of moments of Granville and Soundararajan to its limit. Our work is the first to apply this method to prove weighted Erd\H{o}s--Kac theorems with general additive functions and multiplicative weights. 
We obtain uniform estimates for moments of strength comparable to that of the original estimate of Granville and Soundararajan for $\omega(n)$. With our emphasis on the strength of the method, we have refrained ourselves from pursuing the most general theorems at the risk of complicating our exposition. Despite this compromise, our results retain some of its own advantages over the results in \cite{EG,Tenen1,Tenen2}. Our approach is elementary and flexible, and it can be applied to handle certain arithmetic functions of special interests which were only studied previously by different methods. Some examples are discussed in the comment below Corollary \ref{corWeightedKS}.\\
\\
\noindent {\bf Definitions and notation.} We introduce some terminologies and notation that will be adopted throughout this paper without further clarification. Given any real or complex valued functions $f(x)$ and $g(x)$ with a common domain $ \mathcal{D}\subseteq \R$ , we shall use Landau's big-$O$ notation $f(x)=O(g(x))$ and Vinogradov's notation $f(x)\ll g(x)$ interchangeably to mean that there exists an absolute constant $C>0$ such that $|f(x)|\le C|g(x)|$ for all $x\in \mathcal{D}$. Likewise, we shall use the notation $f(x)\gg g(x)$ interchangeably with $g(x)=O(f(x))$ . If $f(x)=O(g(x))$ and $g(x)=O(f(x))$ hold simultaneously, then we adopt the short-hand notation $f(x)\asymp g(x) $. If $\mathcal{D}$ contains a neighborhood of $\infty$, then we write $f(x)=o(g(x))$ when $f(x)/g(x)\to0$ as $x\to\infty$ and $f(x)\sim g(x)$ when $f(x)/g(x)\to1$ as $x\to\infty$. We shall occasionally make use of the characteristic function $\epsilon_{a,b}$ of the condition $a\ne b$
for any $a,b\in\R$. Equivalently, $\epsilon_{a,b}=1-\delta_{a,b}$, where $\delta_{a,b}$ is the Kronecker delta function. 
\par Throughout, the letter $p$ always denotes a prime, and we write $\pi(x)$ for the prime counting function, namely, $\pi(x)=\sum_{p\leq x}1$. For any $x\in\R$, we write $\lfloor x\rfloor$ for the integer part of $x$ and $\lceil x\rceil$ for the least integer $\ge x$. For every $n\in\N$, denote by $P^-(n)$ and $P^+(n)$ the least and the greatest prime factor of $n$, respectively, with the convention that $P^-(1)=\infty$ and $P^+(1)=1$. We say that $n\in\Z\setminus\{0\}$ is {\it squareful} if for any prime $p\mid n$, one has $p^2\mid n$. Given any prime power $p^{\nu}$, the relation $p^{\nu}\parallel n$ means that $p^{\nu}\mid n$ but $p^{\nu+1}\nmid n$. In addition, we denote by $R_n$ the {\it radical} of $n$, i.e., $R_n\colonequals \rad(n)=\prod_{p\mid n}p$. Finally, we write $\binom{m}{m_1,...,m_k}\colonequals\frac{m!}{m_1!\cdots m_k!}$ for the multinomial coefficient of shape $(m_1,...,m_k)$ of size $m=m_1+\cdots +m_k$.

\medskip
\section{Main Results}\label{S:MR}
The weights $\alpha\colon\N\to \R_{\ge0}$ that we shall consider in this paper form a nice subclass $\mathcal{M}^{\ast}$ of nonnegative multiplicative functions, nice in the sense that there exist constants $A_0,\beta,\sigma_0>0$, $\vartheta_0\ge0$, $\varrho_0\in[0,1)$ and $r\in(0,1)$, such that the following conditions hold:
\begin{alignat}{2}
&\text{(i)~~}&&\alpha(p^{\nu})\ll p^{(\varrho_0+\sigma_0-1)\nu},\label{thmMT2Equ0}\\
&\text{(ii)~~}&&\sum_{p\le x}\frac{\alpha(p)\log p}{p^{\sigma_0-1}}=\beta x+O\left(\frac{x}{(\log x)^{A_0}}\right),\label{thmMT2Equ1}\\
&\text{(iii)~~}&&\sideset{}{'}\sum_{p}\left(\frac{\alpha(p)^2}{p^{2(r+\sigma_0-1)}}+\sum_{\nu\ge2}\frac{\alpha(p^{\nu})}{p^{(r+\sigma_0-1)\nu}}\right)<\infty,\label{thmMT2Equ2}\\
&\text{(iv)~~}&&\sum_{\nu\ge1}\frac{\nu\alpha(p^{\nu})}{p^{\sigma_0\nu}}\ll\frac{(\log\log(p+1))^{\vartheta_0}}{p},\label{thmMT2Equ3}
\end{alignat}
where the sum $\sum_{p}^{'}$ is over all but finitely many primes $p$. 
It is not hard to verify that $\mathcal{M}^{\ast}$ is closed under Dirichlet convolution.
Despite some overlaps between our class $\mathcal{M}^{\ast}$ and the class of multiplicative functions studied by Elboim and Gorodetsky \cite{EG}, neither of them strictly contains the other. On the one hand, the multiplicative function $\alpha$, defined by $\alpha(p)=1$ for all primes $p$ and $\alpha(p^{\nu})=p^{\nu/3}$ for all prime powers $p^{\nu}$ with $\nu\ge2$, falls into $\mathcal{M}^{\ast}$ but is not covered by the first part of \cite[Theorem 1.1]{EG}. On the other hand, the condition (iv) implies the more restrictive growth condition $\alpha(p)/p^{\sigma_0-1}\ll(\log\log(p+1))^{\vartheta_0}$, which is not required in \cite[Theorem 1.1]{EG}.
\par The class $\mathcal{M}^{\ast}$ contains many familiar multiplicative functions, including the $\kappa$-fold divisor function $d_{\kappa}(n)$ for $\kappa>0$, the sum-of-divisors function $\sigma_{\lambda}(n)$ for $\lambda>-1$, Euler's totient function $\varphi(n)$, the characteristic function $\mu(n)^2$ of square-free numbers, and the function $r_2(n)/4$, 
where $\mu(n)$ is the M\"{o}bius function and $r_2(n)\colonequals\#\{(a,b)\in\Z^2\colon n=a^2+b^2\}$. Less obvious examples include $\rho_g(n)$, which denotes the number of zeros of a nonconstant irreducible polynomial $g\in\Z[x]$ in $\Z/n\Z$, and Ramanujan's $\tau$-function $\tau(n)$, which may be defined as the $n$th Fourier coefficient of the modular discriminant. We leave the verification of these claims to the interested reader.
\par Let $\alpha\in\mathcal{M}^{\ast}$ with parameters $A_0,\beta,\sigma_0,\vartheta_0,\varrho_0,r$ and set
\[S(x)=S_{\alpha}(x)\colonequals\sum_{n\leq x}\alpha(n).\] 
For any additive function $f\colon\N\to\R$, define 
\begin{align*}
A(x)&=A_{\alpha,f}(x)\colonequals\sum_{p\leq x}\alpha(p)\frac{f(p)}{p^{\sigma_0}},\\
B(x)&=B_{\alpha,f}(x)\colonequals\sum_{p\leq x}\alpha(p)\frac{f(p)^2}{p^{\sigma_0}}.
\end{align*}
One may think of $n$ as a random variable defined on the sample space $\N\cap[1,x]$ with the natural probability measure induced by $\alpha$, that is to say, $\text{Prob}(n=k)=\alpha(k)/S(x)$ for every $k\in\N\cap[1,x]$. We shall show, by estimating the weighted $m$th moment defined by
\[M(x;m)=M_{\alpha,f}(x;m)\colonequals S(x)^{-1}\sum_{n\leq x}\alpha(n)(f(n)-A(x))^m\]
for every $m\in\N$, that for certain additive functions $f$, the limiting distribution of the normalization $(f(n)-A(x))/\sqrt{B(x)}$ is standard Gaussian. To state our results in a coherent manner, we set $\chi_m\colonequals (1+(-1)^m)/2$, the characteristic function of even integers, and
\[C_m\colonequals\frac{m!}{2^{m/2}\Gamma(m/2+1)}\]
for $m\in\N$, where $\Gamma$ is the Gamma function. One quickly notes that $C_m=\mu_m=(m-1)!!$ for $m$ even. The numbers $C_m$ play a nonnegligible role in the error terms of our uniform estimates for $M(x;m)$. Our first result is the following theorem. 
\begin{thm}\label{thmMT2}
Let $f\colon\N\to\R$ be a strongly additive function with $|f(p)|\le M$ for all $p$, where $M>0$ is constant, and let $\alpha\in\mathcal{M}^{\ast}$ with parameters $A_0,\beta,\sigma_0,\vartheta_0,\varrho_0,r$. If $\beta=1$ and $0<h_0<(3/2)^{2/3}$ is arbitrary, and if $B(x)\to\infty$ as $x\to\infty$, then we have
\[M(x;m)=C_mB(x)^{\frac{m}{2}}\left(\chi_m+O\left(\frac{Mm^{\frac{3}{2}}}{\sqrt{B(x)}}\right)\right)\]
uniformly for all sufficiently large $x$ and all $1\le m\le h_0(B(x)/M^2)^{1/3}$. If $\beta\ne1$ and if $B(x)/(\log\log\log x)^2\to \infty$ as $x\to\infty$, then we have
\[M(x;m)=C_mB(x)^{\frac{m}{2}}\left(\chi_m+O\left(\frac{Mm^{\frac{3}{2}}\log\log\log x}{\sqrt{B(x)}}\right)\right)\]
uniformly for all sufficiently large $x$ and all $1\le m\ll B(x)^{1/3}/(\log\log\log x)^{2/3}$. The implicit constants in the error terms of both asymptotic formulas above depend at most on the explicit and implicit constants in the hypotheses except for $M$.
\end{thm}

\begin{rmk}\label{rmk1.1}
It may be worth pointing out that as in Theorem \ref{thmMT2}, the implicit constants in the estimates appearing in the rest of the paper depend at most on the explicit and implicit constants in the hypotheses unless stated otherwise. 
\end{rmk}

\par In the case where $\alpha=1$ and $f=\omega$, we recover \cite[Theorem 1]{GS} with a slightly wider range $1\le m\le h_0(\log\log x)^{1/3}$ compared to the original range $1\le m\le (\log\log x)^{1/3}$. 
Though Theorem \ref{thmMT2} is formulated for strongly additive functions, similar things can be said about additive functions whose values at prime powers do not grow too rapidly and are hence not expected to contribute very much. A simple example of such functions is $\Omega(n)$. Since $\Omega(p^{\nu})=\nu$ for all $p^{\nu}$, one can show, by establishing a weighted version of (\ref{S1Equ:Omea-omega}), that $\Omega(n)$ does not differ from its cousin $\omega(n)$ very much for most values of $n$, and so they are expected to have the same distribution. More generally, we shall prove the following variant of Theorem \ref{thmMT2} for additive functions. For simplicity's sake, we shall focus on a subclass of the multiplicative functions in $\mathcal{M}^{\ast}$.

\begin{thm}\label{thmMT2V}
Let $f\colon\N\to\R$ be an additive function such that  $f(p^{\nu})\le M\nu^{\kappa}$ for all prime powers $p^{\nu}$, where $M>0$ and $\kappa\ge0$ are constant. Let $\alpha\colon\N\to\R_{\geq0}$ be a multiplicative function, and suppose that there exist constants $A_0,\beta,\sigma_0>0$, $\vartheta_0\ge0$, $\varrho_0\in[0,1/2)$  and $\lambda\in(0,2^{1-2\varrho_0})$, such that $\alpha(n)$ satisfies \emph{(\ref{thmMT2Equ1})}, \emph{(\ref{thmMT2Equ3})}, and the condition that $\alpha(p^{\nu})=O((\lambda p^{\varrho_0+\sigma_0-1})^{\nu})$ for all prime powers $p^{\nu}$. If $\beta=1$ and $0<h_0<(3/2)^{2/3}$ is arbitrary, and if $B(x)\to\infty$ as $x\to\infty$, then we have
\[M(x;m)=C_mB(x)^{\frac{m}{2}}\left(\chi_m+O\left(\frac{Mm^{\kappa+\frac{3}{2}}}{\sqrt{B(x)}}\right)\right)\]
uniformly for all sufficiently large $x$ and all $m\in\N$ satisfying $m\le h_0(B(x)/M^2)^{1/3}$ and $m\ll B(x)^{1/(2\kappa+3)}$. If $\beta\ne1$ and if $B(x)/(\log\log\log x)^2\to \infty$ as $x\to\infty$, then we have
\[M(x;m)=C_mB(x)^{\frac{m}{2}}\left(\chi_m+O\left(\frac{Mm^{\frac{3}{2}}\left(\log\log\log x+m^{\kappa}\right)}{\sqrt{B(x)}}\right)\right)\]
uniformly for all sufficiently large $x$ and all 
\[1\le m\ll \min\left(B(x)^{1/(2\kappa+3)},\frac{B(x)^{1/3}}{(\log\log\log x)^{2/3}}\right).\]
The implicit constants in the error terms of both asymptotic formulas above depend at most on the explicit and implicit constants in the hypotheses except for $M$.
\end{thm}

It is easy to see that if $\alpha(p^{\nu})=O((\lambda p^{r_0+\sigma_0-1})^{\nu})$ for all prime powers $p^{\nu}$, where $\sigma_0>0$, $r_0\in[0,1/2)$ and $\lambda\in(0,2^{1-2r_0})$ are given constants, then conditions (i) and (iii) are automatically fulfilled with any fixed $\max(r_0+\log_2\lambda,0)\le\varrho_0<1$, $r_0+\max(1/2,\log_2\lambda)<r<1$, and the same parameter $\sigma_0$. Indeed, we shall derive Theorem \ref{thmMT2V} as a corollary of Theorem \ref{thmMT2}.
\par Let $g\in\Z[x]$ be a nonconstant irreducible polynomial, and recall that for every $n\in\N$, $\rho_g(n)$ denotes the number of zeros of $g$ in $\Z/n\Z$. More generally, if $g\in\Q[x]$ is a nonconstant irreducible polynomial, we may extend the definition above by setting $\rho_g(n)=0$ if $\gcd(n,c_g)>1$, where $c_g\in\N$ is the least positive integer such that $c_gg(x)\in\Z[x]$, and insisting that $\rho_g(n)$ be the number of zeros of $g(x)$ (or equivalently, $c_gg(x)$) in $\Z/n\Z$ when $\gcd(n,c_g)=1$. Extended this way with the convention that $\rho_g(1)=1$, the function $\rho_g(n)$ remains multiplicative. By \cite[Lemma 1]{Hal3}, $\rho_g$ is bounded on prime powers and satisfies
\[\sum_{p\le x}\frac{\rho_g(p)}{p}=\log\log x+M_{\rho_g}+O\left(\frac{1}{\log x}\right).\]
Given a strongly additive function $f\colon\N\to\R$, we define 
\begin{align*}
A_{f,g}(x)&\colonequals\sum_{p\leq x}\rho_g(p)\frac{f(p)}{p},\\
B_{f,g}(x)&\colonequals\sum_{p\leq x}\rho_g(p)\frac{f(p)^2}{p}.
\end{align*}
For simplicity's sake, suppose that $g(\N)\subseteq\N$. In the case $g\in\Z[x]$, Halberstam \cite[Theorem 3]{Hal2} showed that if $B_{f,g}(x)\to \infty$ as $x\to\infty$, and if  $f(p)=o\left(\sqrt{B_{f,g}(p)}\right)$\footnote{Halberstam \cite{Hal2} wrote that for $g(x)=x$ this pair of conditions contain the condition that $f(p)=o((\log p)^{\epsilon})$ for every given $\epsilon>0$. However, this claim is incorrect. In fact, a simple counterexample may be constructed as follows. Let $\mathcal{P}$ be an arbitrary infinite subset of odd primes such that $\sum_{p\in\mathcal{P}}1/p<\infty$, and put $\mathcal{P}(x)\colonequals\mathcal{P}\cap[3,x]$. Define $f(p)=\sqrt{\log\log p}$ for $p\in\mathcal{P}$ and $f(p)=1$ for $p\notin\mathcal{P}$. From partial summation it follows that $\sum_{p\in\mathcal{P}(x)}f(p)^2/p=o(\log\log x)$. Then one sees readily that $f(p)=o((\log p)^{\epsilon})$ for any given $\epsilon>0$, while $f(p)\sim \sqrt{B(p)}$ for large $p\in\mathcal{P}$.}, then given $m\in\N$,
\[\frac{1}{x}\sum_{n\leq x}\left(f(g(n))-A_{f,g}(x)\right)^m=(\mu_m+o(1))B_{f,g}(x)^{\frac{m}{2}}\]
holds. Under the stronger condition $f(p)=O(1)$, Theorem \ref{thmMT2} leads to a weighted version of this result in the case $g(n)=n$. The remaining cases are captured by the following theorem.
\begin{thm}\label{thmWeightedHal}
Let $f\colon\N\to\R$ be a strongly additive function with $|f(p)|\le M$ for all $p$, where $M>0$ is constant, and let $g\in\Q[x]$ be a nonconstant irreducible polynomial such that $g(0)\ne0$ and $g(\N)\subseteq\N$. Let $\alpha\in\mathcal{M}^{\ast}$ with parameters $A_0,\beta,\sigma_0,\vartheta_0,\varrho_0,r$, and fix $0<h_0<(3/2)^{2/3}$ and $h_0'>0$. For any $q\in\N$ and $a\in\Z$ coprime to $q$, define 
\[\Delta_{\alpha}(x;q,a)\colonequals \sum_{\substack{n\le x\\n\equiv a\,(\emph{mod}\, q)}}\alpha(n)-\frac{1}{\varphi(q)}\sum_{\substack{n\le x\\\gcd(n,q)=1}}\alpha(n).\]
If there exist a constant $\epsilon_0>0$ and a function $\delta(x)\in(0,1]$ with $\delta(x)^2B_{f,g}(x)\to\infty$ as $x\to\infty$, such that 
\begin{equation}\label{Equ:Delta_{alpha}(x;q,a)}
\sum_{\substack{q:\,\omega(q)\le m\\ P^+(q)\le x^{\delta(x)/m}}}\mu(q)^2\sum_{a\in\mathcal{Z}^{\ast}_g(q)}|\Delta_{\alpha}(x;q,a)|\ll S(x)\exp\left({-(\log\log x)^{1/3+\epsilon_0}}\right)
\end{equation}
uniformly for all sufficiently large $x$ and all 
\begin{equation}\label{Equ:mh_0h_1}
1\le m\le\min\left(h_0(B_{f,g}(x)/M^2)^{1/3},h_0'(\delta(x)^2B_{f,g}(x)/M^2)^{1/3}\right),
\end{equation}
where $\mathcal{Z}^{\ast}_g(q)\colonequals\{n\in(\Z/q\Z)^{\times}\colon g(n)\equiv0\,(\emph{mod}\, q)\}$ is the zero locus of $g$ in $(\Z/q\Z)^{\times}$, then
\[S(x)^{-1}\sum_{n\leq x}\alpha(n)\left(f(g(n))-A_{f,g}(x)\right)^m=C_mB_{f,g}(x)^{\frac{m}{2}}\left(\chi_m+O\left(\frac{Mm^{\frac{3}{2}}}{\delta(x)\sqrt{B_{f,g}(x)}}\right)\right)\]
uniformly for all sufficiently large $x$ and all $m\in\N$ in the same range \eqref{Equ:mh_0h_1}, where the implicit constant in the error term depends at most on the explicit and implicit constants in the hypotheses except for $M$.
\end{thm}
Theorem \ref{thmWeightedHal} is applicable to a large class of nonnegative multiplication functions $\alpha(n)$, including $d_{k}(n)$ for $k\in\N$ and $r_2(n)/4$ \cite{BFI,Ngu}. Despite the great generality of (\ref{Equ:Delta_{alpha}(x;q,a)}), it is oftentimes more convenient to work with the stronger variant
\begin{equation}\label{Equ:Delta_{alpha}(x;q,a)1}
\sum_{\substack{q\le x^{\delta(x)}\\\omega(q)\ll(\delta(x)^2\log\log x)^{1/3}}}\mu(q)^2\rho_g(q)\max_{a\in(\Z/q\Z)^{\times}}|\Delta_{\alpha}(x;q,a)|\ll S(x)\exp\left({-(\log\log x)^{1/3+\epsilon_0}}\right)
\end{equation}
This condition may be viewed as an inequality of the {\it Bombieri--Vinogradov type}, which ensures that the values of $\alpha(n)$ are well distributed as $n$ varies over the reduced residue classes $a\,(\text{mod}\,q)$ for most values of $q$ and $a$. In view of \cite[Proposition 4]{GS}, such a condition arises naturally from a sieving process for the sequence $\{g(n)\}_{n\ge1}$. For this process to work, we need information about the average size of $\alpha(n)$ subject to the constraint $d\mid g(n)$ for smooth square-free $d\in\N\cap[1,x^{\delta}]$. If $d$ is also free of small prime factors up to some constant depending on $g$, then this constraint amounts to the congruences $n\equiv a\,(\text{mod}\, d)$ for $a\in\mathcal{Z}^{\ast}_g(d)$. So, \eqref{Equ:Delta_{alpha}(x;q,a)1} reduces the sieving of $\{g(n)\}_{n\ge1}$ to that of $\N$. 
\par We shall only sketch the proof of Theorem \ref{thmWeightedHal}, since it is similar to, and in fact, much easier than that of Theorem \ref{thmMT2}. The argument used in the proof may also be modified to study the joint distribution of $f(n+h_1)$ and $f(n+h_2)$ with any fixed integers $h_1\ne h_2$.
\par It is not hard to see that the condition $f(p)=O(1)$ in Theorem \ref{thmMT2} can be relaxed, especially when we do not pursue uniformity in $m$ in the asymptotics for the $m$th moment. For instance, in the case $\alpha=1$ Delange and Halberstam showed \cite[Theorem 1]{DH} that if $f\colon\N\to\R$ is a strongly additive function such that $B(x)\to\infty$ as $x\to\infty$, $f(p)=O(\sqrt{B(p)})$ for all primes $p$, and
\begin{equation}\label{S1EquKS}
\sum_{\substack{p\le x\\|f(p)|>\epsilon \sqrt{B(x)}}}\frac{f(p)^2}{p}=o(B(x))
\end{equation}
for any given $\epsilon>0$, then 
\[\frac{1}{x}\sum_{n\leq x}\left(f(n)-A(x)\right)^m=(\mu_m+o(1))B(x)^{\frac{m}{2}}\]
for every fixed $m\in\N$. The implication of this result on the distribution of $f$ is slightly weaker than the Kubilius--Shapiro theorem \cite[Theorem A]{Sha} in that the latter asserts that the distribution of an additive function $f\colon\N\to\R$ with an unbounded variance $B(x)$ which satisfies \eqref{S1EquKS} for every given $\epsilon>0$ is necessarily Gaussian with mean $A(x)$ and variance $B(x)$. On the other hand, Delange and Halberstam noted that their result no longer holds if one removes the assumption $f(p)=O(\sqrt{B(p)})$, which incidentally exposes the limitation of the method of moments compared to the method evolved by Erd\H{o}s and Kac. Regardless, it will be clear in the sequel that the proof of Theorem \ref{thmMT2} makes it possible for us to obtain the following natural extension of the result of Delange and Halberstam.
\begin{thm}\label{thmWeightedDH}
Let $f\colon\N\to\R$ be a strongly additive function, and let $\alpha\in\mathcal{M}^{\ast}$ with parameters $A_0,\beta,\sigma_0,\vartheta_0,\varrho_0,r$. Define
\[B^{\ast}(x)\colonequals\begin{cases}
	~B(x),&\text{~~~~~if $\beta=1$},\\
	~B(x)/(\log\log\log x)^2,&\text{~~~~~if $\beta\ne1$},
\end{cases}\]
and suppose $B^{\ast}(x)\to\infty$ as $x\to\infty$. If there exists a constant $K>0$ such that $f(n)=o(\sqrt{B(x)})$ for all squarefree $n\in\N\cap[1,x]$ composed of prime factors $p$ with $|f(p)|>K\sqrt{B^{\ast}(x)}$, and if
\[\sum_{\substack{p\le x\\|f(p)|>\epsilon \sqrt{B^{\ast}(x)}}}\alpha(p)\frac{f(p)^2}{p^{\sigma_0}}=o(B^{\ast}(x))\]
for any given $\epsilon>0$, then $M(x;m)=(\mu_m+o(1))B(x)^{\frac{m}{2}}$ for every fixed $m\in\N$. 
\end{thm}

The proof of Theorem \ref{thmWeightedDH}, which we shall only outline, is based on the proofs of Theorem \ref{thmMT2} and \cite[Theorem 1]{DH}. We shall also obtain as a corollary the following analogue of the Kubilius--Shapiro theorem \cite[Theorem C]{Sha}.

\begin{cor}\label{corWeightedKS}
Under the notation and hypotheses in Theorem \ref{thmWeightedDH}, we have
\[\lim\limits_{x\to\infty}S(x)^{-1}\sum_{\substack{n\leq x\\f(n)\leq A(x)+V\sqrt{B(x)}}}\alpha(n)=\Phi(V)\]
for any given $V\in\R$. The same is true if $f$ is merely additive.
\end{cor}
It is clear that Theorem \ref{thmWeightedDH} implies Corollary \ref{corWeightedKS} when $f$ is strongly additive. To handle the general case where $f$ is merely additive, we shall establish a weighted version of \cite[Theorem B]{Sha} which shows that when it comes to the distribution problem, there is no essential difference between strongly additive functions and general additive functions, and thus the distribution of an additive function $f$ is determined solely by its values at primes.
\par Corollary \ref{corWeightedKS} has many interesting applications. For instance, it implies at once that if $h\colon\N\to\R$ is any {\it completely additive} function, i.e., $h(mn)=h(m)+h(n)$ for all $m,n\in\N$, and if $k>1$ is a positive integer, then the distribution of $h(d_k(n))$ weighted by $\alpha(n)$ is Gaussian with mean $h(k)\beta \log\log x$ and variance $h(k)^2\beta \log\log x$, provided $h(k)\ne0$. In \cite{Elliott1} Elliott proved a weighted Erd\H{o}s--Kac theorem concerning Ramanujan's $\tau$-function. In Remark \ref{rmk8.1} we describe how his result may be derived from Corollary \ref{corWeightedKS}. Analogues on elliptic holomorphic newforms of weight at least 2 can be obtained in the same way. In a similar fashion, one can also show that if the weight $\alpha$ in Corollary \ref{corWeightedKS} satisfies the additional condition that $\alpha(p)\sim\beta p^{\sigma_0-1}$ for all but a subset $E$ of primes $p$, where $\#(E\cap[2,x])=o(x(\log\log x)^{2-\vartheta_0}/(\log x)^3)$ as $x\to\infty$, then the distribution of $\Omega(\varphi(n))$ weighted by $\alpha(n)$ is Gaussian with mean $\beta(\log\log x)^2/2$ and variance $\beta(\log\log x)^3/3$, generalizing an old result of Erd\H{o}s and Pomerance \cite[Theorem 3.1]{EP} in an easy manner.

\begin{rmk}\label{rmk1.2}
The condition that $f(p)=o((\log p)^{\epsilon})$ for any given $\epsilon>0$, mentioned by Halberstam \cite{Hal2}, does not imply  (\ref{S1EquKS}) in general. To see this, assume for the moment that there exists an infinite subset $\mathcal{P}$ of primes such that
\begin{equation}\label{S1Equ:s(x)}
s_{\mathcal{P}}(x)\colonequals\sum_{p\in\mathcal{P}\cap[17,x]}\frac{1}{p}=\frac{\log\log x}{\log\log\log x}+c+o(1)
\end{equation}
for sufficiently large $x$, where $c\in\R$ is some constant. Define $f(p)=(\log p)^{1/(2\log\log\log p)}$ for $p\in\mathcal{P}$ and $f(p)=1$ for $p\notin\mathcal{P}$. Clearly, $f(p)=o((\log p)^{\epsilon})$ for any given $\epsilon>0$. It is easily seen by partial summation that \[\sum_{p\in\mathcal{P}\cap[17,x]}\frac{f(p)^2}{p}=\int_{17^-}^{x}(\log t)^{1/\log\log\log t}\,ds_{\mathcal{P}}(t)=(1+o(1))(\log x)^{1/\log\log\log x},\]
which implies that
\[B(x)=\sum_{p\in\mathcal{P}\cap[17,x]}\frac{f(p)^2}{p}+O(\log\log x)=(1+o(1))(\log x)^{1/\log\log\log x}.\]
Let $y=x^{1/\log\log x}$ and $\epsilon=1/2$. Since
\begin{align*}
\log\log y&=\log\log x-\log\log\log x,\\
\log\log\log y&=\left(1+O\left(\frac{1}{\log\log x}\right)\right)\log\log\log x,
\end{align*}
we have
\[(\log y)^{1/\log\log\log y}=\exp\left(\frac{\log\log x}{\log\log\log x}-1+O\left(\frac{1}{\log\log\log x}\right)\right).\]
It follows that 
\[f(p)^2> (\log y)^{1/\log\log\log y}>\frac{1}{3}(\log x)^{1/\log\log\log x}>\epsilon^2B(x)\]
for $p\in\mathcal{P}\cap(y,x]$ when $x$ is sufficiently large. Hence, we have
\begin{align*}
\sum_{\substack{p\le x\\|f(p)|>\epsilon \sqrt{B(x)}}}\frac{f(p)^2}{p}\ge\sum_{p\in\mathcal{P}\cap(y,x]}\frac{f(p)^2}{p}>\frac{1}{2}(\log x)^{1/\log\log\log x}>\frac{1}{3}B(x).
\end{align*}
It remains to construct a set $\mathcal{P}$ with the desired property (\ref{S1Equ:s(x)}). The following inductive approach was suggested by Prof. Pomerance. Note first that $\sum_{p\le x}1/p=\log\log x+O(1)$ grows slightly faster than our target $u(x)\colonequals\log\log x/\log\log\log x$, according to Mertens' second theorem \cite[Theorem 427]{HW}. Moreover, if $p<p'$ are large consecutive primes, then $u(p')-u(p)=o(1/\log p)$, by Bertrand's postulate. Let 17 be the first prime in $\mathcal{P}$. Suppose that we have already selected for $\mathcal{P}$ the primes up to $q$, where $q$ is prime. We put the next prime $q'$ in $\mathcal{P}$ if $s_{\mathcal{P}}(q)<u(q)$ and leave it out of $\mathcal{P}$ otherwise.  Then the running sum $s_{\mathcal{P}}(x)$ changes by at most $1/q$ as $x$ moves from $q$ to $q'$, while the target $u(x)$ changes by at most $o(1/\log q)$ as $x$ moves from $q$ to $q'$.  Thus, the difference $s_{\mathcal{P}}(x)-u(x)$ can be kept within $o(1/\log x)$. In particular, (\ref{S1Equ:s(x)}) holds for $\mathcal{P}$ with $c=0$.
\end{rmk}

\noindent {\bf Overview of the proof of Theorem \ref{thmMT2}.} Before embarking on the proofs of our results, we describe briefly the main steps in the proof of Theorem \ref{thmMT2}. The starting point is an approximation to moments used by Granville, Soundararajan, Khan, Milinovich and Subedi. Though the underlying idea is the same, we need a more complicated version of this approximation (see Lemma \ref{lemApprox}) due to the more general nature of our multiplicative weights $\alpha$. To utilize it, we first need to develop an asymptotic formula for the mean value of $\alpha(n)$ with $n\le x$ restricted to squarefree integers $a\in\N\cap[1,x]$ (see Lemma \ref{lemMV2}). An important feature of this formula is that it holds uniformly for all squarefree integer $a\in\N\cap[1,x]$, which is key to both applying the moment approximation and making the moment estimates uniform. This formula will serve as a substitute for the one on $d_k$ developed by Khan, Milinovich and Subedi. Unlike their proof, which is based on Perron's formula and the complete sub-multiplicativity of $d_k$, our proof uses the mean value estimate for $\alpha$ supplied by \cite[Theorem 2.1]{BT} and is completely elementary. This is carried out in the next section.
\par After applying the moment approximation, we find that the estimation of the main contribution can be worked out as in \cite{GS} and \cite{KMS}. It is the estimation of the error terms that is more involved in our case. In particular, the estimation of the error term in the moment approximation provided by Lemma \ref{lemApprox} in Section \ref{S:Moments} requires separate treatments according as $\beta=1$ or $\beta\ne1$. Besides, since the error term in our asymptotic formula for the mean value of $\alpha(n)$ over $a\mid n$ supplied by Lemma \ref{lemMV2} in Section \ref{S:MVT} is weaker than what one can obtain for the special weight $d_k(n)$ by complex analytic approaches, we need to handle the case $\beta\in(0,1)$ with some special care and make a careful selection of parameters accordingly in order to minimize the error terms. With these new technicalities taken care of, we obtain the desired uniform estimates for moments stated in Theorem \ref{thmMT2}.

\medskip
\section{Mean Values of Multiplicative Functions}\label{S:MVT}
Without loss of generality, we may assume $A_0\in(0,1)$ in the sequel. In addition, we shall also make use of the asymptotic formula
\begin{equation}\label{thmMT2Equ4}
\sum_{p\le x}\frac{\alpha(p)}{p^{\sigma_0}}=\beta\log\log x+M_{\alpha}+O\left((\log x)^{-A_0}\right)
\end{equation}
with some constant $M_{\alpha}\in\R$, which follows immediately from (\ref{thmMT2Equ1}) via partial summation. In view of our assumption that $f(p) =O(1)$, this formula implies trivially that $B(x)\ll\log\log x$. Moreover, if we define, for every prime $p$,
\[\psi_0(p)\colonequals\sum_{\nu\ge2}\frac{\alpha(p^{\nu})}{p^{\sigma_0\nu}},\]
then we infer from (\ref{thmMT2Equ0}), (\ref{thmMT2Equ2}) and (\ref{thmMT2Equ3}) that \[\psi_0(p)\ll\frac{(\log\log(p+1))^{\vartheta_0}}{p}\] 
and that $\sum_{p}\psi_0(p)<\infty$. 
\begin{lem}\label{lemSumLog}
Let $\alpha\colon\N\to\R_{\geq0}$ be a multiplicative function satisfying \emph{(\ref{thmMT2Equ0})} and \emph{(\ref{thmMT2Equ3})} with some $\sigma_0,\vartheta_0>0$ and $\varrho_0\in[0,1)$. Fix $h\in\R$, $\epsilon_0\in(0,1)$ and $c_0\in[1,\epsilon_0^{-1})$, and define
\[I_{\alpha,h}(x;a)\colonequals\sum_{\substack{q\le x\\R_q=a}}\frac{\alpha(q)}{q^{\sigma_0}}\left(\log \frac{3x}{q}\right)^{h},\]
where $a\in\N\cap[1,x]$ is squarefree. Then there exists a constant $\delta_0>0$ such that uniformly for all sufficiently large $x$, any $\delta\in[\delta_0\log\log x/\log x,1]$, and any squarefree $a\in\N\cap[1,x]$ with $\omega(a)\le(1-\varrho_0)\epsilon_0\delta^{-1}$, we have
\[I_{\alpha,h}(x;a)=\left(\tilde{\lambda}_{\alpha}(a)+O\left(\frac{2^{O(\omega(a))}}{\log x}\left(\frac{1}{x^{c_0\delta\omega(a)}}+\frac{\epsilon_{h,0}L(a)\log P^+(a)}{a}\right)\right)\right)(\log x)^h,\]
where
\begin{align*}
\tilde{\lambda}_{\alpha}(a)&\colonequals\prod_{p\mid a}\sum_{\nu=1}^{\infty}\frac{\alpha(p^{\nu})}{p^{\sigma_0\nu}},\\
L(a)&\colonequals\prod_{p\mid a}(\log\log(p+1))^{\vartheta_0}.
\end{align*}
\end{lem}
\begin{proof} 
Let $\delta\in(0,1]$ and fix $c_1\in(c_0,\epsilon_0^{-1})$. Put $\delta_1\colonequals (1-\varrho_0)^{-1}c_1\delta$ and $y\colonequals x^{k\delta_1}$. For any squarefree $a=p_1\cdots p_k\in\N\cap[1,x]$ with $p_1<\cdots<p_k\le x$ and $k\le(1-\varrho_0)\epsilon_0\delta^{-1}$, we have $k\delta_1\le c_1\epsilon_0<1$ and 
\[I_{\alpha,h}(x;a)=\sum_{\substack{p_1^{\nu_1}\cdots p_k^{\nu_k}\le x\\\nu_1,...,\nu_k\ge1}}\frac{\alpha(p_1^{\nu_1})\cdots\alpha(p_k^{\nu_k})}{p_1^{\sigma_0\nu_1}\cdots p_k^{\sigma_0\nu_k}}\left(\log \frac{3x}{p_1^{\nu_1}\cdots p_k^{\nu_k}}\right)^{h}.\]
On the one hand, we see that
\begin{align*}
&\hspace*{4.8mm}\sum_{\substack{p_1^{\nu_1}\cdots p_k^{\nu_k}\le y\\\nu_1,...,\nu_k\ge1}}\frac{\alpha(p_1^{\nu_1})\cdots\alpha(p_k^{\nu_k})}{p_1^{\sigma_0\nu_1}\cdots p_k^{\sigma_0\nu_k}}\left(\log \frac{3x}{p_1^{\nu_1}\cdots p_k^{\nu_k}}\right)^{h}\\
&=\sum_{\substack{p_1^{\nu_1}\cdots p_k^{\nu_k}\le y\\\nu_1,...,\nu_k\ge1}}\frac{\alpha(p_1^{\nu_1})\cdots\alpha(p_k^{\nu_k})}{p_1^{\sigma_0\nu_1}\cdots p_k^{\sigma_0\nu_k}}(\log3x)^h\left(1+O\left(\frac{\epsilon_{h,0}}{\log 3x}\sum_{i=1}^k\nu_i\log p_i\right)\right)\\
&=\sum_{\substack{p_1^{\nu_1}\cdots p_k^{\nu_k}\le y\\\nu_1,...,\nu_k\ge1}}\frac{\alpha(p_1^{\nu_1})\cdots\alpha(p_k^{\nu_k})}{p_1^{\sigma_0\nu_1}\cdots p_k^{\sigma_0\nu_k}}(\log x)^h+O\left(\frac{2^{O(k)}\epsilon_{h,0}L(a)\log p_k}{a}(\log x)^{h-1}\right),
\end{align*}
by (\ref{thmMT2Equ3}). From (\ref{thmMT2Equ0}) it follows that
\[\sum_{\substack{p_1^{\nu_1}\cdots p_k^{\nu_k}\le y\\\nu_1,...,\nu_k\ge1}}\frac{\alpha(p_1^{\nu_1})\cdots\alpha(p_k^{\nu_k})}{p_1^{\sigma_0\nu_1}\cdots p_k^{\sigma_0\nu_k}}=\tilde{\lambda}_{\alpha}(a)+O\left(2^{O(k)}\sum_{\substack{p_1^{\nu_1}\cdots p_k^{\nu_k}>y \\\nu_1,...,\nu_k\ge1}}\frac{1}{p_1^{(1-\varrho_0)\nu_1}\cdots p_k^{(1-\varrho_0)\nu_k}}\right).\]
The sum in the error term above may be split into two sums according as $p_2^{\nu_2}\cdots p_{k}^{\nu_{k}}\le y$ or $p_2^{\nu_2}\cdots p_{k}^{\nu_{k}}>y$. In the first sum we must have $p_1^{\nu_1}>y/\left(p_2^{\nu_2}\cdots p_{k}^{\nu_{k}}\right)$. Thus summing over $\nu_1$ and then over $\nu_2,...,\nu_k$, we see that the first sum is 
\[\ll\frac{1}{y^{1-\varrho_0}}\sum_{\substack{p_2^{\nu_2}\cdots p_{k}^{\nu_{k}}\le y \\\nu_2,...,\nu_{k}\ge1}}1\le\frac{2^{O(k)}(\log x)^{k-1}}{x^{c_1k\delta}(\log p_2)\cdots(\log p_k)}.\]
The second sum is simply
\[\sum_{\substack{p_2^{\nu_2}\cdots p_{k}^{\nu_{k}}>y \\\nu_2,...,\nu_{k}\ge1}}\frac{1}{p_2^{(1-\varrho_0)\nu_2}\cdots p_{k}^{(1-\varrho_0)\nu_{k}}}\sum_{\nu_1\ge1}\frac{1}{p_1^{(1-\varrho_0)\nu_1}}\ll\sum_{\substack{p_2^{\nu_2}\cdots p_{k}^{\nu_{k}}>y \\\nu_2,...,\nu_{k}\ge1}}\frac{1}{p_2^{(1-\varrho_0)\nu_2}\cdots p_{k}^{(1-\varrho_0)\nu_{k}}}.\]
It follows that 
\[\sum_{\substack{p_1^{\nu_1}\cdots p_k^{\nu_k}>y \\\nu_1,...,\nu_k\ge1}}\frac{1}{p_1^{(1-\varrho_0)\nu_1}\cdots p_k^{(1-\varrho_0)\nu_k}}\ll\sum_{\substack{p_2^{\nu_2}\cdots p_{k}^{\nu_{k}}>y \\\nu_2,...,\nu_{k}\ge1}}\frac{1}{p_2^{(1-\varrho_0)\nu_2}\cdots p_{k}^{(1-\varrho_0)\nu_{k}}}+\frac{2^{O(k)}(\log x)^{k-1}}{x^{c_1k\delta}(\log p_2)\cdots(\log p_k)}.\]
Repeating this argument, we obtain
\[\sum_{\substack{p_1^{\nu_1}\cdots p_k^{\nu_k}>y \\\nu_1,...,\nu_k\ge1}}\frac{1}{p_1^{(1-\varrho_0)\nu_1}\cdots p_k^{(1-\varrho_0)\nu_k}}\le\frac{2^{O(k)}(\log x)^{k-1}}{x^{c_1k\delta}(\log p_2)\cdots(\log p_k)},\]
from which we deduce
\begin{equation}\label{Equ1:sumalphalog}
\sum_{\substack{p_1^{\nu_1}\cdots p_k^{\nu_k}\le y\\\nu_1,...,\nu_k\ge1}}\frac{\alpha(p_1^{\nu_1})\cdots\alpha(p_k^{\nu_k})}{p_1^{\sigma_0\nu_1}\cdots p_k^{\sigma_0\nu_k}}=\tilde{\lambda}_{\alpha}(a)+O\left(\frac{2^{O(k)}(\log x)^{k-1}}{x^{c_1k\delta}(\log p_2)\cdots(\log p_k)}\right).
\end{equation}
On the other hand, we have
\begin{align*}
\sum_{x_1<p^{\nu}\le x_2}\frac{\alpha(p^{\nu})}{p^{\sigma_0\nu}}\left(\log \frac{3x_2}{p^{\nu}}\right)^{h}&\ll\sum_{\substack{\log_px_1<\nu\le \log_px_2\\\nu\in\Z}}\frac{1}{p^{(1-\varrho_0)\nu}}\left(\log \frac{3x_2}{p^{\nu}}\right)^{h}\\
&=-\int_{\log_px_1}^{\log_px_2}\left(\log \frac{3x_2}{p^t}\right)^{h}\,d\left(\sum_{\substack{t<\nu\le\log_px_2\\\nu\in\Z}}\frac{1}{p^{(1-\varrho_0)t}}\right)
\end{align*}
uniformly for all primes $p$ and all $0<x_1\le x_2$. Using integration by parts, we see that the integral above is equal to 
\[-(\log(3x_2/x_1))^h\sum_{\substack{\log_px_1<\nu\le\log_px_2\\\nu\in\Z}}\frac{1}{p^{(1-\varrho_0)\nu}}-\int_{\log_px_1}^{\log_px_2}\left(\sum_{\substack{t<\nu\le\log_px_2\\\nu\in\Z}}\frac{1}{p^{(1-\varrho_0)t}}\right)\,d\left(\log \frac{3x_2}{p^t}\right)^{h}.\]
Since
\[\sum_{\substack{t<\nu\le\log_px_2\\\nu\in\Z}}\frac{1}{p^{(1-\varrho_0)t}}<\frac{1}{p^{(1-\varrho_0)(\lfloor t\rfloor+1)}}\cdot\frac{p^{1-\varrho_0}}{p^{1-\varrho_0}-1}\ll\frac{1}{p^{(1-\varrho_0)t}},\]
we have
\[(\log(3x_2/x_1))^h\sum_{\substack{\log_px_1<\nu\le\log_px_2\\\nu\in\Z}}\frac{1}{p^{(1-\varrho_0)\nu}}\ll\frac{(\log(3x_2/x_1))^h}{x_1^{1-\varrho_0}}\]
and 
\begin{align*}
\int_{\log_px_1}^{\log_px_2}\left(\sum_{\substack{t<\nu\le\log_px_2\\\nu\in\Z}}\frac{1}{p^{(1-\varrho_0)t}}\right)\,d\left(\log \frac{3x_2}{p^t}\right)^{h}&\ll\epsilon_{h,0}\log p\int_{\log_px_1}^{\log_px_2}\frac{1}{p^{(1-\varrho_0)t}}\left(\log \frac{3x_2}{p^t}\right)^{h-1}\,dt\\
&=\frac{\epsilon_{h,0}}{(3x_2)^{1-\varrho_0}}\int_{\log3}^{\log(3x_2/x_1)}t^{h-1}e^{(1-\varrho_0)t}\,dt\\
&\ll\frac{\epsilon_{h,0}}{(3x_2)^{1-\varrho_0}}(\log(3x_2/x_1))^{h-1}\left(\frac{3x_2}{x_1}\right)^{1-\varrho_0}\\
&=\frac{\epsilon_{h,0}(\log(3x_2/x_1))^{h-1}}{x_1^{1-\varrho_0}}.
\end{align*}
Hence, it follows that
\begin{equation}\label{Equ:sumalphalog}
\sum_{x_1<p^{\nu}\le x_2}\frac{\alpha(p^{\nu})}{p^{\sigma_0\nu}}\left(\log \frac{3x_2}{p^{\nu}}\right)^{h}\ll\frac{(\log(3x_2/x_1))^h}{x_1^{1-\varrho_0}}
\end{equation}
uniformly for all primes $p$ and all $0<x_1\le x_2$. This inequality implies immediately 
\begin{align*}
\sum_{\substack{y<p_1^{\nu_1}\cdots p_k^{\nu_k}\le x\\\nu_1,...,\nu_k\ge1}}\frac{\alpha(p_1^{\nu_1})\cdots\alpha(p_k^{\nu_k})}{p_1^{\sigma_0\nu_1}\cdots p_k^{\sigma_0\nu_k}}\left(\log \frac{3x}{p_1^{\nu_1}\cdots p_k^{\nu_k}}\right)^{h}&\le\frac{2^{O(k)}(\log x)^h}{y^{1-\varrho_0}}\sum_{\substack{p_2^{\nu_2}\cdots p_{k}^{\nu_{k}}\le x \\\nu_2,...,\nu_{k}\ge1}}1\\
&\le\frac{2^{O(k)}(\log x)^{k+h-1}}{x^{c_1k\delta}(\log p_2)\cdots(\log p_k)}.
\end{align*}
Lemma \ref{lemSumLog} now follows upon combining the above with (\ref{Equ1:sumalphalog}) and taking $\delta_0=1/(c_1-c_0)$ with the range $\delta\ge\delta_0 \log \log x/\log x$ in mind.
\end{proof}
Let $\alpha$ be a multiplicative function as in Theorem \ref{thmMT2} with $A_0\in(0,1)$. Suppose first that (\ref{thmMT2Equ2}) holds with the restricted sum $\sum_{p}^{'}$ replaced by the full sum $\sum_p$. For $\sigma_0=1$ De la Bret\`{e}che and Tenenbaum \cite[Theorem 2.1]{BT}
showed
\[\sum_{n\le x}\alpha(n)=\frac{1}{\Gamma(\beta)}\prod_{p}\left(1-\frac{1}{p}\right)^{\beta}\sum_{\nu=0}^{\infty}\frac{\alpha(p^{\nu})}{p^{\nu}}x(\log x)^{\beta-1}\left(1+O\left(\frac{1}{(\log x)^{A_0}}\right)\right),\]
where the implicit constant in the error term depends at most on the explicit and implicit constants in the hypotheses. For the general case where $\sigma_0>0$ is arbitrary, it is easy to show, by applying the above to $\alpha(n)/n^{\sigma_0-1}$ and employing partial summation as in the proof of \cite[Corollary 3.3]{EG}, that
\begin{equation}\label{Equ:S(x)}
S(x)=\sum_{n\le x}\alpha(n)=\lambda_{\alpha}x^{\sigma_0}(\log x)^{\beta-1}\left(1+O\left(\frac{1}{(\log x)^{A_0}}\right)\right),
\end{equation}
where
\begin{equation}\label{Equ:lambda_alpha}
\lambda_{\alpha}\colonequals\frac{1}{\sigma_0\Gamma(\beta)}\prod_{p}\left(1-\frac{1}{p}\right)^{\beta}\sum_{\nu=0}^{\infty}\frac{\alpha(p^{\nu})}{p^{\sigma_0\nu}}.
\end{equation}
Suppose now that (\ref{thmMT2Equ2}) holds with the restricted sum $\sum_{p}^{'}$ being the sum $\sum_{p>Q_0}$, where $Q_0\ge1$ is some constant. Let $P_0\colonequals\prod_{p\le Q_0}p$ and $\textbf{1}_{P_0}(n)$ the indicator function of the set $\{n\in\N\colon \gcd(n,P_0)=1\}$. Then $\alpha(n)\textbf{1}_{P_0}(n)$ is a nonnegative multiplicative function satisfying (\ref{thmMT2Equ0})--(\ref{thmMT2Equ3}) with the sum $\sum_{p}^{'}$ in (\ref{thmMT2Equ2}) replaced by the full sum $\sum_p$. In particular, (\ref{Equ:S(x)}) is applicable to $\alpha(n)\textbf{1}_{P_0}(n)$. Thus, we obtain
\begin{equation}\label{Equ:S(x)P_0}
\sum_{\substack{n\le x\\\gcd(n,P_0)=1}}\alpha(n)=\lambda_{\alpha}(P_0)x^{\sigma_0}(\log 3x)^{\beta-1}\left(1+O\left(\frac{1}{(\log 3x)^{A_0}}\right)\right),
\end{equation}
where 
\[\lambda_{\alpha}(P_0)\colonequals\frac{1}{\sigma_0\Gamma(\beta)}\prod_{p\le Q_0}\left(1-\frac{1}{p}\right)^{\beta}\cdot\prod_{p>Q_0 }\left(1-\frac{1}{p}\right)^{\beta}\sum_{\nu=0}^{\infty}\frac{\alpha(p^{\nu})}{p^{\sigma_0\nu}}.\]
Examining the proof of Lemma \ref{lemSumLog}, we find that for every given $h\in\R$, 
\[
\sum_{\substack{q\le x\\R_q\mid P_0}}\frac{\alpha(q)}{q^{\sigma_0}}\left(\log \frac{3x}{q}\right)^{h}=\prod_{p\le Q_0}\sum_{\nu=0}^{\infty}\frac{\alpha(p^{\nu})}{p^{\sigma_0\nu}}(\log x)^h\left(1+O\left(\frac{1}{\log x}\right)\right)\]
for all sufficiently large $x$. Combining this with (\ref{Equ:S(x)P_0}) gives
\[S(x)=\sum_{\substack{q\le x\\R_q\mid P_0}}\alpha(q)\sum_{\substack{n'\le x/q\\\gcd(n',P_0)=1}}\alpha(n')=\lambda_{\alpha}x^{\sigma_0}(\log x)^{\beta-1}\left(1+O\left(\frac{1}{(\log x)^{A_0}}\right)\right),\]
which is the same as (\ref{Equ:S(x)}). 
\par For our applications, we will need an asymptotic formula for
\[S(x;a)=S_{\alpha}(x;a)\colonequals \sum_{\substack{n\le x\\\gcd(n,a)=1}}\alpha(n)\]
uniform in $a\in\N\cap[1,x]$. One may be tempted to apply (\ref{Equ:S(x)}) to the function $\alpha(n)\textbf{1}_{a}(n)$, where $\textbf{1}_{a}(n)$ is the indicator function of the set $\{n\in\N\colon \gcd(n,a)=1\}$. However, it is not immediately clear whether the implied constant in the error term obtained via this naive approach is independent of $a\in\N\cap[1,x]$. Fortunately, the following lemma provides the desired estimate for $S(x;a)$ under the hypotheses (i)--(iv).
\begin{lem}\label{lemS(x;a)}
For any $\alpha\in\mathcal{M}^{\ast}$ with parameters $A_0,\beta,\sigma_0,\vartheta_0,\varrho_0,r$, we have
\[S(x;a)=x^{\sigma_0}(\log x)^{\beta-1}\left(\lambda_{\alpha}(a)+O\left(\frac{1}{(\log x)^{A_0}}\right)\right)\]
uniformly for all sufficiently large $x$ and all $a\in\N\cap[1,x]$, where 
\[\lambda_{\alpha}(a)\colonequals\frac{1}{\sigma_0\Gamma(\beta)}\prod_{p\mid a}\left(1-\frac{1}{p}\right)^{\beta}\cdot\prod_{p\nmid a}\left(1-\frac{1}{p}\right)^{\beta}\sum_{\nu=0}^{\infty}\frac{\alpha(p^{\nu})}{p^{\sigma_0\nu}},\]
The implicit constant in the error term depends at most on the explicit and implicit constants in the hypotheses.
\end{lem}
\begin{proof}
Let $a\in\N\cap[1,x]$. For simplicity of notation, we write $\sum^{a}$ for sums in which the indices take values coprime to $a$. As we have demonstrated above, there is no loss of generality by assuming that $\sigma_0=1$ and that (\ref{thmMT2Equ2}) holds with the restricted sum $\sum_{p}^{'}$ replaced by the full sum $\sum_p$. Note that
\[0<\lambda_{\alpha}(a)=\lambda_{\alpha}\prod_{p\mid a}\left(\sum_{\nu=0}^{\infty}\frac{\alpha(p^{\nu})}{p^{\nu}}\right)^{-1}\le\lambda_{\alpha}.\]
To estimate $S(x,a)$, we start by connecting it with 
\[T(x;a)=T_{\alpha}(x;a)\colonequals \sideset{}{^a}\sum_{n\le x}\frac{\alpha(n)}{n}.\]
It is clear from (\ref{Equ:S(x)}) that $S(x;a)\le S(x;1)\ll x(\log x)^{\beta-1}$ and $T(x;a)\le T(x,1)\ll(\log x)^{\beta}$. Moreover, it is shown in the the proof of \cite[Theorem 2.1]{BT} that 
\begin{equation}\label{Equ:T(x;1)}
T(x;1)=\left(1+O\left(\frac{1}{\log x}\right)\right)\frac{\lambda_{\alpha}}{\beta}(\log x)^{\beta}.
\end{equation}
Following the proof of \cite[Theorem 2.1]{BT}, we find
\begin{align}\label{Equ:S(x;a)logx}
S(x;a)\log x&=\sideset{}{^a}\sum_{n\le x}\alpha(n)\log n+\sideset{}{^a}\sum_{n\le x}\alpha(n)\log\frac{x}{n}\nonumber\\
&=\sideset{}{^a}\sum_{k\le x}\alpha(k)\sideset{}{^a}\sum_{\substack{p^{\nu}\le x/k\\p\nmid k}}\alpha(p^{\nu})\log p^{\nu}+\int_{1^-}^{x}\frac{S(t;a)}{t}\,dt\nonumber\\
&=\sideset{}{^a}\sum_{k\le x}\alpha(k)\sideset{}{^a}\sum_{p\le x/k}\alpha(p)\log p+O\left(\sum_{k\le x}\alpha(k)\sum_{\substack{p\le x/k\\p\mid k}}\alpha(p)\log p\right)\nonumber\\
&\hspace*{3mm}+O\left(\sum_{k\le x}\alpha(k)\sum_{\substack{p^{\nu}\le x/k\\\nu\ge2}}\alpha(p^{\nu})\log p^{\nu}\right)+O\left(x(\log x)^{\beta-1}\right)\nonumber\\
&=\beta xT(x,a)-\sideset{}{^a}\sum_{k\le x}\alpha(k)\sum_{\substack{p\le x/k\\p\mid a}}\alpha(p)\log p+O\left(x\sum_{k\le x}\frac{\alpha(k)}{k(\log(3x/k))^{A_0}}\right)\nonumber\\
&\hspace*{3mm}+O\left(x(\log x)^{\beta-1}\right).
\end{align}
By partial summation we have
\begin{align}\label{Equ:alpha(k)/k(log(3x/k))^A}
\sum_{k\le x}\frac{\alpha(k)}{k(\log(3x/k))^{A_0}}&=\frac{S(x)}{x(\log 3)^{A_0}}+\int_{1^-}^{x}\frac{\log(3x/t)-A_0}{t^2(\log(3x/t))^{A_0+1}}S(t)\,dt\nonumber\\
&\ll(\log x)^{\beta-1}+\int_{1}^{x}\frac{(\log3t)^{\beta-1}}{t(\log(3x/t))^{A_0}}\,dt\nonumber\\
&=(\log x)^{\beta-1}+\int_{0}^{\log x}\frac{(\log3+t)^{\beta-1}}{(\log3x-t)^{A_0}}\,dt\nonumber\\
&\ll(\log x)^{\beta-1}+\frac{1}{(\log x)^{A_0}}\int_{0}^{(\log x)/2}(\log3+t)^{\beta-1}\,dt\nonumber\\
&\hspace*{3mm}+(\log x)^{\beta-1}\int_{(\log x)/2}^{\log x}\frac{1}{(\log3x-t)^{A_0}}\,dt\nonumber\\
&\ll(\log x)^{\beta-A_0}.
\end{align}
Let $x_1\colonequals x/(\log x)^2$. For $k\le x_1$ we see that 
\[\sum_{\substack{p\le x/k\\p\mid a}}\alpha(p)\log p\ll (\log\log x)^{\vartheta_0}\log a\ll\frac{x}{k(\log(x/k))^{A_0}},\]
so that 
\begin{equation}\label{Equ:k<=x_1}
\sideset{}{^a}\sum_{k\le x_1}\alpha(k)\sum_{\substack{p\le x/k\\p\mid a}}\alpha(p)\log p\ll x\sum_{k\le x}\frac{\alpha(k)}{k(\log(3x/k))^{A_0}}\ll x(\log x)^{\beta-A_0}.
\end{equation}
On the other hand, we have by (\ref{thmMT2Equ1}) that
\begin{align}\label{Equ:k>x_1}
\sideset{}{^a}\sum_{x_1<k\le x}\alpha(k)\sum_{\substack{p\le x/k\\p\mid a}}\alpha(p)\log p\ll x\sum_{x_1<k\le x}\frac{\alpha(k)}{k}&\ll x\left((\log x)^{\beta}-(\log x_1)^{\beta}+O((\log x)^{\beta-1})\right)\nonumber\\
&\ll x(\log x)^{\beta-1}\log\log x,
\end{align}
where we have used (\ref{Equ:T(x;1)}) to estimate the sum over $k$ and the mean value theorem to get 
\[(\log x)^{\beta}-(\log x_1)^{\beta}=\beta\xi^{\beta-1}\log\frac{x}{x_1}\ll(\log x)^{\beta-1}\log\log x\]
for some $\xi\in(\log x_1,\log x)$. Combining (\ref{Equ:k<=x_1}) with (\ref{Equ:k>x_1}), we obtain
\[\sideset{}{^a}\sum_{k\le x}\alpha(k)\sum_{\substack{p\le x/k\\p\mid a}}\alpha(p)\log p\ll x(\log x)^{\beta-A_0}.\]
Inserting this and (\ref{Equ:alpha(k)/k(log(3x/k))^A}) into (\ref{Equ:S(x;a)logx}) yields
\begin{equation}\label{Equ:S(x;a)}
S(x;a)=\frac{\beta x}{\log x}T(x;a)+O\left(x(\log x)^{\beta-1-A_0}\right)
\end{equation}
uniformly for all sufficiently large $x$ and all $a\in\N\cap[1,x]$.
\par It remains to estimate $T(x,a)$. To this end, we repeat the argument above with $\alpha(n)$ replaced by $\alpha(n)/n$. From (\ref{thmMT2Equ1}) it follows that
\begin{equation}\label{thmMT2Equ5}
\sum_{p\le x}\frac{\alpha(p)}{p}\log p=\beta\log x+O\left((\log x)^{1-A_0}\right).
\end{equation}
Thus, we have
\begin{align*}
T(x;a)\log x&=\sideset{}{^a}\sum_{n\le x}\frac{\alpha(n)}{n}\log n+\sideset{}{^a}\sum_{n\le x}\frac{\alpha(n)}{n}\log\frac{x}{n}\\
&=\sideset{}{^a}\sum_{k\le x}\frac{\alpha(k)}{k}\sideset{}{^a}\sum_{\substack{p^{\nu}\le x/k\\p\nmid k}}\frac{\alpha(p^{\nu})}{p^{\nu}}\log p^{\nu}+U(x;a)\\
&=\sideset{}{^a}\sum_{k\le x}\frac{\alpha(k)}{k}\sideset{}{^a}\sum_{p\le x/k}\frac{\alpha(p)}{p}\log p+O\left(\sideset{}{^a}\sum_{k\le x}\frac{\alpha(k)}{k}\sum_{\substack{p\le x/k\\p\mid k}}\frac{\alpha(p)}{p}\log p\right)\\
&\hspace*{3mm}+O\left(\sideset{}{^a}\sum_{k\le x}\frac{\alpha(k)}{k}\sum_{\substack{p^{\nu}\le x/k\\\nu\ge2}}\frac{\alpha(p^{\nu})}{p^{\nu}}\log p^{\nu}\right)+U(x;a)\\
&=(\beta+1)U(x;a)-\sideset{}{^a}\sum_{k\le x}\frac{\alpha(k)}{k}\sum_{\substack{p\le x/k\\p\mid a}}\frac{\alpha(p)}{p}\log p+O\left((\log x)^{1-A_0}T(x;a)\right),
\end{align*}
where 
\begin{equation}\label{Equ:U(x;a)}
U(x;a)\colonequals\sideset{}{^a}\sum_{n\le x}\frac{\alpha(n)}{n}\log\frac{x}{n}=\int_{1^-}^{x}\frac{T(t;a)}{t}\,dt.
\end{equation}
In view of (\ref{thmMT2Equ5}), we have
\begin{align*}
\sum_{\substack{p\le x/k\\p\mid a}}\frac{\alpha(p)}{p}\log p&\le\sum_{p\le (\log x)^2}\frac{\alpha(p)}{p}\log p+\sum_{\substack{(\log x)^2<p\le x\\p\mid a}}\frac{\alpha(p)}{p}\log p\\
&\ll\log\log x+(\log\log x)^{\vartheta_0}\sum_{\substack{(\log x)^2<p\le x\\p\mid a}}\frac{\log p}{p}\\
&\ll\log\log x+(\log\log x)^{\vartheta_0}\omega(a)\frac{\log\log x}{(\log x)^2}\ll\log\log x,
\end{align*}
so that 
\[\sideset{}{^a}\sum_{k\le x}\frac{\alpha(k)}{k}\sum_{\substack{p\le x/k\\p\mid a}}\frac{\alpha(p)}{p}\log p\ll(\log\log x)T(x;a).\]
It follows that
\[T(x;a)\log x=(\beta+1)U(x;a)+O\left((\log x)^{1-A_0}T(x;a)\right).\]
Hence, there exists a function $\epsilon(x;a)$ such that $\epsilon(x;a)=O((\log x)^{-A_0})$ and 
\begin{equation}\label{Equ:T(x;a)}
T(x;a)=\frac{1}{1-\epsilon(x;a)}\cdot\frac{\beta+1}{\log x}U(x;a)
\end{equation}
uniformly for all sufficiently large $x$ and all $a\in\N\cap[1,x]$.
\par Finally, we estimate $U(x;a)$ and $T(x;a)$ by following the proof of \cite[Theorem A]{Song}. For $y\ge2$ and $a\in\N\cap[1,y]$, let
\[V(y;a)\colonequals\log\left(\frac{\beta+1}{(\log y)^{\beta+1}}U(y;a)\right).\]
In light of (\ref{Equ:U(x;a)}) and (\ref{Equ:T(x;a)}), we have
\begin{align*}
\frac{d}{dy}V(y;a)&=-\frac{\beta+1}{y\log y}+\frac{1}{U(y;a)}\cdot\frac{d}{dy}U(y;a)\\
&=-\frac{\beta+1}{y\log y}+\frac{T(y;a)}{U(y;a)y}\\
&=\frac{\beta+1}{y\log y}\cdot\frac{\epsilon(y;a)}{1-\epsilon(y;a)}\ll\frac{1}{y(\log y)^{A_0+1}}
\end{align*}
uniformly for all sufficiently large $y$ and all $a\in\N\cap[1,y]$, which implies that
\[V_a\colonequals\int_{2}^{\infty}\frac{d}{dy}V(y;a)\,dy<\infty.\]
Since
\[V(x;a)-V(2;a)=V_a-\int_{x}^{\infty}\frac{d}{dy}V(y;a)\,dy=V_a+O\left((\log x)^{-A_0}\right)\]
uniformly for all sufficiently large $x$ and all $a\in\N\cap[1,x]$, it follows that
\[\frac{\beta+1}{(\log x)^{\beta+1}}U(x;a)=\exp(V(x;a))=\exp(V_a+V(2;a))\left(1+O\left((\log x)^{-A_0}\right)\right).\]
Combining this estimate with (\ref{Equ:T(x;a)}), we infer 
\begin{equation}\label{Equ:Asymp T(x;a)}
T(x;a)=\exp(V_a+V(2;a))(\log x)^{\beta}\left(1+O\left((\log x)^{-A_0}\right)\right)
\end{equation}
uniformly for all sufficiently large $x$ and all $a\in\N\cap[1,x]$. The leading coefficient can be made explicit by arguing as in the proof of \cite[Theorem A]{Song}. Alternatively, we can also take advantage of (\ref{Equ:S(x)P_0}). Fixing $a\in\N$, we have by (\ref{Equ:S(x)P_0}) with $\sigma_0=1$ that
\[T(x;a)=\frac{\lambda_{\alpha}(a)}{\beta}(\log x)^{\beta}\left(1+O\left((\log x)^{-A_0}\right)\right)\]
for all sufficiently large $x$. Comparing this with (\ref{Equ:Asymp T(x;a)}) shows that $\exp(V_a+V(2;a))=\lambda_{\alpha}(a)/\beta$. Carrying this back into (\ref{Equ:Asymp T(x;a)}), we obtain
\[T(x;a)=\frac{\lambda_{\alpha}(a)}{\beta}(\log x)^{\beta}\left(1+O\left((\log x)^{-A_0}\right)\right)\]
uniformly for all sufficiently large $x$ and all $a\in\N\cap[1,x]$. Inserting the above into (\ref{Equ:S(x;a)}) completes the proof of the lemma.
\end{proof}

\par The next result, which is key to the computation of moments, is a direct corollary of Lemmas \ref{lemSumLog} and \ref{lemS(x;a)}.
\begin{lem}\label{lemMV2}
Fix $\epsilon_0\in(0,1)$, and let $\alpha\in\mathcal{M}^{\ast}$ with parameters $A_0,\beta,\sigma_0,\vartheta_0,\varrho_0,r$. Then there exist constants $\delta_0>0$ and $Q_0\ge2$, such that uniformly for all sufficiently large $x$, any $\delta\in[\delta_0 \log \log x/\log x,1]$, and any square-free $a\in \N\cap[1,x]$ with $\omega(a)\le(1-\rho_0)\epsilon_0\delta^{-1}$, $P^-(a)>Q_0$ and $P^+(a)\le x^{\delta}$, we have
\[\sum_{\substack{n\le x\\a\mid n}}\alpha(n)=\lambda_{\alpha}\left(F(\sigma_0,a)+O\left(\frac{2^{O(\omega(a))}L(a)}{a}\left(\frac{1}{(\log x)^{A_0}}+\frac{\epsilon_{\beta,1}\log P^+(a)}{\log x}\right)\right)\right)x^{\sigma_0}(\log x)^{\beta-1},\]
where $L(a)$ is defined as in Lemma \ref{lemSumLog},
\[F(\sigma_0,a)\colonequals\prod_{p\mid a}\left(1-\left(\sum_{\nu=0}^{\infty}\alpha(p^{\nu})p^{-\sigma_0\nu }\right)^{-1}\right),\]
and $\lambda_{\alpha}$ is defined by \emph{(\ref{Equ:lambda_alpha})}. 
\end{lem}
\begin{proof}
Suppose that $\delta_0>0$ is a constant for which Lemma \ref{lemSumLog} holds when $c_0=1$ and $h\in\{\beta-1,\beta-1-A_0\}$. Let $Q_0\ge2$ be such that 
\[\sum_{\nu=1}^{\infty}\frac{\alpha(p^{\nu})}{p^{\sigma_0\nu}}\le\frac{1}{2}\]
for all $p>Q_0$. Then we have
\begin{equation}\label{Equ:F(sigma_0,p)}
F(\sigma_0,p)=\sum_{\nu=1}^{\infty}\frac{\alpha(p^{\nu})}{p^{\sigma_0\nu}}+O\left(\left(\sum_{\nu=1}^{\infty}\frac{\alpha(p^{\nu})}{p^{\sigma_0\nu}}\right)^2\right)=\frac{\alpha(p)}{p^{\sigma_0}}+O\left(\psi_0(p)+\frac{\alpha(p)^2}{p^{2\sigma_0}}\right)
\end{equation}
for all $p>Q_0$. For any square-free integer $a\in[1,x]$ with $\omega(a)\le (1-\varrho_0)\epsilon_0\delta^{-1}$, $P^-(a)>Q_0$ and $P^+(a)\le x^{\delta}$, we have by Lemma \ref{lemS(x;a)} that
\begin{equation}\label{lemMV2Equ1}
\sum_{\substack{n\le x\\\gcd(n,a)=1}}\alpha(n)=x^{\sigma_0}(\log 3x)^{\beta-1}\left(\lambda_{\alpha}(a)+O\left(\frac{1}{(\log 3x)^{A_0}}\right)\right).
\end{equation}
Note that
\[\sum_{\substack{n\le x\\a\mid n}}\alpha(n)=\sum_{\substack{q\le x\\R_q=a}}\alpha(q)\sum_{\substack{n'\le x/q\\\gcd(n',a)=1}}\alpha(n').\]
By (\ref{lemMV2Equ1}), the main term of the inner sum contributes
\[\lambda_{\alpha}(a)x^{\sigma_0}\sum_{\substack{q\le x\\R_q=a}}\alpha(q)\left(\log \frac{3x}{q}\right)^{\beta-1},\]
which, by Lemma \ref{lemSumLog}, is equal to
\begin{align*}
&\hspace*{4.8mm}\lambda_{\alpha}(a)x^{\sigma_0}\left(\tilde{\lambda}_{\alpha}(a)+O\left(\frac{2^{O(\omega(a))}}{\log x}\left(\frac{1}{x^{\delta\omega(a)}}+\frac{\epsilon_{\beta,1}L(a)\log P^+(a)}{a}\right)\right)\right)(\log x)^{\beta-1}\\
&=\lambda_{\alpha}\left(F(\sigma_0,a)+O\left(\frac{2^{O(\omega(a))}}{\log x}\left(\frac{1}{a}+\frac{\epsilon_{\beta,1}L(a)\log P^+(a)}{a}\right)\right)\right)x^{\sigma_0}(\log x)^{\beta-1},
\end{align*}
since $a\le x^{\delta\omega(a)}$. Analogously, the contribution from the error term of the inner sum is
\begin{align*}
&\ll\lambda_{\alpha}\left(F(\sigma_0,a)+\frac{2^{O(\omega(a))}L(a)\log P^+(a)}{a\log x}\right)x^{\sigma_0}(\log x)^{\beta-1-A_0}\\
&\ll\frac{\lambda_{\alpha}2^{O(\omega(a))}L(a) }{a}x^{\sigma_0}(\log x)^{\beta-1-A_0},
\end{align*}
where we have used the estimate $F(\sigma_0,a)\ll2^{O(\omega(a))}L(a)/a$, which follows directly from (\ref{thmMT2Equ3}) and (\ref{Equ:F(sigma_0,p)}). Combining these estimates completes the proof of Lemma \ref{lemMV2}.
\end{proof}

\begin{rmk}\label{rmk2.1}
We point out that the lower bound $Q_0$ for $\omega(a)$ in the lemma above is by and large an artificial thing, whose value is insignificant for our applications. However, we need it because (\ref{Equ:F(sigma_0,p)}) may not hold for small primes. As we shall see later, having such a lower bound also frees us from dealing with minor contributions from small primes.
\end{rmk}

\medskip
\section{Computing Moments}\label{S:Moments}
By rescaling the strongly additive function $f$ in Theorem \ref{thmMT2}, we may assume, without loss of generality, that $|f(p)|\leq1$ for all primes $p$. Note that $0\le F(\sigma_0,p)<1$ for all primes $p$. For every $p$ we define $f_p\colon\N\to\R$ by
\[f_p(n)\colonequals\begin{cases}
	~f(p)(1-F(\sigma_0,p)),&\text{~~~~~if $p\mid n$},\\
	~-f(p)F(\sigma_0,p),&\text{~~~~~otherwise}.
\end{cases}\]
Given any $q\in\N$ we may also extend $f_p$ via complete multiplicativity by setting
\[f_q(n)\colonequals \prod_{p^{\nu}\parallel q}f_{p}(n)^{\nu}.\]
It is clear that $|f_q(n)|\le 1$. The following result provides an approximation of the moments of $f$ in terms of those of $f_p$.
\begin{lem}\label{lemApprox}
Let $\alpha\in\mathcal{M}^{\ast}$ with parameters $A_0,\beta,\sigma_0,\vartheta_0,\varrho_0,r$, and let $f\colon\N\to\R$ be a strongly additive function with $|f(p)|\le 1$ for all $p$. Then there exists a constant $Q_0\ge2$, such that
\[\sum_{n\leq y}\alpha(n)(f(n)-A(x))^m=\sum_{n\leq y}\alpha(n)\left(\sum_{Q_0<p\leq z}f_p(n)\right)^m+O\left(E(y,z,w;m)\right)\]
holds uniformly for all sufficiently large $x\ge z$, any $y\ge1$, and all $m\in\N$, where 
\[E(y,z,w;m)\colonequals\sum_{\substack{a+b+c=m\\0\le a<m\\b,c\ge0}}\binom{m}{a,b,c}2^{O(m-a)}\left(\log(v+2)\right)^{c}\sum_{n\leq y}\alpha(n)\left|\sum_{Q_0<p\leq z}f_p(n)\right|^a\omega(n;z,w)^{b},\]
$v\colonequals \log x/\log z$, $w\colonequals x^{1/\log(v+2)}$, and
\[\omega(n;z,w)\colonequals \sum_{\substack{z<p\le w\\p\mid n}}1.\]
\end{lem}
\begin{proof}
Let $Q_0\ge2$ be a constant for which (\ref{Equ:F(sigma_0,p)}) holds. Suppose that $z>Q_0$ is sufficiently large. By (\ref{thmMT2Equ3}), (\ref{Equ:F(sigma_0,p)}) and the fact that $\sum_{p}\psi_0(p)<\infty$, we find 
\[\sum_{Q_0<p\leq x}f(p)F(\sigma_0,p)=A(x)+O(1).\]
We compute
\begin{align*}
f(n)-A(x)&=\sum_{\substack{p\mid n\\p>Q_0}}f(p)-\sum_{Q_0<p\leq x}f(p)F(\sigma_0,p)+O(1)\\
&=\sum_{\substack{Q_0<p\leq z\\p\mid n}}f(p)+\sum_{\substack{p>z\\p\mid n}}f(p)-\sum_{Q_0<p\leq z}f(p)F(\sigma_0,p)-\sum_{z<p\leq x}f(p)F(\sigma_0,p)+O(1)\\
&=\sum_{Q_0<p\leq z}f_p(n)+\sum_{\substack{p>z\\p\mid n}}f(p)-\sum_{z<p\leq x}f(p)F(\sigma_0,p)+O(1).
\end{align*}
By (\ref{thmMT2Equ4}) we have
\[\left|\sum_{z<p\leq x}f(p)F(\sigma_0,p)\right|\le\sum_{z<p\leq x}\frac{\alpha(p)}{p^{\sigma_0}}+O(1)\le\beta\log(v+2)+O(1).\]
Since
\[\sum_{\substack{p>z\\p\mid n}}|f(p)|\le\sum_{\substack{z<p\le x\\p\mid n}}1<\omega(n;z,w)+\log(v+2),\]
it follows that
\[f(n)-A(x)=\sum_{Q_0<p\leq z}f_p(n)+O\left(\omega(n;z,w)+\log(v+2)\right).\]
We have therefore proved
\[\sum_{n\leq y}\alpha(n)(f(n)-A(x))^m=\sum_{n\leq y}\alpha(n)\left(\sum_{Q_0<p\leq z}f_p(n)+O\left(\omega(n;z,w)+\log(v+2)\right)\right)^m.\]
Opening the $m$th power on the right-hand side by means of the multinomial theorem completes the proof of Lemma \ref{lemApprox}.
\end{proof}

Let $z=x^{1/v}$ and $w=x^{1/\log(v+2)}$ be as in Lemma \ref{lemApprox}, where $v\ge1$ is a function of $x$ and $m$ to be chosen later. Fix $\epsilon_0\in(0,1)$ and $\eta_0\in(0,1]$, and suppose that $y\in[x^{\eta_0},x]$. Under the hypotheses in Theorem \ref{thmMT2}, we seek to estimate the weighted moments
\[\sum_{n\le y}\alpha(n)\left(\sum_{Q_0<p\le z}f_p(n)\right)^m,\]
appearing in Lemma \ref{lemApprox}. Expanding out the $m$th power we see that
\begin{equation}\label{Equ:Main1}
\sum_{n\leq y}\alpha(n)\left(\sum_{Q_0<p\le z}f_p(n)\right)^m=\sum_{Q_0<p_1,...,p_m\le z}\sum_{n\leq y}\alpha(n)f_{p_1\cdots p_m}(n).
\end{equation}
This suggests that we study the sum 
\[\sum_{n\le y}\alpha(n)f_{q}(n)\]
for $q\in\N$ with $\omega(q)\le m$, $P^-(q)>Q_0$ and $P^+(q)\le z$. A key observation is that $f_q(n)=f_q(\gcd(n,R_q))$. From this we deduce 
\[\sum_{n\leq y}\alpha(n)f_{q}(n)=\sum_{a\mid R_q}f_q(a)\sum_{\substack{n\leq y\\\gcd(n,R_q)=a}}\alpha(n)=\sum_{ab\mid R_q}f_q(a)\mu(b)\sum_{\substack{n\leq y\\ab\mid n}}\alpha(n).\]
Note that $\log y/\log z\in[\eta_0v,v]$. By Lemma \ref{lemMV2}, there exists a constant $v_0>0$, independent of $Q_0$ and $\eta_0$, such that 
\begin{equation}\label{Equ:Sumf_q}
\sum_{n\leq y}\alpha(n)f_{q}(n)=\lambda_{\alpha}\left(G(\sigma_0,q)+O\left(2^{O(m)}E_y(q)\right)\right)y^{\sigma_0}(\log y)^{\beta-1}
\end{equation}
holds uniformly for all sufficiently large $x$, any $y\in[x^{\eta_0},x]$ and $v\in[\eta_0^{-1},v_0\log x/\log\log x]$, and all $m\le(1-\varrho_0)\epsilon_0\log y/\log z$, where 
\begin{align*}
G(\sigma_0,q)&\colonequals\sum_{ab\mid R_q}f_q(a)\mu(b)F(\sigma_0,ab),\\
E_y(q)&\colonequals\sum_{ab\mid R_q}\frac{|f_q(a)|L(ab)}{ab}\left(\frac{1}{(\log y)^{A_0}}+\frac{\epsilon_{\beta,1}\log P^+(ab)}{\log y}\right).
\end{align*}
Combining (\ref{Equ:Sumf_q}) with (\ref{Equ:Main1}) gives
\begin{equation}\label{Equ:Main2}
\sum_{n\leq y}\alpha(n)\left(\sum_{Q_0<p\leq z}f_p(n)\right)^m=\lambda_{\alpha}\left(G(z)+O\left(2^{O(m)}D(y,z)\right)\right)y^{\sigma_0}(\log y)^{\beta-1},
\end{equation}
where 
\begin{align*}
G(z)&\colonequals\sum_{Q_0<p_1,...,p_m\leq z}G(\sigma_0,p_1\cdots p_m),\\
D(y,z)&\colonequals\sum_{Q_0<p_1,...,p_m\leq z}E_y(p_1\cdots p_m).
\end{align*}

\medskip
\section{Estimation of $G(z)$ and $D(y,z)$}\label{S:G(z)&D(y,z)}
It is easy to see that $G(\sigma_0,q)$ is multiplicative as a function of $q$. Indeed, given any $q_1,q_2\in\N$ with $\gcd(q_1,q_2)=1$, we have
\begin{align*}
G(\sigma_0,q_1)G(\sigma_0,q_2)&=\sum_{\substack{a_1b_1\mid R_{q_1}\\ a_2b_2\mid R_{q_2}}}f_{q_1}(a_1)f_{q_2}(a_2)\mu(b_1)\mu(b_2)F(\sigma_0,a_1b_1)F(\sigma_0,a_2b_2)\\
&=\sum_{\substack{a_1b_1\mid R_{q_1}\\ a_2b_2\mid R_{q_2}}}f_{q_1}(a_1a_2)f_{q_2}(a_1a_2)\mu(b_1b_2)F(\sigma_0,a_1a_2b_1b_2)\\
&=\sum_{\substack{a_1b_1\mid R_{q_1}\\ a_2b_2\mid R_{q_2}}}f_{q_1q_2}(a_1a_2)\mu(b_1b_2)F(\sigma_0,a_1a_2b_1b_2)\\
&=\sum_{ab\mid R_{q_1q_2}}f_{q_1q_2}(a)\mu(b)F(\sigma_0,ab)=G(\sigma_0,q_1q_2).
\end{align*}
Furthermore, we have
\begin{align*}
G(\sigma_0,p^\nu)&=f_{p^\nu}(1)+f_{p^{\nu}}(p)F(\sigma_0,p)-f_{p^{\nu}}(1)F(\sigma_0,p)\\
&=(-f(p)F(\sigma_0,p))^{\nu}+(f(p)(1-F(\sigma_0,p)))^{\nu}F(\sigma_0,p)-(-f(p)F(\sigma_0,p))^{\nu}F(\sigma_0,p)\\
&=f(p)^{\nu}F(\sigma_0,p)(1-F(\sigma_0,p))\left((-1)^{\nu}F(\sigma_0,p)^{\nu-1}+(1-F(\sigma_0,p))^{\nu-1}\right)
\end{align*}
for all prime powers $p^{\nu}$. Note that $G(\sigma_0,p)=0$, $|G(\sigma_0,p^\nu)|\leq1/4$, and $G(\sigma_0,p^\nu)\geq0$ when $2\mid\nu$. In addition, we have by (\ref{Equ:F(sigma_0,p)}) that
\begin{equation}\label{Equ:G(sigma_0,p^2)}
G(\sigma_0,p^2)=f(p)^2F(\sigma_0,p)(1-F(\sigma_0,p))=\alpha(p)\frac{f(p)^2}{p^{\sigma_0}}+O\left(\psi_0(p)+\frac{\alpha(p)^2}{p^{2\sigma_0}}\right)
\end{equation}
and that
\begin{equation}\label{Equ:G(sigma_0,p^nu)}
|G(\sigma_0,p^\nu)|\le |f(p)|^{\nu}F(\sigma_0,p)\le\alpha(p)\frac{f(p)^2}{p^{\sigma_0}}+O\left(\psi_0(p)+\frac{\alpha(p)^2}{p^{2\sigma_0}}\right)
\end{equation}
for all $p^{\nu}$ with $p>Q_0$ and $\nu\ge2$.
\par Now we proceed to estimate $G(z)$ in the main term of (\ref{Equ:Main2}). Recall that $y\in [x^{\eta_0},x]$ and $z=x^{1/v}$. We shall suppose in this section that $1\le m\le \min(v,h_0B(x)^{1/3})$, $\log(v+2)=o(B(x))$, and $m\log(v+2)\ll B(x)$, where $0<h_0<(3/2)^{2/3}$ is any given constant, and obtain a uniform treatment for $G(z)$ and $D(y,z)$ under this more general assumption. Since $G(\sigma_0,q)$ is multiplicative in $q$ and $G(\sigma_0,p)=0$ for all $p>Q_0$, we have
\begin{equation}\label{Equ:G(z)}
G(z)=\sum_{\substack{Q_0<p_1,...,p_m\leq z\\p_1\cdots p_m \text{~squareful}}}G(\sigma_0,p_1\cdots p_m).
\end{equation}
When $2\mid m$, the main contribution arises from
\begin{equation}\label{Equ:MainMoment}
\frac{m!}{(m/2)!2^{m/2}}\sum_{\substack{Q_0<p_1,...,p_{m/2}\leq z\\p_1,...,p_{m/2} \text{~distinct}}}G(\sigma_0,p_1^2\cdots p_{m/2}^2)=C_m\sum_{\substack{Q_0<p_1,...,p_{m/2}\leq z\\p_1,...,p_{m/2} \text{~distinct}}}\prod_{i=1}^{m/2}G(\sigma_0,p_i^2),
\end{equation}
since the number of ways to partition a set of $m$ elements into $m/2$ two-element equivalence classes is 
\[\frac{m!}{(m/2)!2^{m/2}}=\frac{m!}{m!!}=C_m.\]
The sum on the right-hand side of (\ref{Equ:MainMoment}) can be rewritten as
\[\sum_{\substack{Q_0<p_1,...,p_{m/2-1}\leq z\\p_1,...,p_{m/2-1} \text{~distinct}}}\prod_{i=1}^{m/2-1}G(\sigma_0,p_i^2)\sum_{\substack{Q_0<p_{m/2}\leq z\\p_{m/2}\ne p_1,...,p_{m/2-1}}}G(\sigma_0,p_{m/2}^2).\]
By (\ref{Equ:G(sigma_0,p^2)}) and (\ref{thmMT2Equ4}), the inner sum over $p_{m/2}$ is equal to 
\[\sum_{Q_0<p\leq z}G(\sigma_0,p^2)-\sum_{i=1}^{m/2-1}G(\sigma_0,p_i^2)=B(z)+O\left(\sum_{Q_0<p\le q_{N}}\frac{\alpha(p)}{p^{\sigma_0}}\right)=B(z)+O\left(\log\log(m+2)\right),\]
where $N=m/2+\pi(Q_0)$ and $q_{N}$ is the $N$th prime. Repeating this argument we obtain
\[\sum_{\substack{Q_0<p_1,...,p_{m/2}\leq z\\p_1,...,p_{m/2} \text{~distinct}}}\prod_{i=1}^{m/2}G(\sigma_0,p_i^2)=\left(B(z)+O\left(\log\log(m+2)\right)\right)^{m/2}.\]
But
\[B(x)-B(z)=\sum_{z<p\leq x}\alpha(p)\frac{f(p)^2}{p^{\sigma_0}}\le\beta\log(v+2)+O(1).\]
Hence when $m$ is even, the main contribution to $G(z)$ is given by
\[C_m\left(B(x)+O(\log(v+2))\right)^{m/2}=C_mB(x)^{\frac{m}{2}}\left(1+O\left(mB(x)^{-1}\log(v+2)\right)\right).\]
The remaining contribution to $G(z)$ comes from
\begin{equation}\label{Equ:MinorMoment}
\sum_{s<m/2}\sum_{Q_0<p_1<\cdots<p_s\leq z}\sum_{\substack{k_1+\cdots+k_s=m\\k_1,...,k_s\geq2}}\binom{m}{k_1,...,k_s}\prod_{i=1}^{s}G(\sigma_0,p_i^{k_i}).
\end{equation}
Since (\ref{Equ:MinorMoment}) vanishes when $m\le 2$, we may suppose $m\ge3$. By (\ref{Equ:G(sigma_0,p^nu)}) we see that
\begin{align*}
\prod_{i=1}^{s}\left|G(\sigma_0,p_i^{k_i})\right|&\leq\prod_{i=1}^{s}\left(\alpha(p_i)\frac{f(p_i)^2}{p_i^{\sigma_0}}+O\left(\psi_0(p_i)+\frac{\alpha(p_i)^2}{p_i^{2\sigma_0}}\right)\right).
\end{align*}
Thus, we have
\[\sum_{Q_0<p_1<\cdots<p_s\leq z}\prod_{i=1}^{s}\left|G(\sigma_0,p_i^{k_i})\right|\le\frac{1}{s!}(B(x)+O(1))^s=\frac{1}{s!}B(x)^s\left(1+O\left(sB(x)^{-1}\right)\right)\ll\frac{B(x)^s}{s!}.\]
Since
\[\sum_{\substack{k_1+\cdots+k_s=m\\k_1,...,k_s\geq2}}\binom{m}{k_1,...,k_s}\le\frac{m!}{2^s}\sum_{\substack{k_1+\cdots+k_s=m\\k_1,...,k_s\geq2}}1=\frac{m!}{2^s}\binom{m-s-1}{s-1},\footnote{We have corrected the binomial coefficient in \cite[Equ. (11)]{GS}.}\]
(\ref{Equ:MinorMoment}) is
\[\ll m!\sum_{s<m/2}\frac{1}{s!2^s }\binom{m-s-1}{s-1}B(x)^s.\]
To estimate the sum above, we put $m_1\colonequals \lfloor(m-1)/2\rfloor$ and observe that
\begin{align*}
\sum_{s<m/2}\frac{1}{s!2^s }\binom{m-s-1}{s-1}B(x)^s&=B(x)^{m_1}\sum_{s\le m_1}\frac{1}{s!2^s }\binom{m-s-1}{s-1}B(x)^{s-m_1}\\
&\le B(x)^{m_1}m^{-3m_1}\sum_{s\le m_1}\frac{1}{s!2^s }\binom{m-s-1}{s-1}h_0^{3(m_1-s)}m^{3s},
\end{align*}
where we have used the assumption that $B(x)\ge m^3/h_0^3$ with some $0<h_0<(3/2)^{2/3}$. Let 
\[e_m\colonequals\begin{cases}
	~1,&\text{~~~~~if $2\mid m$},\\
	~1/2,&\text{~~~~~otherwise}.
\end{cases}\]
Then $m_1=m/2-e_m$. Note that 
\begin{align*}
m^{-3m_1}\sum_{s\le m/4}\frac{1}{s!2^s }\binom{m-s-1}{s-1}h_0^{3(m_1-s)}m^{3s}&\le m^{-3m_1}\sum_{s\le m/4}\frac{1}{s!(s-1)!}\left(\frac{9}{4}\right)^{m_1-s}\left(\frac{m^4}{2}\right)^s\\
&\ll m^{-3m_1}\left(\frac{9}{4}\right)^{m_1}\left(\frac{m^4}{2}\right)^{m/4}\ll\frac{C_m}{m!}m^{3e_m},
\end{align*}
since 
\[\frac{C_m}{m!}=\frac{1}{2^{m/2}\Gamma(m/2+1)}\asymp m^{-\frac{m+1}{2}}e^{\frac{m}{2}}\]
by Stirling's formula. Next, we have
\begin{align*}
m^{-3m_1}\sum_{m/4<s\le m/3}\frac{1}{s!2^s}\binom{m-s-1}{s-1}h_0^{3(m_1-s)}m^{3s}&\le 2^{O(m)}m^{-3m_1}\sum_{m/4<s\le m/3}\frac{1}{s!(s-1)!}\left(\frac{m^4}{2}\right)^s\\
&\le \frac{2^{O(m)}m^{-3m_1}}{m^{m/2}}\sum_{m/4<s\le m/3}\left(\frac{m^4}{2}\right)^s\\
&\le\frac{2^{O(m)}m^{-3m_1}}{m^{m/2}}m^{4m/3}\\
&=2^{O(m)}m^{-2m/3+3e_m}\ll\frac{C_m}{m!}m^{3e_m}.
\end{align*}
Finally, we observe that
\begin{align*}
&\hspace*{5.3mm} m^{-3m_1}\sum_{m/3<s\le m_1}\frac{1}{s!2^s}\binom{m-s-1}{s-1}h_0^{3(m_1-s)}m^{3s}\\
&=m^{-3m_1}\sum_{m/3<s\le m_1}\frac{1}{s!2^s}\binom{m-s-1}{m-2s}h_0^{3(m_1-s)}m^{3s}\\
&\le m^{-3m_1}\sum_{m/3<s\le m_1}\frac{1}{s!2^s}(m-s)^{m-2s}h_0^{3(m_1-s)}m^{3s}\\
&\le \frac{m^{-3m_1}}{m_1!}\sum_{m/3<s\le m_1}\frac{m_1!}{s!2^s}\left(\frac{2m}{3}\right)^{m-2s}h_0^{3(m_1-s)}m^{3s}\\
&\le \frac{m^{-3m_1}}{m_1!}\sum_{m/3<s\le m_1}\frac{1}{2^s}\left(\frac{m}{2}\right)^{m_1-s}\left(\frac{2m}{3}\right)^{m-2s}h_0^{3(m_1-s)}m^{3s}\\
&\ll \frac{m^{m-2m_1}}{m_1!2^{m/2}}\sum_{m/3<s\le m_1}\left(\frac{2h_0^{3/2}}{3}\right)^{m-2s}\ll\frac{C_m}{m!}m^{3e_m}.
\end{align*}
Collecting the estimates above, we see that the contribution to $G(z)$ from (\ref{Equ:MinorMoment}) is 
\[\ll C_m m^{3e_m}B(x)^{m_1}=C_m B(x)^{\frac{m}{2}}\left(\frac{m^{3}}{B(x)}\right)^{e_m}\le C_mB(x)^{\frac{m}{2}} \frac{m^{\frac{3}{2}}}{\sqrt{B(x)}}.\]
We can therefore conclude that
\begin{equation}\label{Equ:Asymp for G(z)}
G(z)=C_mB(x)^{\frac{m}{2}}\left(\chi_m\left(1+O\left(\frac{m\log(v+2)}{B(x)}\right)\right)+O\left(\frac{m^{\frac{3}{2}}}{\sqrt{B(x)}}\right)\right).
\end{equation}
\par Next, we estimate $D(y,z)$ in the error term of (\ref{Equ:Main2}). By definition, we have
\[D(y,z)=\sum_{s\leq m}\sum_{Q_0<p_1<\cdots<p_s\leq z}\sum_{\substack{k_1+\cdots+k_s=m\\k_1,...,k_s\in\N}}\binom{m}{k_1,...,k_s}E_y\left(p_1^{k_1}\cdots p_s^{k_s}\right).\]
Let 
\[H(\sigma_0,q)\colonequals\sum_{ab\mid R_q}\frac{|f_q(a)|L(ab)}{ab}.\]
Then $H(\sigma_0,q)$ is multiplicative in $q$. Moreover, we have
\[E_y(q)\le H(\sigma_0,q)\left(\frac{1}{(\log y)^{A_0}}+\frac{\epsilon_{\beta,1}\log P^+(q)}{\log y}\right).\]
It follows that $D(y,z)\le D_1(y,z)+\epsilon_{\beta,1}D_2(y,z)$, where
\begin{align*}
D_1(y,z)&\colonequals\frac{1}{(\log y)^{A_0}}\sum_{s\leq m}\sum_{Q_0<p_1<\cdots<p_s\leq z}\sum_{\substack{k_1+\cdots+k_s=m\\k_1,...,k_s\in\N}}\binom{m}{k_1,...,k_s}\prod_{i=1}^{s}H(\sigma_0,p_i^{k_i}),\\
D_2(y,z)&\colonequals\frac{1}{\log y}\sum_{s\leq m}\sum_{Q_0<p_1<\cdots<p_s\leq z}\log p_s\sum_{\substack{k_1+\cdots+k_s=m\\k_1,...,k_s\in\N}}\binom{m}{k_1,...,k_s}\prod_{i=1}^{s}H(\sigma_0,p_i^{k_i}).
\end{align*}
By Mertens' theorems \cite[Theorems 425, 427]{HW} we have, for any $t\ge3$, that
\begin{equation}\label{Equ:sum(loglog(p+1))^{theta_0}/p}
\sum_{p\le t}\frac{(\log\log(p+1))^{\vartheta_0}}{p}=\frac{1}{\vartheta_0+1}(\log\log t)^{\vartheta_0+1}+O(1)
\end{equation}
and that
\begin{equation}\label{Equ:sum(loglog(p+1))^{theta_0}logp/p}
\sum_{p\le t}\frac{(\log\log(p+1))^{\vartheta_0}\log p}{p}=\left(1+O\left(\frac{1}{\log t}+\frac{\vartheta_0}{\log\log t}\right)\right)(\log\log t)^{\vartheta_0}\log t.
\end{equation}
Furthermore, let 
\[T_{n}(t)\colonequals\sum_{k=0}^{n}{n\brace k}t^k\] 
denote the $n$th Touchard polynomial, where 
\[{n\brace k}\colonequals\frac{1}{k!}\sum_{\substack{n_1+\cdots+n_k=n\\n_1,...,n_k\in\N}}\binom{n}{n_1,...,n_k}\]
is the $k$th Stirling number of the second kind of size $n$. The sequence $\{T_n(t)\}_{n=0}^{\infty}$ of the Touchard polynomials is known to satisfy the recurrence relation
\[T_{n+1}(t)=t\sum_{i=0}^{n}\binom{n}{i}T_{i}(t),\]
from which one verifies readily by induction that 
\begin{equation}\label{Equ:Touchard}
T_n(t)\le\left(t+\frac{n-1}{2}\right)^n
\end{equation}
for all $n\ge1$ and $t\ge0$. Since
\begin{align*}
H(\sigma_0,p^{\nu})&=|f(p)|F(\sigma_0,p)\left(1+\frac{L(p)}{p}\right)+\frac{|f(p)|L(p)}{p}(1-F(\sigma_0,p))\\
&=|f(p)|\left(F(\sigma_0,p)+\frac{(\log\log(p+1))^{\vartheta_0}}{p}\right)
\end{align*}
for any prime powers $p^{\nu}$ with $p>Q_0$, we obtain, from (\ref{Equ:F(sigma_0,p)}), (\ref{Equ:sum(loglog(p+1))^{theta_0}/p}), (\ref{Equ:sum(loglog(p+1))^{theta_0}logp/p}) and (\ref{Equ:Touchard}), that
\begin{align*}
D_1(y,z)&\le\frac{2^{O(m)}}{(\log x)^{A_0}}\sum_{s\leq m}\frac{1}{s!}(\log\log z)^{s(\vartheta_0+1)}\sum_{\substack{k_1+\cdots+k_s=m\\k_1,...,k_s\in\N}}\binom{m}{k_1,...,k_s}\\
&\le\frac{2^{O(m)}}{(\log x)^{A_0}}T_m\left((\log\log z)^{\vartheta_0+1}\right)\le\frac{2^{O(m)}}{(\log x)^{A_0}}(\log\log x)^{m(\vartheta_0+1)},
\end{align*}
and that
\begin{align*}
D_2(y,z)&\le\frac{2^{O(m)}\log z}{\log x}\sum_{s\leq m}\frac{1}{(s-1)!}(\log\log z)^{s(\vartheta_0+1)-1}\sum_{\substack{k_1+\cdots+k_s=m\\k_1,...,k_s\in\N}}\binom{m}{k_1,...,k_s}\\
&\le\frac{2^{O(m)}}{v\log\log z}T_m\left((\log\log z)^{\vartheta_0+1}\right)\le\frac{2^{O(m)}}{v}(\log\log x)^{m(\vartheta_0+1)-1}.
\end{align*}
Hence, we conclude that
\begin{equation}\label{Equ:D(y,z)}
D(y,z)\le2^{O(m)}(\log\log x)^{m(\vartheta_0+1)-1}\left(\frac{\log\log x}{(\log x)^{A_0}}+\frac{\epsilon_{\beta,1}}{v}\right).
\end{equation}

\medskip
\section{Estimation of $E(y,z,w;m)$}\label{S:E(y,z,w;m)}
In this section, we seek to bound the function $E(y,z,w;m)$ introduced in Lemma \ref{lemApprox} under the assumptions in Theorem \ref{thmMT2}. We start with the case $\beta=1$. Suppose that $1\le m\le h_0B(x)^{1/3}$, where $0<h_0<(3/2)^{2/3}$ is any given constant. Recall that $y\in [x^{\eta_0},x]$, $z=x^{1/v}$ and $w=x^{1/\log(v+2)}$. With the choice $v=(1-\varrho_0)^{-1}\epsilon_0^{-1}\eta_0^{-1}m$, we clearly have $v\in[\eta_0^{-1},v_0\log x/\log\log x]$ and $m\le(1-\varrho_0)\epsilon_0\log y/\log z$. Inputting (\ref{Equ:Asymp for G(z)}) and (\ref{Equ:D(y,z)}) into (\ref{Equ:Main2}), we obtain
\begin{equation}\label{Equ:Main3,beta=1}
\sum_{n\leq y}\alpha(n)\left(\sum_{Q_0<p\leq z}f_p(n)\right)^m=\lambda_{\alpha}C_mB(x)^{\frac{m}{2}}\left(\chi_m+O\left(\frac{m^{\frac{3}{2}}}{\sqrt{B(x)}}\right)\right)y^{\sigma_0}.
\end{equation}
The key lies in the estimation of the sum
\begin{equation}\label{Equ:omega_z(n)^{b}}
\sum_{n\leq y}\alpha(n)\left|\sum_{Q_0<p\leq z}f_p(n)\right|^a\omega(n;z,w)^{b}.
\end{equation}
In the present case, we may simply use the trivial bound $\omega(n;z,w)\ll v\ll m$, so that (\ref{Equ:omega_z(n)^{b}}) is bounded above by
\[2^{O(b)}m^b\sum_{n\leq y}\alpha(n)\left|\sum_{Q_0<p\leq z}f_p(n)\right|^a.\]
It is clear that we can use (\ref{Equ:Main3,beta=1}) to handle the sum above. If $a$ is even, then this sum is $\ll \lambda_{\alpha}C_aB(x)^{\frac{a}{2}}y^{\sigma_0}$; if $a$ is odd, then it is 
\begin{align*}
&\le\left(\sum_{n\le y}\alpha(n)\left|\sum_{Q_0<p\leq z}f_p(n)\right|^{a-1}\right)^{1/2}\left(\sum_{n\le y}\alpha(n)\left|\sum_{Q_0<p\leq z}f_p(n)\right|^{a+1}\right)^{1/2}\\
&\ll\lambda_{\alpha}\sqrt{C_{a-1}C_{a+1}}B(x)^{\frac{a}{2}}y^{\sigma_0}
\end{align*}
by the Cauchy--Schwarz inequality. The sequence $\{C_{\ell}\}_{\ell=1}^{\infty}$ is strictly increasing, which can be easily seen from the identity
\[\frac{C_{\ell+1}}{C_\ell}=\frac{\ell+1}{\sqrt{2}}\cdot\frac{\Gamma(\ell/2+1)}{\Gamma((\ell+1)/2+1)}=\sqrt{2}\cdot\frac{\Gamma(\ell/2+1)}{\Gamma((\ell/2+1/2)}\]
and the fact that $\Gamma(y)$ is strictly increasing on $[3/2,\infty)$. Moreover, we have by Stirling's formula that
\[\frac{C_{\ell}}{C_{\ell+1}}\ll\frac{1}{\ell+1}\cdot\frac{((\ell+1)/2)^{\ell/2+1}e^{-(\ell+1)/2}}{(\ell/2)^{(\ell+1)/2}e^{-\ell/2}}\ll\frac{1}{\sqrt{\ell+1}},\]
which implies that 
\[C_{a}\le2^{O(m-a)}C_m\sqrt{\frac{a!}{m!}}\le2^{O(m-a)}C_m\sqrt{\frac{a^a}{m^m}}\le\frac{2^{O(m-a)}C_m}{\left(\sqrt{m}\right)^{m-a}}\]
for all $0\le a\le m$. Hence, (\ref{Equ:omega_z(n)^{b}}) is bounded above by
\[\frac{2^{O(m-a)}\lambda_{\alpha} C_m m^b}{\left(\sqrt{m}\right)^{m-a}}B(x)^{\frac{a}{2}}y^{\sigma_0}\le2^{O(m-a)}\lambda_{\alpha} C_m\left(\sqrt{m}\right)^{m-a}B(x)^{\frac{a}{2}}y^{\sigma_0}.\]
Inputting this inequality into the definition of $E(y,z,w;m)$, we conclude that 
\begin{equation}\label{Equ:E(y,z,w;m),beta=1}
E(y,z,w;m)\le \lambda_{\alpha}C_my^{\sigma_0}\sum_{a=0}^{m-1}\binom{m}{a}B(x)^{\frac{a}{2}}\left(O\left(\sqrt{m}\right)\right)^{m-a}\ll \lambda_{\alpha}C_m m^{\frac{3}{2}}B(x)^{\frac{m-1}{2}}y^{\sigma_0}.
\end{equation}
\par Now we consider the case $\beta\ne1$. Suppose that $1\le m\ll B(x)^{1/3}/(\log\log\log x)^{2/3}$ and that $B(x)/(\log\log\log x)^2\to \infty$ as $x\to\infty$. In this case we take $v=(\log\log x)^{m(\vartheta_0+2)}$, so that $v\in[2\eta_0^{-1},v_0\log x/\log\log x]$ and $m\le(1-\varrho_0)\epsilon_0\log t/\log z$ for any $t\in[x^{\eta_0/2},x]$ when $x$ is sufficiently large. Inserting (\ref{Equ:Asymp for G(z)}) and (\ref{Equ:D(y,z)}) into (\ref{Equ:Main2}) leads to 
\begin{equation}\label{Equ:Main3}
\sum_{n\leq t}\alpha(n)\left(\sum_{Q_0<p\leq z}f_p(n)\right)^m=\lambda_{\alpha}C_mB(x)^{\frac{m}{2}}\left(\chi_m+O\left(\frac{m^{\frac{3}{2}}}{\sqrt{B(x)}}\right)\right)t^{\sigma_0}(\log t)^{\beta-1}
\end{equation}
uniformly for all $t\in[x^{\eta_0/2},x]$. Again, we need to estimate (\ref{Equ:omega_z(n)^{b}}) uniformly for $y\in[x^{\eta_0},x]$. Note that (\ref{Equ:omega_z(n)^{b}}) can be rewritten as
\[\sum_{k=1}^{b}\sum_{z<p_1<\cdots<p_k\le w}\sum_{\substack{l_1+\cdots+l_k=b\\l_1,...,l_k\ge1}}\binom{b}{l_1,...,l_k}\sum_{\substack{n\le y\\ p_1\cdots p_{k}\mid n}}\alpha(n)\left|\sum_{Q_0<p\leq z}f_p(n)\right|^a.\]
Observe that
\begin{align*}
\sum_{\substack{n\le y\\ p_1\cdots p_{k}\mid n}}\alpha(n)\left|\sum_{Q_0<p\leq z}f_p(n)\right|^a&=\sum_{\substack{q\le y\\R_q=p_1\cdots p_k}}\alpha(q)\sum_{\substack{n'\le y/q\\\gcd(n',q)=1}}\alpha(n')\left|\sum_{Q_0<p\leq z}f_p(n')\right|^a\\
&\le\sum_{\substack{q\le y\\R_q=p_1\cdots p_k}}\alpha(q)\sum_{n\le y/q}\alpha(n)\left|\sum_{Q_0<p\leq z}f_p(n)\right|^a,
\end{align*}
since $p_1,...,p_k>p$.  If $q=p_1^{\nu_1}\cdots p_k^{\nu_k}>\sqrt{y}$ with given $z<p_1<\cdots<p_k\le w$, then we have the trivial estimate
\[\sum_{n\le y/q}\alpha(n)\left|\sum_{Q_0<p\leq z}f_p(n)\right|^a\le\pi(z)^a\sum_{n\le 3y/q}\alpha(n)\ll\lambda_{\alpha}\pi(z)^a\left(\frac{y}{q}\right)^{\sigma_0}\left(\log\frac{3y}{q}\right)^{\beta-1}\]
by (\ref{Equ:S(x)}) and the fact that $|f_p(n)|\le1$. By the proof of Lemma \ref{lemSumLog}, and particularly by (\ref{Equ:sumalphalog}), we find that
\begin{align*}
\sum_{\substack{\sqrt{y}<q\le y\\R_q=p_1\cdots p_k}}\frac{\alpha(q)}{q^{\sigma_0}}\left(\log\frac{3y}{q}\right)^{\beta-1}&=\sum_{\substack{\sqrt{y}<p_1^{\nu_1}\cdots p_k^{\nu_k}\le y\\\nu_1,...,\nu_k\ge1}}\frac{\alpha(p_1^{\nu_1})\cdots\alpha(p_k^{\nu_k})}{p_1^{\sigma_0\nu_1}\cdots p_k^{\sigma_0\nu_k}}\left(\log\frac{3y}{p_1^{\nu_1}\cdots p_k^{\nu_k}}\right)^{\beta-1}\\
&\ll\frac{2^{O(k)}(\log y)^{k+\beta-2}}{\left(\sqrt{y}\right)^{1-\varrho_0}},
\end{align*}
from which it follows that
\[\sum_{\substack{\sqrt{y}<q\le y\\R_q=p_1\cdots p_k}}\alpha(q)\sum_{n\le y/q}\alpha(n)\left|\sum_{Q_0<p\leq z}f_p(n)\right|^a\ll\lambda_{\alpha}\pi(z)^a\frac{2^{O(k)}y^{\sigma_0}(\log y)^{k+\beta-2}}{\left(\sqrt{y}\right)^{1-\varrho_0}}.\]
Summing the above over all $z<p_1<\cdots<p_k\le w$ yields immediately
\begin{align}\label{Equ:sum for sqrt{y}<q<=y}
&\hspace*{4.8mm}\sum_{z<p_1<\cdots<p_k\le w}\sum_{\substack{\sqrt{y}<q\le y\\R_q=p_1\cdots p_k}}\alpha(q)\sum_{n\le y/q}\alpha(n)\left|\sum_{Q_0<p\leq z}f_p(n)\right|^a\nonumber\\
&\le\lambda_{\alpha}\pi(z)^a\pi(w)^k\frac{2^{O(k)}y^{\sigma_0}(\log y)^{k+\beta-2}}{k!\left(\sqrt{y}\right)^{1-\varrho_0}}\le\frac{\lambda_{\alpha}y^{\sigma_0}(\log y)^{\beta-1}}{k!\left(\sqrt[3]{y}\right)^{1-\varrho_0}}
\end{align}
for sufficiently large $x$, since $y\in[x^{\eta_0},x]$, $a+k\le m\ll(\log\log x)^{1/3}/(\log\log\log x)^{2/3}$, and
\[\pi(z)^a\pi(w)^k\le \left(\frac{w}{\log w}+O\left(\frac{w}{(\log w)^2}\right)\right)^m\ll\left(\frac{w}{\log w}\right)^m\le \frac{x^{1/\log\log\log x}(m\log\log\log x)^m}{(\log x)^m}.\]
If $q=p_1^{\nu_1}\cdots p_k^{\nu_k}\le\sqrt{y}$, then $x^{\eta_0/2}\le\sqrt{y}\le y/q\le y\le x$. Thus, we can apply (\ref{Equ:Main3}) with $t=y/q$ to handle 
\[\sum_{n\le y/q}\alpha(n)\left|\sum_{Q_0<p\leq z}f_p(n)\right|^a.\]
If $a$ is even, then this sum is 
\[\ll\lambda_{\alpha}C_aB(x)^{\frac{a}{2}}\left(\frac{y}{q}\right)^{\sigma_0}\left(\log\frac{y}{q}\right)^{\beta-1}\le\frac{2^{O(m-a)}\lambda_{\alpha}C_m}{\left(\sqrt{m}\right)^{m-a}}B(x)^{\frac{a}{2}}\left(\frac{y}{q}\right)^{\sigma_0}\left(\log\frac{y}{q}\right)^{\beta-1};\]
if $a$ is odd, then it is 
\begin{align*}
&\le\left(\sum_{n\le y/q}\alpha(n)\left|\sum_{Q_0<p\leq z}f_p(n)\right|^{a-1}\right)^{1/2}\left(\sum_{n\le y/q}\alpha(n)\left|\sum_{Q_0<p\leq z}f_p(n)\right|^{a+1}\right)^{1/2}\\
&\ll\lambda_{\alpha}\sqrt{C_{a-1}C_{a+1}}B(x)^{\frac{a}{2}}\left(\frac{y}{q}\right)^{\sigma_0}\left(\log\frac{y}{q}\right)^{\beta-1}\\
&\le\frac{2^{O(m-a)}\lambda_{\alpha}C_m}{\left(\sqrt{m}\right)^{m-a}}B(x)^{\frac{a}{2}}\left(\frac{y}{q}\right)^{\sigma_0}\left(\log\frac{y}{q}\right)^{\beta-1}
\end{align*}
by Cauchy--Schwarz. It follows that
\begin{align*}
&\hspace*{4.8mm}\sum_{\substack{q\le \sqrt{y}\\R_q=p_1\cdots p_k}}\alpha(q)\sum_{n\le y/q}\alpha(n)\left|\sum_{Q_0<p\leq z}f_p(n)\right|^a\\
&\le\frac{2^{O(m-a)}\lambda_{\alpha}C_my^{\sigma_0}}{\left(\sqrt{m}\right)^{m-a}}B(x)^{\frac{a}{2}}
\sum_{\substack{q\le \sqrt{y}\\R_q=p_1\cdots p_k}}\frac{\alpha(q)}{q^{\sigma_0}}\left(\log\frac{y}{q}\right)^{\beta-1}\\
&\le\frac{2^{O(m-a)}\lambda_{\alpha}C_m}{\left(\sqrt{m}\right)^{m-a}}B(x)^{\frac{a}{2}}y^{\sigma_0}(\log y)^{\beta-1}\prod_{i=1}^{k}\sum_{\nu=1}^{\infty}\frac{\alpha(p_i^{\nu})}{p_i^{\sigma_0\nu}}\\
&=\frac{2^{O(m-a)}\lambda_{\alpha}C_m}{\left(\sqrt{m}\right)^{m-a}}B(x)^{\frac{a}{2}}y^{\sigma_0}(\log y)^{\beta-1}\prod_{i=1}^{k}\left(\frac{\alpha(p_i)}{p_i^{\sigma_0}}+\psi_0(p_i)\right)
\end{align*}
for all $0\le a<m$. Since (\ref{thmMT2Equ4}) implies that
\[\sum_{z<p_1<\cdots<p_k\le w}\prod_{i=1}^{k}\left(\frac{\alpha(p_i)}{p_i^{\sigma_0}}+\psi_0(p_i)\right)\le\frac{1}{k!}\left(\sum_{z<p\le w}\left(\frac{\alpha(p)}{p^{\sigma_0}}+\psi_0(p)\right)\right)^k\le\frac{2^{O(k)}}{k!}(\log v)^k,\]
we obtain
\begin{align}\label{Equ:sum for y<=sqrt{y}}
&\hspace*{4.8mm}\sum_{z<p_1<\cdots<p_k\le w}\sum_{\substack{q\le \sqrt{y}\\R_q=p_1\cdots p_k}}\alpha(q)\sum_{n\le y/q}\alpha(n)\left|\sum_{Q_0<p\leq z}f_p(n)\right|^a\nonumber\\
&\le\frac{2^{O(m-a)}\lambda_{\alpha}C_m}{k!\left(\sqrt{m}\right)^{m-a}}(\log v)^kB(x)^{\frac{a}{2}}y^{\sigma_0}(\log y)^{\beta-1}.
\end{align}
Combining (\ref{Equ:sum for y<=sqrt{y}}) with (\ref{Equ:sum for sqrt{y}<q<=y}) and extending the inner sum over $q$ to the entire range, we conclude that
\begin{align*}
&\hspace*{4.8mm}\sum_{z<p_1<\cdots<p_k\le w}\sum_{\substack{q\le y\\R_q=p_1\cdots p_k}}\alpha(q)\sum_{n\le y/q}\alpha(n)\left|\sum_{Q_0<p\leq z}f_p(n)\right|^a\\
&\le\frac{2^{O(m-a)}\lambda_{\alpha}C_m}{k!\left(\sqrt{m}\right)^{m-a}}(\log v)^kB(x)^{\frac{a}{2}}y^{\sigma_0}(\log y)^{\beta-1}.
\end{align*}
Hence, (\ref{Equ:omega_z(n)^{b}}) is bounded above by
\begin{align*}
&\hspace*{4.8mm} \frac{2^{O(m-a)}\lambda_{\alpha}C_m}{\left(\sqrt{m}\right)^{m-a}}B(x)^{\frac{a}{2}}y^{\sigma_0}(\log y)^{\beta-1}\sum_{k=1}^{b}\frac{(\log v)^k}{k!}\sum_{\substack{l_1+\cdots+l_k=b\\l_1,...,l_k\ge1}}\binom{b}{l_1,...,l_k}\\
&=\frac{2^{O(m-a)}\lambda_{\alpha}C_m}{\left(\sqrt{m}\right)^{m-a}}B(x)^{\frac{a}{2}}y^{\sigma_0}(\log y)^{\beta-1}\sum_{k=1}^{b}{b\brace k}(\log v)^k\\
&\le\frac{2^{O(m-a)}\lambda_{\alpha}C_m}{\left(\sqrt{m}\right)^{m-a}}B(x)^{\frac{a}{2}}T_{b}(\log v)y^{\sigma_0}(\log y)^{\beta-1}.
\end{align*}
It follows by (\ref{Equ:Touchard}) that the above does not exceed
\[\frac{2^{O(m-a)}\lambda_{\alpha}C_m}{\left(\sqrt{m}\right)^{m-a}}B(x)^{\frac{a}{2}}(\log v)^b y^{\sigma_0}(\log y)^{\beta-1},\]
where we have used the observation that $\log v > 
m\log\log\log x>m\ge b$. In other words, we have shown that
\[\sum_{n\leq y}\alpha(n)\left|\sum_{Q_0<p\leq z}f_p(n)\right|^a\omega(n;z,w)^{b}\le\frac{2^{O(m-a)}\lambda_{\alpha}C_m}{\left(\sqrt{m}\right)^{m-a}}B(x)^{\frac{a}{2}}(\log v)^by^{\sigma_0}(\log y)^{\beta-1}.\]
Inputting this inequality into the definition of $E(y,z,w;m)$, we conclude that 
\begin{align}\label{Equ:E(y,z,w;m)}
E(y,z,w;m)&\le \lambda_{\alpha}C_my^{\sigma_0}(\log y)^{\beta-1}\sum_{a=0}^{m-1}\binom{m}{a}B(x)^{\frac{a}{2}}\left(O\left(\frac{\log v}{\sqrt{m}}\right)\right)^{m-a}\nonumber\\
&\ll \lambda_{\alpha}C_m\sqrt{m}(\log v)B(x)^{\frac{m-1}{2}}y^{\sigma_0}(\log y)^{\beta-1}\nonumber\\
&\ll \lambda_{\alpha}C_m m^{\frac{3}{2}}(\log\log\log x)B(x)^{\frac{m-1}{2}}y^{\sigma_0}(\log y)^{\beta-1}.
\end{align}

\medskip
\section{Deduction of Theorems \ref{thmMT2} and \ref{thmMT2V}}\label{S:Proof of Thms1.1&1.2}
Theorem \ref{thmMT2} now follows immediately upon combining (\ref{Equ:Main3,beta=1}) and (\ref{Equ:Main3}) with (\ref{Equ:E(y,z,w;m),beta=1}) and (\ref{Equ:E(y,z,w;m)}) and invoking Lemma \ref{lemApprox} and (\ref{Equ:S(x)}). In fact, we have shown that the same asymptotic formulas which hold for $M(x;m)$ also hold for 
\begin{equation}\label{Equ:GMomentM(y;m)}
S(y)^{-1}\sum_{n\leq y}\alpha(n)(f(n)-A(x))^m
\end{equation}
uniformly in the range $y\in[x^{\eta_0},x]$, where $\eta_0\in(0,1]$ is any fixed constant.
\par Now we prove Theorem \ref{thmMT2V}. Recall that under the hypotheses in Theorem \ref{thmMT2V}, the multiplicative function $\alpha(n)$ satisfies conditions (i)--(iv). We shall again suppose $A_0\in(0,1)$ throughout the proof. Define the strongly additive function $\tilde{f}\colon\N\to\R$, called the {\it strongly additive contraction of $f$}, by $\tilde{f}(p)=f(p)$ for all primes $p$. Then
\begin{equation}\label{Equ:from tilde{f} to f}
\sum_{n\leq x}\alpha(n)(f(n)-A(x))^m=\sum_{k=0}^{m}\binom{m}{k}\sum_{n\leq x}\alpha(n)\left(\tilde{f}(n)-A(x)\right)^k\left(f(n)-\tilde{f}(n)\right)^{m-k}
\end{equation}
for every $m\in\N$. The term corresponding to $k=m$ can be estimated directly using Theorem \ref{thmMT2}. Hence, it remains to deal with
\begin{equation}\label{Equ:sum(tilde{f}-A(x))^k(f(n)-tilde{f}(n))^l}
\sum_{n\leq x}\alpha(n)\left(\tilde{f}(n)-A(x)\right)^k\left(f(n)-\tilde{f}(n)\right)^l
\end{equation}
for $0\leq k<m$ and $l=m-k$. Note that
\begin{align*}
&\hspace*{4.8mm}\left|\sum_{n\leq x}\alpha(n)\left(\tilde{f}(n)-A(x)\right)^k\left(f(n)-\tilde{f}(n)\right)^{l}\right|\\
&\le\sum_{n\leq x}\alpha(n)\left|\tilde{f}(n)-A(x)\right|^k\left|\sum_{p^{\nu}\parallel n,\nu\geq2}(f(p^{\nu})-f(p))\right|^l\\
&\le\sum_{p_1,...,p_l\le \sqrt{x}}\sum_{\substack{p_1^{\nu_1},..., p_l^{\nu_l}\le x\\\nu_1,...,\nu_l\geq2}}|f(p_1^{\nu_1})-f(p_1)|\cdots |f(p_l^{\nu_l})-f(p_l)|\sum_{\substack{n\le x\\p_1^{\nu_1},...,p_l^{\nu_l}\parallel n}}\alpha(n)\left|\tilde{f}(n)-A(x)\right|^k.
\end{align*}
Since $f(p^{\nu})=O(\nu^{\kappa})$ for all $p^{\nu}$, the last expression above does not exceed
\[2^{O(l)}\sum_{s\leq l}\sum_{p_1<\cdots<p_s\le \sqrt{x}}\sum_{\substack{l_1+\cdots+l_s=l\\l_1,...,l_s\in\N}}\binom{l}{l_1,...,l_s}\sum_{\substack{p_1^{\nu_1}\cdots p_s^{\nu_s}\le x\\\nu_1,...,\nu_s\ge2}}\nu_1^{\kappa l_1}\cdots \nu_s^{\kappa l_s}\sum_{\substack{n\le x\\p_1^{\nu_1},...,p_s^{\nu_s}\parallel n}}\alpha(n)\left|\tilde{f}(n)-A(x)\right|^k.\] 
If we write $n=p_1^{\nu_1}\cdots p_s^{\nu_s}n'$ with $\gcd(n',p_1\cdots p_s)=1$, then it is clear that 
\[\left|\tilde{f}(n)-A(x)\right|^k=\left|\tilde{f}(n')-A(x)+\sum_{i=1}^sf(p_i)\right|^k\le\sum_{a=0}^k\binom{k}{a}\left|\tilde{f}(n')-A(x)\right|^{a}\left|\sum_{i=1}^sf(p_i)\right|^{k-a}.\]
Thus, the innermost sum of $\alpha(n)|\tilde{f}(n)-A(x)|^k$ is
\begin{equation}\label{Equ:sumalpha(n)|tilde{f}(n)-A(x)|^k}
\le\alpha(p_1^{\nu_1})\cdots\alpha(p_s^{\nu_s})\sum_{a=0}^k\binom{k}{a}\left|\sum_{i=1}^sf(p_i)\right|^{k-a}\sum_{n\le x/\left(p_1^{\nu_1}\cdots p_s^{\nu_s}\right)}\alpha(n)\left|\tilde{f}(n)-A(x)\right|^a,
\end{equation}
where we have dropped the superscript of $n$ for simplicity of notation. Since the right-hand side of the above clearly vanishes if $p_1\cdots p_s>\sqrt{x}$, we may assume $p_1\cdots p_s\le\sqrt{x}$ instead. Let $\lambda'\colonequals 1-\varrho_0-\log_2\lambda>\rho_0$, and choose a constant $\max(1/2,\sqrt{\varrho_0/\lambda'})<\delta_0<1$, so that $1-\varrho_0+\delta_0^2\lambda'>1$. Let $x_s\colonequals x/(p_1\cdots p_s)$ and $y_s\colonequals x_s^{\delta_0}$. Then $x_s\ge\sqrt{x}\ge p_1\cdots p_s$. If $p_1^{\nu_1}\cdots p_s^{\nu_s}>p_1\cdots p_sy_s$ with given $p_1<\cdots<p_s$, then we use the trivial estimate
 \begin{align*}
\sum_{n\le x/\left(p_1^{\nu_1}\cdots p_s^{\nu_s}\right)}\alpha(n)\left|\tilde{f}(n)-A(x)\right|^a&\ll2^{O(a)}(\log x)^{a}\sum_{n\le x/\left(p_1^{\nu_1}\cdots p_s^{\nu_s}\right)}\alpha(n)\\
&\ll\lambda_{\alpha}2^{O(a)}(\log x)^{a}\left(\frac{x}{p_1^{\nu_1}\cdots p_s^{\nu_s}}\right)^{\sigma_0}\left(\log\frac{3x}{p_1^{\nu_1}\cdots p_s^{\nu_s}}\right)^{\beta-1}.
 \end{align*}
Thus, (\ref{Equ:sumalpha(n)|tilde{f}(n)-A(x)|^k}) is 
\[\ll\frac{\alpha(p_1^{\nu_1})\cdots\alpha(p_s^{\nu_s})}{p_1^{\sigma_0\nu_1}\cdots p_s^{\sigma_0\nu_s}}\left(\log\frac{3x}{p_1^{\nu_1}\cdots p_s^{\nu_s}}\right)^{\beta-1}\lambda_{\alpha}2^{O(k)}x^{\sigma_0}(\log x)^k.\]
Since $\alpha(p^{\nu})=O((\lambda p^{\varrho_0+\sigma_0-1})^{\nu})$ for all $p^{\nu}$, we have
\begin{align*}
&\hspace*{4.8mm}\sum_{\substack{p_1\cdots p_sy_s<p_1^{\nu_1}\cdots p_s^{\nu_s}\le x\\\nu_1,...,\nu_s\ge2}}\nu_1^{\kappa l_1}\cdots \nu_s^{\kappa l_s}\frac{\alpha(p_1^{\nu_1})\cdots\alpha(p_s^{\nu_s})}{p_1^{\sigma_0\nu_1}\cdots p_s^{\sigma_0\nu_s}}\left(\log\frac{3x}{p_1^{\nu_1}\cdots p_s^{\nu_s}}\right)^{\beta-1}\\
&\le2^{O(l)}\sum_{\substack{p_1\cdots p_sy_s<p_1^{\nu_1}\cdots p_s^{\nu_s}\le x\\\nu_1,...,\nu_s\ge2}}\nu_1^{\kappa l_1}\cdots \nu_s^{\kappa l_s}\left(\frac{\lambda}{p_1^{1-\varrho_0}}\right)^{\nu_1}\cdots\left(\frac{\lambda}{p_s^{1-\varrho_0}}\right)^{\nu_s}\left(\log\frac{3x}{p_1^{\nu_1}\cdots p_s^{\nu_s}}\right)^{\beta-1}\\
&\le\frac{2^{O(l)}}{(p_1\cdots p_s)^{1-\varrho_0}}\sum_{\substack{y_s<p_1^{\nu_1}\cdots p_s^{\nu_s}\le x_s\\\nu_1,...,\nu_s\ge1}}\nu_1^{\kappa l_1}\cdots \nu_s^{\kappa l_s}\left(\frac{\lambda}{p_1^{1-\varrho_0}}\right)^{\nu_1}\cdots\left(\frac{\lambda}{p_s^{1-\varrho_0}}\right)^{\nu_s}\left(\log\frac{3x_s}{p_1^{\nu_1}\cdots p_s^{\nu_s}}\right)^{\beta-1}.
\end{align*}
It is not hard to see that the proof of (\ref{Equ:sumalphalog}) also gives
\[\sum_{z_1<p^{\nu}\le z_2}\left(\frac{\lambda}{p^{1-\varrho_0}}\right)^{\nu}\left(\log \frac{3z_2}{p^{\nu}}\right)^{\beta-1}\ll\frac{(\log(3z_2/z_1))^{\beta-1}}{z_1^{1-\varrho_0-\log_p\lambda}}\]
uniformly for all primes $p$ and all $0<z_1\le z_2$. Thus, we have
\begin{align*}
&\hspace*{4.8mm}\sum_{\substack{y_s<p_1^{\nu_1}\cdots p_s^{\nu_s}\le x_s\\\nu_1,...,\nu_s\ge1}}\nu_1^{\kappa l_1}\cdots \nu_s^{\kappa l_s}\left(\frac{\lambda}{p_1^{1-\varrho_0}}\right)^{\nu_1}\cdots\left(\frac{\lambda}{p_s^{1-\varrho_0}}\right)^{\nu_s}\left(\log\frac{3x_s}{p_1^{\nu_1}\cdots p_s^{\nu_s}}\right)^{\beta-1}\\
&\le2^{O(l)}(\log x)^{\kappa l}\sum_{\substack{y_s<p_1^{\nu_1}\cdots p_s^{\nu_s}\le x_s\\\nu_1,...,\nu_s\ge1}}\left(\frac{\lambda}{p_1^{1-\varrho_0}}\right)^{\nu_1}\cdots\left(\frac{\lambda}{p_s^{1-\varrho_0}}\right)^{\nu_s}\left(\log\frac{3x_s}{p_1^{\nu_1}\cdots p_s^{\nu_s}}\right)^{\beta-1}\\
&\le2^{O(l)}(\log x)^{\kappa l}\sum_{\substack{p_2^{\nu_2}\cdots p_{s}^{\nu_{s}}\le x_s\\\nu_2,...,\nu_s\ge1}}\left(\frac{\lambda}{p_2^{\log_{p_1}\lambda}}\right)^{\nu_2}\cdots\left(\frac{\lambda}{p_s^{\log_{p_1}\lambda}}\right)^{\nu_s}\frac{(\log(3x_s/y_s))^{\beta-1}}{y_s^{1-\varrho_0-\log_{p_1}\lambda}}\\
&\le\frac{2^{O(l)}(\log x)^{(\kappa+1) m+\beta-2}}{x^{\delta_0(1-\delta_0)\lambda'/2}(p_1\cdots p_s)^{\delta_0^2\lambda'}}\le\frac{2^{O(l)}(\log x)^{\beta-1}}{x^{(1-\delta_0)\lambda'/5}(p_1\cdots p_s)^{\delta_0^2\lambda'}},
\end{align*}
where the penultimate inequality follows from the previous line together with the observations that $p_i^{\log p_1\lambda}>\lambda$ for all $2\le i\le s$, that $x^{(1+\delta_0)/2}\ge(p_1\cdots p_s)^{1+\delta_0}$, and that
\[y_s^{1-\varrho_0-\log_{p_1}\lambda}\ge y_s^{\lambda'}=\left(x^{(1-\delta_0)/2}\cdot\frac{x^{(1+\delta_0)/2}}{p_1\cdots p_s}\right)^{\delta_0\lambda'}\ge x^{\delta_0(1-\delta_0)\lambda'/2}(p_1\cdots p_s)^{\delta_0^2\lambda'}.\]
It follows that
\begin{align*}
&\hspace*{4.8mm}\sum_{\substack{p_1\cdots p_sy_s<p_1^{\nu_1}\cdots p_s^{\nu_s}\le x\\\nu_1,...,\nu_s\ge2}}\nu_1^{\kappa l_1}\cdots \nu_s^{\kappa l_s}\frac{\alpha(p_1^{\nu_1})\cdots\alpha(p_s^{\nu_s})}{p_1^{\sigma_0\nu_1}\cdots p_s^{\sigma_0\nu_s}}\left(\log\frac{3x}{p_1^{\nu_1}\cdots p_s^{\nu_s}}\right)^{\beta-1}\\
&\le\frac{2^{O(l)}}{(p_1\cdots p_s)^{1-\varrho_0+\delta_0^2\lambda'}}x^{-(1-\delta_0)\lambda'/5}(\log x)^{\beta-1},
\end{align*}
from which we deduce that
\begin{align}\label{Equ:sum > p_1...p_sy_s}
&\hspace*{4.8mm}\sum_{p_1<\cdots<p_s\le \sqrt{x}}\sum_{\substack{p_1\cdots p_sy_s<p_1^{\nu_1}\cdots p_s^{\nu_s}\le x\\\nu_1,...,\nu_s\ge2}}\nu_1^{\kappa l_1}\cdots \nu_s^{\kappa l_s}\sum_{\substack{n\le x\\p_1^{\nu_1},...,p_s^{\nu_s}\parallel n}}\alpha(n)\left|\tilde{f}(n)-A(x)\right|^k\nonumber\\
&\le2^{O(m)}\lambda_{\alpha}x^{\sigma_0-(1-\delta_0)\lambda'/5}(\log x)^{k+\beta-1}\sum_{p_1<\cdots<p_s\le \sqrt{x}}\frac{1}{(p_1\cdots p_s)^{1-\varrho_0+\delta_0^2\lambda'}}\nonumber\\
&\le\frac{1}{s!}\lambda_{\alpha}x^{\sigma_0-(1-\delta_0)\lambda'/6}(\log x)^{\beta-1}.
\end{align}
On the other hand, if $p_1^{\nu_1}\cdots p_s^{\nu_s}\le p_1\cdots p_sy_s$, then $x^{(1-\delta_0)/2}\le x/(p_1^{\nu_1}\cdots p_s^{\nu_s})\le x$. Thus, we can apply the asymptotic formulas for (\ref{Equ:GMomentM(y;m)}) with $\eta_0=(1-\delta_0)/2$ and $y=x/(p_1^{\nu_1}\cdots p_s^{\nu_s})$, in conjunction with the Cauchy--Schwarz inequality, to estimate the inner sum in (\ref{Equ:sumalpha(n)|tilde{f}(n)-A(x)|^k}).
As a consequence, we have
\[\sum_{n\le x/\left(p_1^{\nu_1}\cdots p_s^{\nu_s}\right)}\alpha(n)\left|\tilde{f}(n)-A(x)\right|^a\ll\frac{2^{O(m-a)}\lambda_{\alpha}C_m}{\left(\sqrt{m}\right)^{m-a}}B(x)^{\frac{a}{2}}\left(\frac{x}{p_1^{\nu_1}\cdots p_s^{\nu_s}}\right)^{\sigma_0}\left(\log\frac{x}{p_1^{\nu_1}\cdots p_s^{\nu_s}}\right)^{\beta-1}.\]
Inserting this into (\ref{Equ:sumalpha(n)|tilde{f}(n)-A(x)|^k}) shows that the sum
\[\sum_{\substack{n\le x\\p_1^{\nu_1},...,p_s^{\nu_s}\parallel n}}\alpha(n)\left|\tilde{f}(n)-A(x)\right|^k\]
is
\begin{align*}
&\le\frac{2^{O(m-k)}\lambda_{\alpha}C_m}{\left(\sqrt{m}\right)^{m-k}}\cdot\frac{\alpha(p_1^{\nu_1})\cdots\alpha(p_s^{\nu_s})}{p_1^{\sigma_0\nu_1}\cdots p_s^{\sigma_0\nu_s}}\left(\sqrt{B(x)}+O\left(\frac{1}{\sqrt{m}}\sum_{i=1}^s|f(p_i)|\right)\right)^kx^{\sigma_0}(\log x)^{\beta-1}\\
&=\frac{2^{O(l)}\lambda_{\alpha}C_m}{m^{l/2}}\cdot\frac{\alpha(p_1^{\nu_1})\cdots\alpha(p_s^{\nu_s})}{p_1^{\sigma_0\nu_1}\cdots p_s^{\sigma_0\nu_s}}B(x)^{\frac{k}{2}}\left(1+O\left(\sqrt{\frac{m}{B(x)}}\right)\right)x^{\sigma_0}(\log x)^{\beta-1}\\
&\le\frac{2^{O(l)}\lambda_{\alpha}C_m}{m^{l/2}}\cdot\frac{\alpha(p_1^{\nu_1})\cdots\alpha(p_s^{\nu_s})}{p_1^{\sigma_0\nu_1}\cdots p_s^{\sigma_0\nu_s}}B(x)^{\frac{k}{2}}x^{\sigma_0}(\log x)^{\beta-1}.
\end{align*}
Note that
\begin{align*}
&\hspace*{4.8mm}\sum_{\substack{p_1^{\nu_1}\cdots p_s^{\nu_s}\le p_1\cdots p_sy_s\\\nu_1,...,\nu_s\ge2}}\nu_1^{\kappa l_1}\cdots \nu_s^{\kappa l_s}\frac{\alpha(p_1^{\nu_1})\cdots\alpha(p_s^{\nu_s})}{p_1^{\sigma_0\nu_1}\cdots p_s^{\sigma_0\nu_s}}\\
&\le2^{O(l)}\sum_{\substack{p_1^{\nu_1}\cdots p_s^{\nu_s}\le p_1\cdots p_sy_s\\\nu_1,...,\nu_s\ge2}}\nu_1^{\kappa l_1}\cdots \nu_s^{\kappa l_s}\left(\frac{\lambda}{p_1^{1-\varrho_0}}\right)^{\nu_1}\cdots\left(\frac{\lambda}{p_s^{1-\varrho_0}}\right)^{\nu_s}\\
&\le\frac{2^{O(l)}}{(p_1\cdots p_s)^{1-\varrho_0}}\sum_{\substack{p_1^{\nu_1}\cdots p_s^{\nu_s}\le y_s\\\nu_1,...,\nu_s\ge1}}\nu_1^{\kappa l_1}\cdots \nu_s^{\kappa l_s}\left(\frac{\lambda}{p_1^{1-\varrho_0}}\right)^{\nu_1}\cdots\left(\frac{\lambda}{p_s^{1-\varrho_0}}\right)^{\nu_s}\\
&\le \frac{2^{O(l)}}{(p_1\cdots p_s)^{1-\varrho_0}}\prod_{i=1}^{s}\Li_{-\lceil\kappa l_i\rceil}\left(\lambda/p_i^{1-\varrho_0}\right),
\end{align*}
where 
\[\Li_{-\ell}(\zeta)\colonequals\sum_{n=1}^{\infty}n^{\ell}\zeta^n\]
is the polylogarithm function of order $-\ell$ and complex argument $\zeta$ with $|\zeta|<1$, where $\ell\ge0$ is any integer. For example, $\Li_{0}(\zeta)=\zeta/(1-\zeta)$ and $\Li_{-1}(\zeta)=\zeta/(1-\zeta)^2$. The function $\Li_{-\ell}(\zeta)$ can be expressed in terms of the Eulerian polynomial $A_{\ell}(\zeta)$:
\[\Li_{-\ell}(\zeta)=\frac{\zeta A_{\ell}(\zeta)}{(1-\zeta)^{\ell+1}},\]
where
\[A_{\ell}(\zeta)\colonequals\sum_{j=0}^{\ell}\left\langle \ell\atop j\right\rangle \zeta^j\]
is the $\ell$th Eulerian polynomial, and 
\[\left\langle \ell\atop j\right\rangle\colonequals\sum_{a=0}^{j}(-1)^a\binom{\ell+1}{a}(j+1-a)^{\ell}\]
is the $j$th Eulerian number of size $\ell$. Combinatorially, it is known that, for every $\ell\ge1$,
\[\left\langle \ell\atop j\right\rangle=\#\{\tau\in S_{\ell}\colon \tau\text{~has~exactly~}j\text{~ascents}\},\]
where $S_{\ell}$ is the set of all permutations of $\{1,...,\ell\}$. Using this combinatorial intepretation one finds that $A_{\ell}(1)=\#S_{\ell}=\ell!$. Since $l_1+\cdots+l_s=l\le m$, we have
\[\prod_{i=1}^{s}\Li_{-\lceil\kappa l_i\rceil}\left(\lambda/p_i^{1-\varrho_0}\right)\le \frac{2^{O(l)}\lceil\kappa l_1\rceil!\cdots \lceil\kappa l_s\rceil!}{(p_1\cdots p_s)^{1-\varrho_0}}=\frac{2^{O(l)}\left(l_1^{l_1}\cdots l_s^{l_s}\right)^{\kappa}}{(p_1\cdots p_s)^{1-\varrho_0}}\le\frac{2^{O(l)}m^{\kappa l}}{(p_1\cdots p_s)^{1-\varrho_0}},\]
by Stirling's formula. Hence, we obtain
\[\sum_{\substack{p_1^{\nu_1}\cdots p_s^{\nu_s}\le p_1\cdots p_sy_s\\\nu_1,...,\nu_s\ge2}}\nu_1^{\kappa l_1}\cdots \nu_s^{\kappa l_s}\frac{\alpha(p_1^{\nu_1})\cdots\alpha(p_s^{\nu_s})}{p_1^{\sigma_0\nu_1}\cdots p_s^{\sigma_0\nu_s}}\le\frac{2^{O(l)}m^{\kappa l}}{(p_1\cdots p_s)^{2(1-\varrho_0)}}.\]
It follows that
\begin{align*}
&\hspace*{4.8mm}\sum_{\substack{p_1^{\nu_1}\cdots p_s^{\nu_s}\le p_1\cdots p_sy_s\\\nu_1,...,\nu_s\ge2}}\nu_1^{\kappa l_1}\cdots \nu_s^{\kappa l_s}\sum_{\substack{n\le x\\p_1^{\nu_1},...,p_s^{\nu_s}\parallel n}}\alpha(n)\left|\tilde{f}(n)-A(x)\right|^k\\
&\le\frac{2^{O(l)}\lambda_{\alpha}C_mm^{\kappa l}}{m^{l/2}(p_1\cdots p_s)^{2(1-\varrho_0)}}B(x)^{\frac{k}{2}}x^{\sigma_0}(\log x)^{\beta-1}.
\end{align*}
Summing the above over $p_1<\cdots<p_s\le\sqrt{x}$, we arrive at
\begin{align*}
&\hspace*{4.8mm}\sum_{p_1<\cdots<p_s\le\sqrt{x}}\sum_{\substack{p_1^{\nu_1}\cdots p_s^{\nu_s}\le p_1\cdots p_sy_s\\\nu_1,...,\nu_s\ge2}}\nu_1^{\kappa l_1}\cdots \nu_s^{\kappa l_s}\sum_{\substack{n\le x\\p_1^{\nu_1},...,p_s^{\nu_s}\parallel n}}\alpha(n)\left|\tilde{f}(n)-A(x)\right|^k\\
&\le2^{O(l)}\lambda_{\alpha}C_m m^{(\kappa-1/2)l}B(x)^{\frac{k}{2}}x^{\sigma_0}(\log x)^{\beta-1}\sum_{p_1<\cdots<p_s\le \sqrt{x}}\frac{1}{(p_1\cdots p_s)^{2(1-\varrho_0)}}\\
&\le\frac{2^{O(l)}}{s!}\lambda_{\alpha}C_m m^{(\kappa-1/2)l} B(x)^{\frac{k}{2}}x^{\sigma_0}(\log x)^{\beta-1},
\end{align*}
since $\varrho_0\in[0,1/2)$. Combining this estimate with (\ref{Equ:sum > p_1...p_sy_s}), we obtain
\begin{align*}
&\hspace*{4.8mm}\sum_{p_1<\cdots<p_s\le\sqrt{x}}\sum_{\substack{p_1^{\nu_1}\cdots p_s^{\nu_s}\le x\\\nu_1,...,\nu_s\ge2}}\nu_1^{\kappa l_1}\cdots \nu_s^{\kappa l_s}\sum_{\substack{n\le x\\p_1^{\nu_1},...,p_s^{\nu_s}\parallel n}}\alpha(n)\left|\tilde{f}(n)-A(x)\right|^k\\
&\le\frac{2^{O(l)}}{s!}\lambda_{\alpha}C_m m^{(\kappa-1/2)l} B(x)^{\frac{k}{2}}x^{\sigma_0}(\log x)^{\beta-1}.
\end{align*}
Therefore, (\ref{Equ:sum(tilde{f}-A(x))^k(f(n)-tilde{f}(n))^l}) is bounded above by
\begin{align*}
&\hspace*{4.8mm}2^{O(l)}\lambda_{\alpha}C_m m^{(\kappa-1/2)l}B(x)^{\frac{k}{2}}x^{\sigma_0}(\log x)^{\beta-1}\sum_{s\leq l}\frac{1}{s!}\sum_{\substack{l_1+\cdots+l_s=l\\l_1,...,l_s\in\N}}\binom{l}{l_1,...,l_s}\\
&\le2^{O(l)}\lambda_{\alpha}C_m m^{(\kappa-1/2)l} B(x)^{\frac{k}{2}}x^{\sigma_0}(\log x)^{\beta-1}T_l(1)\\
&\le2^{O(l)}C_m m^{(\kappa+1/2)l}B(x)^{\frac{k}{2}}S(x),
\end{align*}
which allows us to conclude that
\[\sum_{k=0}^{m-1}\binom{m}{k}\sum_{n\leq x}\alpha(n)\left(\tilde{f}(n)-A(x)\right)^k\left(f(n)-\tilde{f}(n)\right)^{m-k}\ll C_m m^{\kappa+\frac{3}{2}}B(x)^{\frac{m-1}{2}}S(x),\]
provided that in addition, we also have $1\le m\ll B(x)^{1/(2\kappa+3)}$. Inserting the above estimate and the estimate for the term corresponding to $k=m$ into (\ref{Equ:from tilde{f} to f}) completes the proof of Theorem \ref{thmMT2V}. 

\medskip
\section{Proof of Theorem \ref{thmWeightedHal} (sketch)}\label{S:Proof of thmWeightedHal}
Now we outline the proof of Theorem \ref{thmWeightedHal}. The first step is to redefine $f_q(n)$ introduced in Section \ref{S:Moments}. Again, let us suppose that $A_0\in(0,1)$ and that $|f(p)|\le 1$ for all primes $p$. For every $q\in\N$ we define
\begin{align*}
\widebar{F}(\sigma_0,q)&\colonequals\prod_{p\mid q}(1-F(\sigma_0,p)),\\
\widetilde{F}(\sigma_0,q)&\colonequals\frac{\rho_g(q)}{\varphi(q)}\widebar{F}(\sigma_0,q).
\end{align*}
For each prime $p$ we put
\[f_p(n)\colonequals\begin{cases}
	~f(p)(1-\widetilde{F}(\sigma_0,p)),&\text{~~~~~if $p\mid n$},\\
	~-f(p)\widetilde{F}(\sigma_0,p),&\text{~~~~~otherwise}.
\end{cases}\]
And as before, we set
\[f_q(n)\colonequals \prod_{p^{\nu}\parallel q}f_{p}(n)^{\nu}\]
for any $q\in\N$. In addition, let $c_g\in\N$ be the least positive integer such that $c_gg(x)\in\Z[x]$, and let $Q_0>c_g|g(0)|\ge1$ be such that (\ref{Equ:F(sigma_0,p)}) holds. Then for each $q\in\N$ with $P^-(q)>Q_0$ we have $\mathcal{Z}_g(q)=\mathcal{Z}^{\ast}_g(q)\subseteq(\Z/q\Z)^{\times}$ and $\rho_g(q)=\#\mathcal{Z}^{\ast}_g(q)$, where $\mathcal{Z}_g(q)$ denotes the zero locus of $g$ in $\Z/q\Z$. In particular, we have $0\le\rho_g(q)\le\varphi(q)$, which implies that $0\le\widetilde{F}(\sigma_0,q)\le1$ and that $|f_q(n)|\le 1$ for all $n\in\N$.
\par Next, we need an analogue of Lemma \ref{lemApprox}. Let $x$ be sufficiently large and set $z\colonequals x^{\delta(x)/m}>Q_0$. Then we have
\[\sum_{Q_0<p\leq x}f(p)\widetilde{F}(\sigma_0,p)=A_{f,g}(x)+O(1)\]
by (\ref{thmMT2Equ0}), (\ref{Equ:F(sigma_0,p)}), and the facts that $\rho_g$ is bounded on prime powers and that $\sum_{p}\psi_0(p)<\infty$. It is easily seen that
\[f(g(n))-A_{f,g}(x)=\sum_{Q_0<p\leq z}f_p(g(n))+\sum_{\substack{p>z\\p\mid g(n)}}f(p)-\sum_{z<p\leq x}f(p)\widetilde{F}(\sigma_0,p)+O(1).\]
Note that 
\[\sum_{z<p\leq x}f(p)\widetilde{F}(\sigma_0,p)=A_{f,g}(x)-A_{f,g}(z)+O(1)\ll \log\left(\frac{m}{\delta(x)}+1\right).\]
Since $1\le g(n)\ll n^{d_g}$ uniformly for all $n\in\N$, where $d_g\colonequals \deg g\ge1$, we have
\[\sum_{\substack{p>z\\p\mid g(n)}}f(p)\ll \frac{m}{\delta(x)}.\]
It follows that
\[\sum_{n\leq x}\alpha(n)(f(g(n))-A_{f,g}(x))^m=\sum_{n\leq x}\alpha(n)\left(\sum_{Q_0<p\leq z}f_p(g(n))\right)^m+O(E_g(x;m)),\]
where 
\[E_g(x;m)\colonequals \sum_{k=0}^{m-1}\binom{m}{k}2^{O(m-k)}\left(m\delta(x)^{-1}\right)^{m-k}\sum_{n\leq x}\alpha(n)\left|\sum_{p\leq z}f_p(g(n))\right|^k.\]
\par Now we turn to the estimation of
\[\sum_{n\le x}\alpha(n)\left(\sum_{Q_0<p\le z}f_p(g(n))\right)^m=\sum_{Q_0<p_1,...,p_m\le z}\sum_{n\leq x}\alpha(n)f_{p_1\cdots p_m}(g(n)).\]
Let $q\in\N\cap[1,x^{\delta(x)}]$ with $\omega(q)\le m$, $P^-(q)>Q_0$ and $P^+(q)\le z$. Then we have
\[\sum_{n\leq x}\alpha(n)f_{q}(g(n))=\sum_{ab\mid R_q}f_q(a)\mu(b)\sum_{\substack{n\leq x\\ab\mid g(n)}}\alpha(n)=\sum_{ab\mid R_q}f_q(a)\mu(b)\sum_{c\in\mathcal{Z}^{\ast}_g(ab)}\sum_{\substack{n\leq x\\n\equiv c\,(\text{mod}\,ab)}}\alpha(n).\]
Thus in place of Lemma \ref{lemMV2}, we need to input in our analysis the information about the distribution of values of $\alpha(n)$ with $n$ restricted to reduced residue classes. By Lemma \ref{lemS(x;a)}, the innermost sum differs from
\begin{align*}
\frac{1}{\varphi(ab)}\sum_{\substack{n\le x\\\gcd(n,ab)=1}}\alpha(n)&=\frac{1}{\varphi(ab)}x^{\sigma_0}(\log x)^{\beta-1}\left(\lambda_{\alpha}(ab)+O\left(\frac{1}{(\log x)^{A_0}}\right)\right)\\
&=\frac{1}{\varphi(ab)}x^{\sigma_0}(\log x)^{\beta-1}\left(\lambda_{\alpha}\widebar{F}(\sigma_0,ab)+O\left(\frac{1}{(\log x)^{A_0}}\right)\right)\\
&=\frac{1}{\varphi(ab)}\lambda_{\alpha}x^{\sigma_0}(\log x)^{\beta-1}\left(\widebar{F}(\sigma_0,ab)+O\left(\frac{1}{(\log x)^{A_0}}\right)\right)
\end{align*}
by the amount $\Delta_{\alpha}(x;ab,c)$. Hence, we obtain
\[\sum_{n\leq x}\alpha(n)f_{q}(g(n))=\lambda_{\alpha}\left(\widetilde{G}_1(\sigma_0,q)+O\left(\frac{\widetilde{G}_2(\sigma_0,q)}{(\log x)^{A_0}}\right)\right)x^{\sigma_0}(\log x)^{\beta-1}+J_q(x),\]
where 
\begin{align*}
\widetilde{G}_1(\sigma_0,q)&\colonequals \sum_{ab\mid R_q}f_q(a)\mu(b)\widetilde{F}(\sigma_0.ab),\\
\widetilde{G}_2(\sigma_0,q)&\colonequals \sum_{ab\mid R_q}\frac{\rho_g(ab)}{\varphi(ab)}|f_q(a)|,\\
J_q(x)&\colonequals\sum_{ab\mid R_q}f_q(a)\mu(b)\sum_{c\in\mathcal{Z}^{\ast}_g(ab)}\Delta_{\alpha}(x;ab,c).
\end{align*}
It is clear that $\widetilde{G}_1$ and $\widetilde{G}_2$ are both multiplicative in $q$. Easy calculation shows that
\[\widetilde{G}_1(\sigma_0,p^\nu)=f(p)^{\nu}\widetilde{F}(\sigma_0,p)\left(1-\widetilde{F}(\sigma_0,p)\right)\left((-1)^{\nu}\widetilde{F}(\sigma_0,p)^{\nu-1}+(1-\widetilde{F}(\sigma_0,p))^{\nu-1}\right)\]
for any prime power $p^{\nu}$. In particular, we have $\widetilde{G}_1(\sigma_0,p)=0$, $|\widetilde{G}_1(\sigma_0,p^\nu)|\leq1/4$, and $\widetilde{G}_1(\sigma_0,p^\nu)\geq0$ when $2\mid\nu$.
Moreover, we have that
\[\widetilde{G}_1(\sigma_0,p^2)=f(p)^2\widetilde{F}(\sigma_0,p)\left(1-\widetilde{F}(\sigma_0,p)\right)=\rho_g(p)\frac{f(p)^2}{p}+O\left(\frac{F(\sigma_0,p)}{p}+\frac{1}{p^2}\right),\]
and that
\[|\widetilde{G}_1(\sigma_0,p^\nu)|\le |f(p)|^{\nu}\widetilde{F}(\sigma_0,p)\le \rho_g(p)\frac{f(p)^2}{\varphi(p)}=\rho_g(p)\frac{f(p)^2}{p}+O\left(\frac{1}{p^{2}}\right)\]
for all $p^{\nu}$ with $p>Q_0$ and $\nu\ge2$. In addition, the quantity
\[\left|\sum_{Q_0<p_1,...,p_m\le z}J_{p_1\cdots p_m}(x)\right|\]
is bounded above by
\[\sum_{\substack{q:\,\omega(q)\le m\\ P^+(q)\le z}}\mu(q)^2\sum_{c\in\mathcal{Z}^{\ast}_g(q)}|\Delta_{\alpha}(x;q,c)|\sum_{a\mid q}\sum_{s\le m}\sum_{\substack{k_1+\cdots+k_s=m\\k_1,...,k_s\in\N}}\binom{m}{k_1,...,k_s}\sum_{\substack{Q_0<p_1<\cdots<p_{s}\leq z\\ q\mid p_1\cdots p_s}}\left|f_{p_1^{k_1}\cdots p_s^{k_s}}(a)\right|.\]
The inner sum over $s$ is 
\begin{align*}
&\le\sum_{s=\omega(q)}^m\sum_{\substack{k_1+\cdots+k_s=m\\k_1,...,k_s\in\N}}\binom{m}{k_1,...,k_s}\frac{1}{(s-\omega(q))!}\left(\sum_{\substack{Q_0<p\le z\\ p\nmid q}}|f_p(a)|\right)^{s-\omega(q)}\\
&\le 2^{O(m)}\sum_{s=\omega(q)}^m\sum_{\substack{k_1+\cdots+k_s=m\\k_1,...,k_s\in\N}}\binom{m}{k_1,...,k_s}\frac{1}{(s-\omega(q))!}(\log\log x)^{s-\omega(q)}\\
&\le 2^{O(m)}T_m(\log\log x)\\
&\le 2^{O(m)}(\log\log x)^m.
\end{align*}
It follows by \eqref{Equ:Delta_{alpha}(x;q,a)} that
\begin{align*}
\left|\sum_{Q_0<p_1,...,p_m\le z}J_{p_1\cdots p_m}(x)\right|&\le2^{O(m)}(\log\log x)^m\sum_{\substack{q:\,\omega(q)\le m\\ P^+(q)\le z}}\mu(q)^2\sum_{c\in\mathcal{Z}^{\ast}_g(q)}|\Delta_{\alpha}(x;q,c)|\\
&\le2^{O(m)}(\log\log x)^mS(x)\exp\left({-(\log\log x)^{1/3+\epsilon_0}}\right),
\end{align*}
which is $o(S(x))$. These observations allow us to conclude the proof of Theorem \ref{thmWeightedHal} by arguing as in Sections \ref{S:G(z)&D(y,z)} and \ref{S:E(y,z,w;m)}. It is also clear from the last inequality above that the bound $S(x)\exp\left({-(\log\log x)^{1/3+\epsilon_0}}\right)$ in \eqref{Equ:Delta_{alpha}(x;q,a)} can be weakened to a complicated one involving $\delta(x)$ and $B_{f,g}(x)$.

\medskip
\section{Proofs of Theorem \ref{thmWeightedDH} and Corollary \ref{corWeightedKS} (sketch)}\label{S:Proof of Thm1.4&Cor1.5}
Now we outline the proof of Theorem \ref{thmWeightedDH}, which borrows the ideas from the proofs of Theorem \ref{thmMT2} and \cite[Theorem 1]{DH} with proper modifications. Let $0<\epsilon<\min(1,K)$, and take $z\colonequals x^{1/v}$ and 
\[w\colonequals \begin{cases}
	~x^{1/\log(v+2)},&\text{~~~~~if $\beta=1$},\\
	~x^{1/(\epsilon\log(v+2))},&\text{~~~~~if $\beta\ne1$},
\end{cases}\]
where we recall that $v\asymp m$ when $\beta=1$ and $v=(\log\log x)^{m(\vartheta_0+2)}$ when $\beta\ne1$ as chosen in Section \ref{S:E(y,z,w;m)}. Having made these choices, we have $\epsilon\log(v+2)\to\infty$ as $x\to\infty$ in the case $\beta\ne1$. Let
\begin{align*}
\mathcal{P}^-_{\epsilon}(x)&\colonequals\left\{p\le x\colon |f(p)|\le\epsilon\sqrt{B^{\ast}(x)}\right\},\\
\mathcal{P}^+_{\epsilon}(x)&\colonequals\left\{p\le x\colon \epsilon\sqrt{B^{\ast}(x)}<|f(p)|\le K\sqrt{B^{\ast}(x)}\right\},\\
\mathcal{P}_{\infty}(x)&\colonequals\left\{p\le x\colon |f(p)|>K\sqrt{B^{\ast}(x)}\right\},
\end{align*}
and put $\mathcal{P}_{K}(x)\colonequals\mathcal{P}^-_{\epsilon}(x)\cup\mathcal{P}^+_{\epsilon}(x)$. We consider the strongly additive function
\[f_{\epsilon}(n;x)\colonequals\sum_{\substack{p\mid n\\p\in\mathcal{P}^-_{\epsilon}(x)}}f(p)+\epsilon_{\beta,1}\sum_{\substack{p\mid n\\p\in\mathcal{P}^+_{\epsilon}(x)\cap(z,x]}}f(p)+\sum_{\substack{p\mid n\\p\in\mathcal{P}_{\infty}(x)}}f(p),\]
where we recall that $\epsilon_{\beta,1}$ takes value 0 if $\beta=1$ and 1 otherwise, and define 
\begin{align*}
A_{\epsilon}(x)&\colonequals\sum_{p\in\mathcal{P}^-_{\epsilon}(x)}\alpha(p)\frac{f(p)}{p^{\sigma_0}},\\
B_{\epsilon}(x)&\colonequals\sum_{p\in\mathcal{P}^-_{\epsilon}(x)}\alpha(p)\frac{f(p)^2}{p^{\sigma_0}}.
\end{align*}
By hypothesis, 
\[B(x)-B_{\epsilon}(x)=\sum_{\substack{p\le x\\|f(p)|>\epsilon\sqrt{B^{\ast}(x)}}}\alpha(p)\frac{f(p)^2}{p^{\sigma_0}}=o(B^{\ast}(x)),\]
and so
\[|A_{\epsilon}(x)-A(x)|\le\frac{1}{\epsilon\sqrt{B^{\ast}(x)}}\sum_{\substack{p\le x\\|f(p)|>\epsilon\sqrt{B^{\ast}(x)}}}\alpha(p)\frac{f(p)^2}{p^{\sigma_0}}=o\left(\epsilon^{-1}\sqrt{B^{\ast}(x)}\right).\]
We expect that the distribution of $f_{\epsilon}(n;x)$ is close to being Gaussian with mean $A(x)$ and variance $B(x)$ when $x$ gets sufficiently large. In what follows, we shall restrict our attention to the case $\beta\ne1$, since the opposite case $\beta=1$ is not only similar but also easier. Looking back at the proof of Lemma \ref{lemApprox}, we find, for sufficiently large $x$, that
\[\sum_{p\in\mathcal{P}^-_{\epsilon}(x)\cap(Q_0,x]}f(p)F(\sigma_0,p)=A_{\epsilon}(x)+O\left(\epsilon\sqrt{B^{\ast}(x)}\right)=A(x)+O\left(\epsilon\sqrt{B^{\ast}(x)}\right),\]
so that
\begin{equation}\label{Equ:f_{epsilon}(n;x)-A(x)}
f_{\epsilon}(n;x)-A(x)=\sum_{p\in\mathcal{P}^-_{\epsilon}(x)\cap(Q_0,z]}f_p(n)+\sum_{\substack{p\mid n\\p\in\mathcal{P}_K(x)\cap(z,w]}}f(p)+O\left(\epsilon\sqrt{B(x)}\right),
\end{equation}
where we have used the hypothesis that $f(n)=o(\sqrt{B(x)})$ for all $n\le x$ whose prime factors $p$ satisfy $|f(p)|>K\sqrt{B^{\ast}(x)}$. This leads to an analogue of Lemma \ref{lemApprox} in which 
the second sum above plays the same role as $\omega(n;z,w)$. To estimate the moments of $f_{\epsilon}(n;x)$, one only needs to recycle the arguments used in the proof of Theorem \ref{thmMT2} and make suitable modifications. For instance, the estimation of 
\[\sum_{n\leq y}\alpha(n)\left(\sum_{p\in\mathcal{P}^-_{\epsilon}(x)\cap(Q_0,z]}f_p(n)\right)^m\]
is essentially the same as that of (\ref{Equ:Main1}) given in Sections \ref{S:Moments} and \ref{S:G(z)&D(y,z)}, except that we use the inequality $|f(p)|\le\epsilon\sqrt{B^{\ast}(x)}$ for $p\in\mathcal{P}^-_{\epsilon}(x)$ in place of the bound $f(p)=O(1)$ throughout the argument. This way, we obtain
\begin{align}\label{Equ:Moment of f_p(n) in P^-}
\sum_{n\leq y}\alpha(n)\left(\sum_{p\in\mathcal{P}^-_{\epsilon}(x)\cap(Q_0,z]}f_p(n)\right)^m&=\lambda_{\alpha}\left(\mu_m+O\left(\frac{\epsilon\log v}{\log\log\log x}\right)\right)B(x)^{\frac{m}{2}}y^{\sigma_0}(\log y)^{\beta-1}\nonumber\\
&=\lambda_{\alpha}(\mu_m+O(\epsilon))B(x)^{\frac{m}{2}}y^{\sigma_0}(\log y)^{\beta-1}
\end{align}
uniformly for $y\in[x^{\eta_0},x]$, where $\eta_0\in(0,1]$ is any given constant. On the other hand, the estimation of the error involving the second sum in (\ref{Equ:f_{epsilon}(n;x)-A(x)}) is essentially the same as that of $E(y,z,w;m)$ in the case $\beta\ne1$ given in Section \ref{S:E(y,z,w;m)}. The only difference is that we now make use of the estimates that $|f(p)|\le K\sqrt{B^{\ast}(x)}$ for all $p\in\mathcal{P}_K(x)$ and that
\[\sum_{p\in\mathcal{P}_K(x)\cap(z,w]}\alpha(p)\frac{|f(p)|^{\nu}}{p^{\sigma_0}}\ll B^{\ast}(x)^{\frac{\nu-1}{2}}\sum_{p\in\mathcal{P}_K(x)\cap(z,w]}\alpha(p)\frac{|f(p)|}{p^{\sigma_0}}\ll\epsilon B^{\ast}(x)^{\frac{\nu}{2}}\log v\ll\epsilon B(x)^{\frac{\nu}{2}}\]
for all $\nu\ge1$, which can be easily seen by considering $p\in\mathcal{P}^-_{\epsilon}(x)$ and $p\in\mathcal{P}^+_{\epsilon}(x)$ separately, in place of the estimates that $f(p)=O(1)$ and that
\[\sum_{z<p\le w}\alpha(p)\frac{|f(p)|^{\nu}}{p^{\sigma_0}}=O(\log v)=O(\log\log\log x),\]
respectively. One shows in this way that the error involving the second sum in (\ref{Equ:f_{epsilon}(n;x)-A(x)}) is $O(\epsilon \lambda_{\alpha}B(x)^{\frac{m}{2}}y^{\sigma_0}(\log y)^{\beta-1})$. Combining this estimate with  (\ref{Equ:Moment of f_p(n) in P^-}) and taking $y=x$ yields
\[S(x)^{-1}\sum_{n\le x}\alpha(n)\left(f_{\epsilon}(n;x)-A(x)\right)^m=(\mu_m+O(\epsilon))B(x)^{\frac{m}{2}}\] 
for every fixed $m\in\N$ and all sufficiently large $x$, where
the implied constant in the error term is independent of $\epsilon$. 
\par To complete the proof of Theorem \ref{thmWeightedDH} for the case $\beta\ne1$, it is sufficient to show
\begin{equation}\label{Equ:f(n)-f_{epsilon}(n;x)1}
S(x)^{-1}\sum_{n\le x}\alpha(n)|f(n)-f_{\epsilon}(n;x)|^m=O\left(\epsilon B^{\ast}(x)^{\frac{m}{2}}\right)
\end{equation}
for every given $\epsilon\in(0,1)$ and $m\in\N$, where
the implicit constant in the error term is independent of $\epsilon$. Since the case where $m$ is odd follows from the case where $m$ is even by Cauchy--Schwarz, we need only to consider the latter case. The proof of this case is largely the same as that of \cite[Lemma 2]{DH}, except for the slight complication in the possible case $\beta\in(0,1)$. When $m$ is even, we have
\[S(x)^{-1}\sum_{n\le x}\alpha(n)|f(n)-f_{\epsilon}(n;x)|^m=S(x)^{-1}\sum_{n\le x}\alpha(n)\sum_{\substack{p_1,...,p_m\mid n\\p_1,...,p_m\in\mathcal{P}^+_{\epsilon}(x)\cap[2,z]}}f(p_1)\cdots f(p_m),\]
which, after grouping terms according to the distinct primes among $p_1,...,p_m$, becomes
\begin{equation}\label{Equ:f(n)-f_{epsilon}(n;x)2}
S(x)^{-1}\sum_{s\le m}\sum_{\substack{p_1<\cdots<p_s\le z\\p_1,...,p_s\in\mathcal{P}^+_{\epsilon}(x)}}\sum_{\substack{k_1+\cdots+k_s=m\\k_1,...,k_s\in\N}}\binom{m}{k_1,...,k_s}f(p_1)^{k_1}\cdots f(p_s)^{k_s}\sum_{\substack{n\le x\\p_1\cdots p_s\mid n}}\alpha(n).
\end{equation}
By (\ref{Equ:S(x)}) we have
\[\sum_{\substack{n\le x\\p_1\cdots p_s\mid n}}\alpha(n)=\sum_{\substack{q\le x\\R_q=p_1\cdots p_s}}\alpha(q)\sum_{\substack{n'\le x/q\\\gcd(n',q)=1}}\alpha(n')\ll \lambda_{\alpha}x^{\sigma_0}\sum_{\substack{q\le x\\R_q=p_1\cdots p_s}}\frac{\alpha(q)}{q^{\sigma_0}}\left(\log\frac{3x}{q}\right)^{\beta-1}.\]
Appealing to (\ref{Equ:sumalphalog}) we derive
\begin{align*}
\sum_{\substack{q\le x\\R_q=p_1\cdots p_s}}\frac{\alpha(q)}{q^{\sigma_0}}\left(\log\frac{3x}{q}\right)^{\beta-1}&\ll(\log x)^{\beta-1}\sum_{\substack{q\le \sqrt{x}\\R_q=p_1\cdots p_s}}\frac{\alpha(q)}{q^{\sigma_0}}+\sum_{\substack{\sqrt{x}<q\le x\\R_q=p_1\cdots p_s}}\frac{\alpha(q)}{q^{\sigma_0}}\left(\log\frac{3x}{q}\right)^{\beta-1}\\
&\ll(\log x)^{\beta-1}\prod_{i=1}^{s}\sum_{\nu=1}^{\infty}\frac{\alpha(p_i^{\nu})}{p_i^{\sigma_0\nu}}+\frac{(\log x)^{s+\beta-2}}{\left(\sqrt{x}\right)^{1-\varrho_0}}\\
&=(\log x)^{\beta-1}\prod_{i=1}^{s}\left(\frac{\alpha(p_i)}{p_i^{\sigma_0}}+\psi_0(p_i)\right)+\frac{(\log x)^{s+\beta-2}}{\left(\sqrt{x}\right)^{1-\varrho_0}}.
\end{align*}
These estimates together with (\ref{Equ:S(x)}) imply that (\ref{Equ:f(n)-f_{epsilon}(n;x)2}) is $\ll\Sigma_1+\Sigma_2$, where
\begin{align*}
\Sigma_1&\colonequals \sum_{s\le m}\sum_{\substack{p_1<\cdots<p_s\le z\\p_1,...,p_s\in\mathcal{P}^+_{\epsilon}(x)}}\sum_{\substack{k_1+\cdots+k_s=m\\k_1,...,k_s\in\N}}\binom{m}{k_1,...,k_s}\left|f(p_1)^{k_1}\cdots f(p_s)^{k_s}\right|\prod_{i=1}^{s}\left(\frac{\alpha(p_i)}{p_i^{\sigma_0}}+\psi_0(p_i)\right),\\
\Sigma_2&\colonequals\frac{(\log x)^{m-1}}{\left(\sqrt{x}\right)^{1-\varrho_0}}\sum_{s\le m}\sum_{\substack{p_1<\cdots<p_s\le z\\p_1,...,p_s\in\mathcal{P}^+_{\epsilon}(x)}}\sum_{\substack{k_1+\cdots+k_s=m\\k_1,...,k_s\in\N}}\binom{m}{k_1,...,k_s}\left|f(p_1)^{k_1}\cdots f(p_s)^{k_s}\right|.
\end{align*}
Since $f(p)\le K\sqrt{B^{\ast}(x)}$ for all $p\in\mathcal{P}^+_{\epsilon}(x)$, we have
\[\Sigma_2\ll\frac{(\log x)^{m-1}}{\left(\sqrt{x}\right)^{1-\varrho_0}}\pi(z)^mB^{\ast}(x)^{\frac{m}{2}}=o\left(B^{\ast}(x)^{\frac{m}{2}}\right)\ll\epsilon B^{\ast}(x)^{\frac{m}{2}}.\]
To bound $\Sigma_1$, we observe
\[\left|f(p_1)^{k_1}\cdots f(p_s)^{k_s}\right|\ll B^{\ast}(x)^{\frac{m-s}{2}}|f(p_1)\cdots f(p_s)|.\]
Thus, we have
\begin{align*}
\Sigma_1&\le\sum_{s\le m}B^{\ast}(x)^{\frac{m-s}{2}}\frac{1}{s!}\left(\sum_{\substack{p\le z\\ p\in\mathcal{P}^+_{\epsilon}(x)}}\left(\alpha(p)\frac{|f(p)|}{p^{\sigma_0}}+\psi_0(p)\right)\right)^s\sum_{\substack{k_1+\cdots+k_s=m\\k_1,...,k_s\in\N}}\binom{m}{k_1,...,k_s}\\
&=\sum_{s\le m}B^{\ast}(x)^{\frac{m-s}{2}}\frac{1}{s!}\left(o\left(\epsilon^{-1}\sqrt{B^{\ast}(x)}\right)\right)^s\sum_{\substack{k_1+\cdots+k_s=m\\k_1,...,k_s\in\N}}\binom{m}{k_1,...,k_s}\ll\epsilon B^{\ast}(x)^{\frac{m}{2}}.
\end{align*}
Combining these estimates completes the proof of (\ref{Equ:f(n)-f_{epsilon}(n;x)1}) in the case $\beta\ne1$. 
\par As we mentioned in Section \ref{S:MR}, Corollary \ref{corWeightedKS} is an immediate consequence of Theorem \ref{thmWeightedDH} when $f$ is strongly additive. The transition to the general additive case is then accomplished by applying the following analogue of \cite[Theorem B]{Sha}. And this is the only place where we need to make use of characteristic functions.
\begin{lem}\label{lem:f and tilde{f}}
Let $f\colon\N\to\R$ be an additive function, and let $\alpha\in\mathcal{M}^{\ast}$ with parameters $A_0,\beta,\sigma_0,\vartheta_0,\varrho_0,r$. Denote by $\tilde{f}$ the strongly additive contraction of $f$. Suppose that $B(x)\to\infty$ as $x\to\infty$. Then $X_N(n)\colonequals(f(n)-A(N))/\sqrt{B(N)}$ possesses a limiting distribution function with respect to the natural probability measure induced by $\alpha$ if and only if $\widetilde{X}_N(n)\colonequals(\tilde{f}(n)-A(N))/\sqrt{B(N)}$ does, in which case they share the same limiting distribution function.
\end{lem}
\begin{proof}
As before, we shall assume $A_0\in(0,1)$. For each $N\in\N$, the distribution functions of $X_N(n)$ and $\widetilde{X}_N(n)$ are given by 
\begin{align*}
\Phi_{N}(V)&=S(N)^{-1}\sum_{\substack{n\le N\\X_N\le V}}\alpha(n),\\
\widetilde{\Phi}_{N}(V)&=S(N)^{-1}\sum_{\substack{n\le N\\\widetilde{X}_N\le V}}\alpha(n),
\end{align*}
respectively. We have to show that $\Phi_{N}(V)$ converges weakly to a distribution function as $N\to\infty$ if and only if $\widetilde{\Phi}_{N}(V)$ does, in which case they converge weakly to the same distribution function. Note that the characteristic functions of $X_N(n)$ and  $\widetilde{X}_N(n)$ are 
\begin{align*}
\varphi_{N}(t)&=S(N)^{-1}\sum_{n\le N}\alpha(n)e^{itX_N(n)},\\
\widetilde{\varphi}_{N}(t)&=S(N)^{-1}\sum_{n\le N}\alpha(n)e^{it\widetilde{X}_N(n)},
\end{align*}
respectively. By L\'{e}vy's continuity theorem \cite[Theorem III.2.6]{TenenText}, it suffices to show 
\begin{equation}\label{Equ:phi_{N}-tilde{phi}_{N}}
\lim\limits_{N\to\infty}\left(\varphi_{N}(t)-\widetilde{\varphi}_{N}(t)\right)=0
\end{equation}
for any given $t\in\R$. To prove this, let us fix $t\in\R$ and let $\epsilon\in(0,1/(2|t|+1))$ be arbitrary. Denote by $J_{\epsilon}(N)$ the greatest integer not exceeding $\sqrt{N}$ such that the inequality $|f(n)|\le\epsilon\sqrt{B(N)}$ holds for all $1\le n\le J_{\epsilon}(N)$. Since $B(N)\nearrow\infty$ as $N\to\infty$, we have $J_{\epsilon}(N)\nearrow\infty$ as $N\to\infty$. By (\ref{Equ:S(x)}) we have
\begin{align*}
\left|\varphi_{N}(t)-\widetilde{\varphi}_{N}(t)\right|&\le S(N)^{-1}\sum_{n\le N}\alpha(n)\left|\exp\left(it\frac{f(n)-\tilde{f}(n)}{\sqrt{B(N)}}\right)-1\right|\\
&=S(N)^{-1}\sum_{\substack{a\le N\\a \text{~squareful}}}\alpha(a)\left|\exp\left(it\frac{f(a)-f(R_a)}{\sqrt{B(N)}}\right)-1\right|\sum_{\substack{b\le N/a\\b \text{~squarefree}\\\gcd(b,a)=1}}\alpha(b)\\
&\ll S(N)^{-1}\lambda_{\alpha}N^{\sigma_0}\sum_{\substack{a\le N\\a \text{~squareful}}}\frac{\alpha(a)}{a^{\sigma_0}}\left(\log\frac{3N}{a}\right)^{\beta-1}\left|\exp\left(it\frac{f(a)-f(R_a)}{\sqrt{B(N)}}\right)-1\right|.
\end{align*}
From (\ref{thmMT2Equ0}) and (\ref{thmMT2Equ2}) it follows that
\[\sum_{\substack{a=1\\a \text{~squareful}}}^{\infty}\frac{\alpha(a)}{a^{s}}=\prod_{p}\left(1+\sum_{\nu\ge2}\frac{\alpha(p^{\nu})}{p^{\nu s}}\right)\]
is absolutely convergent for $s\in\C$ with $\Re(s)>\max(\varrho_0,r)+\sigma_0-1$. Thus
\[c(\delta)\colonequals\sum_{\substack{a=1\\a \text{~squareful}}}^{\infty}\frac{\alpha(a)}{a^{\sigma_0-\delta}}=\prod_{p}\left(1+\sum_{\nu\ge2}\frac{\alpha(p^{\nu})}{p^{\nu (\sigma_0-\delta)}}\right)<\infty\]
for any $\delta<1-\max(\varrho_0,r)$. Since
\[\left|it\frac{f(a)-f(R_a)}{\sqrt{B(N)}}\right|\le 2\epsilon|t|<1\] 
for all $a\le J_{\epsilon}(N)$, this implies 
\[\sum_{\substack{a\le J_{\epsilon}(N)\\a \text{~squareful}}}\frac{\alpha(a)}{a^{\sigma_0}}\left(\log\frac{3N}{a}\right)^{\beta-1}\left|\exp\left(it\frac{f(a)-f(R_a)}{\sqrt{B(N)}}\right)-1\right|\ll\epsilon|t|(\log N)^{\beta-1}.\]
Now fix $0<\delta<1-\max(\varrho_0,r)$. By partial summation we have
\[\sum_{\substack{a\le x\\a \text{~squareful}}}\frac{\alpha(a)}{a^{\sigma_0}}=c(0)-\int_{x}^{\infty}\frac{1}{t^{\delta}}\,d\left(\sum_{\substack{a\le t\\a \text{~squareful}}}\frac{\alpha(a)}{a^{\sigma_0-\delta}}\right)=c(0)+o\left(x^{-\delta}\right)\]
when $x$ is sufficiently large. It follows that
\begin{align*}
&\hspace*{4.8mm}\sum_{\substack{J_{\epsilon}(N)<a\le N\\a \text{~squareful}}}\frac{\alpha(a)}{a^{\sigma_0}}\left(\log\frac{3N}{a}\right)^{\beta-1}\left|\exp\left(it\frac{f(a)-f(R_a)}{\sqrt{B(N)}}\right)-1\right|\\
&\le2\sum_{\substack{J_{\epsilon}(N)<a\le N\\a \text{~squareful}}}\frac{\alpha(a)}{a^{\sigma_0}}\left(\log\frac{3N}{a}\right)^{\beta-1}\\
&=2\int_{J_{\epsilon}(N)}^{N}\left(\log\frac{3N}{t}\right)^{\beta-1}\,d\left(\sum_{\substack{a\le t\\a \text{~squareful}}}\frac{\alpha(a)}{a^{\sigma_0}}\right)\\
&=o\left(N^{-\delta}\right)+o\left((\log N)^{\beta-1}J_{\epsilon}(N)^{-\delta}\right)+o\left(\int_{J_{\epsilon}(N)}^{N}t^{-1-\delta}\left(\log\frac{3N}{t}\right)^{\beta-2}\,dt\right)
\end{align*}
for sufficiently large $N$. By a change of variable we see that
\begin{align*}
\int_{J_{\epsilon}(N)}^{N}t^{-1-\delta}\left(\log\frac{3N}{t}\right)^{\beta-2}\,dt&=(3N)^{-\delta}\int_{\log 3}^{\log(3N/J_{\epsilon}(N))}e^{\delta t}t^{\beta-2}\,dt\\
&\ll(3N)^{-\delta}\left(\frac{3N}{J_{\epsilon}(N)}\right)^{\delta}\left(\log\frac{3N}{J_{\epsilon}(N)}\right)^{\beta-2}\\
&\ll(\log N)^{\beta-2}J_{\epsilon}(N)^{-\delta}.
\end{align*}
Hence, we have
\[\sum_{\substack{J_{\epsilon}(N)<a\le N\\a \text{~squareful}}}\frac{\alpha(a)}{a^{\sigma_0}}\left(\log\frac{3N}{a}\right)^{\beta-1}\left|\exp\left(it\frac{f(a)-f(R_a)}{\sqrt{B(N)}}\right)-1\right|=o\left((\log N)^{\beta-1}J_{\epsilon}(N)^{-\delta}\right).\]
for sufficiently large $N$. Gathering the estimates above, we obtain
\[\varphi_{N}(t)-\widetilde{\varphi}_{N}(t)\ll\epsilon|t|+o\left(J_{\epsilon}(N)^{-\delta}\right)\]
for sufficiently large $N$, where the implicit constants are independent of $t$, $\epsilon$ and $N$. From this estimate we infer that
\[\limsup_{N\to\infty}\left|\varphi_{N}(t)-\widetilde{\varphi}_{N}(t)\right|=O(\epsilon|t|),\] 
where the implicit constant is independent of $t$ and $\epsilon$. Since $\epsilon\in(0,1/(2|t|+1))$ is arbitrary, we obtain (\ref{Equ:phi_{N}-tilde{phi}_{N}}) as desired.
\end{proof}

\begin{rmk}\label{rmk8.1}
Let $\alpha(n)=\tau(n)^2/n^{11}$, where $\tau$ is Ramanujan's $\tau$-function, and define the additive function $f(n)$ by $f(p^{\nu})=\log\sqrt{\alpha(p^{\nu})}$ if $\alpha(p^{\nu})\ne0$ and $f(p^{\nu})=0$ otherwise, where $p^{\nu}$ is any prime power. Then $\alpha(n)$ satisfies conditions (i)--(iv) with any fixed $A_0>0$, $\beta=1$, $\sigma_0=1$, $\vartheta_0=0$, and any fixed $\varrho_0\in(0,1)$ and $r\in(1/2,1)$. Moreover, we have $\alpha(n)\le d(n)^2$ by Deligne's bound \cite{Deligne}. As alluded to in Section \ref{S:MR}, Elliott \cite{Elliott1} showed, using ideas from probability theory, that the limiting distribution of $(f(n)-A(x))/\sqrt{B(x)}$ with respect to the natural probability measure induced by $\alpha$ is the standard Gaussian distribution. In fact, we can derive his result from Corollary \ref{corWeightedKS} in combination with Lemma \ref{lem:f and tilde{f}} and \cite[Lemma 7]{Elliott1} without difficulty. In comparison to Elliott's probabilistic approach, our approach enables us to get around some of the complications resulting from the analysis of $\tau(n)$.
\par To illustrate this, let us consider the strongly additive function $f_0(n)$ defined by $f_0(p)=\log\sqrt{\alpha(p)}$ if $p\notin E_0$ and $f_0(p)=0$ otherwise, where $E_0\colonequals\{p>2\colon\alpha(p)\le\exp(-2\sqrt[3]{\log\log p})\}$. Denote by $A_0(x)$ and $B_0(x)$ the expected mean and variance of $f_0(n)$ weighted by $\alpha(n)$, respectively. It can be shown \cite[Lemma 7]{Elliott1} that $B(x)\asymp\log\log x$. Since the inequality $t|\log t|\le\sqrt{t}$ holds for all $t\in[0,1]$, we have 
\[\sum_{\substack{p\le x\\p\in E_0}}\alpha(p)\frac{|f(p)|}{p}\le\sum_{\substack{p\le x\\p\in E_0}}\frac{\sqrt{\alpha(p)}}{p}\le \sum_{p>2}\frac{1}{p}\exp\left(-\sqrt[3]{\log\log p}\right)<\infty.\]
It follows that $A_0(x)=A(x)+O(1)$. A similar argument shows that $B_0(x)=B(x)+O(1)\asymp\log\log x$. Thus, $f_0(p)=O(B_0(p)^{1/3})$ for all $p$, which shows that $f_0(n)$ satisfies the hypotheses in Corollary \ref{corWeightedKS}. Hence, the limiting distribution of $(f_0(n)-A(x))/\sqrt{B(x)}$ with respect to the natural probability measure induced by $\alpha$ is the standard Gaussian distribution. 
\par To complete our argument, let $\widetilde{f}$ be the strongly additive contraction of $f$. Then $f_0(n)\ge \widetilde{f}(n)$ for all $n\in\N$. Moreover, Deligne's bound and the fact that $\tau(n)\in\Z$ for all $n\in\N$ imply that $-(11\log p)/2\le f(p)\le\log 2$ whenever $\alpha(p)\ne0$. Since
\[\sum_{\nu\ge1}\frac{\alpha(p^{\nu})}{p^{\nu}}\le\frac{\alpha(p)}{p}+\sum_{\nu\ge2}\frac{(\nu+1)^2}{p^{\nu}}=\frac{\alpha(p)}{p}+O\left(\frac{1}{p^2}\right),\]
we have
\begin{align*}
S(x)^{-1}\sum_{n\le x}\alpha(n)\left(f_0(n)-\widetilde{f}(n)\right)&=S(x)^{-1}\sum_{\substack{p\le x\\p\in E_0}}|f(p)|\sum_{\substack{n\le x\\p\mid n}}\alpha(n)\\
&=S(x)^{-1}\sum_{\substack{p\le x\\p\in E_0}}|f(p)|\sum_{\nu\ge1}\alpha(p^{\nu})\sum_{\substack{n'\le x/p^{\nu}\\p\nmid n'}}\alpha(n')\\
&\ll S(x)^{-1}\lambda_{\alpha}x\sum_{\substack{p\le x\\p\in E_0}}|f(p)|\sum_{\nu\ge1}\frac{\alpha(p^{\nu})}{p^{\nu}}\\
&\ll \sum_{\substack{p\le x\\p\in E_0}}\alpha(p)\frac{|f(p)|}{p}+O\left(\sum_{p>2}\frac{\log p}{p^2}\right)\ll 1.
\end{align*}
This estimate is sufficient for us to conclude that the limiting distribution of $(\widetilde{f}(n)-A(x))/\sqrt{B(x)}$ with respect to the natural probability measure induced by $\alpha$ is also the standard Gaussian distribution. By Lemma \ref{lem:f and tilde{f}}, the same is true for $(f(n)-A(x))/\sqrt{B(x)}$.
\end{rmk}

\medskip
\section{Concluding Remarks}\label{S:Final}
Although in the present paper we only focused on the subclass $\mathcal{M}^{\ast}$ of multiplicative functions, it is also of interest to consider weight functions $\alpha(n)$ which satisfy certain Landau--Selberg--Delange type conditions. Given more information about $\alpha(n)$ and its associated Dirichlet series $F(s)=\sum_{n=1}^{\infty}\alpha(n)n^{-s}$, better results are obtainable in some circumstances. Below we give a brief description of the method in the special case where $F(s)$ is close to an integral power of the Riemann zeta-function $\zeta(s)$. 
\par For a complex number $s\in\C$, we write $\sigma=\Re(s)$ and $t=\Im(s)$. Let $\alpha\colon\N\to\R_{\geq0}$ be a multiplicative function whose Dirichlet series $F(s)=\sum_{n=1}^{\infty}\alpha(n)n^{-s}$ is absolutely convergent for $s\in\C$ with $\sigma>\sigma_0$, where $\sigma_0>0$ is constant. Suppose that there exist constants $\beta\in\N$, $0<\theta_0<\sigma_0$, $B>0$, and $0<\delta<1$, such that $H_{\beta}(s)\colonequals F(s)\zeta(s-\sigma_0+1)^{-\beta}$ has an analytic continuation in the half plane $\sigma\ge\theta_0$ with \[\lim\limits_{s\to\sigma_0}F(s)(s-\sigma_0)^{\beta}>0,\]
and such that $|H_{\beta}(s)|\leq B(1+|t|)^{1-\delta}$ for all $s\in\C$ with $\sigma\ge\theta_0$. It is clear that $F(s)$ has (absolute) abscissa of convergence $\sigma_0$. Adapting the argument used in the proof of \cite[Lemma 2.1]{KMS} or \cite[Theorem II.5.2]{TenenText}, one can show that there exists some constant $\epsilon_0>0$ such that 
\begin{equation}\label{Equ:S(x) via S-D}
S(x)=\frac{1}{\sigma_0}\Res_{s=\sigma_0}\left(\frac{F(s)x^s}{s-\sigma_0+1}\right)-x^{\sigma_0}(\log x)^{\beta-1}\sum_{k=1}^{\beta-1}\sum_{j=0}^{k-1}c_{j,k}\frac{\mu_j(\beta)}{(\log x)^k}+O\left(Bx^{\theta}\right)
\end{equation}
uniformly for all $x\geq3$ and $\theta\in(\sigma_0-\epsilon_0,\sigma_0)$, where
\begin{align*}
\mu_k(\beta)&\colonequals\frac{1}{k!}\cdot\frac{d^k}{ds^k}\left(\frac{F(s)(s-\sigma_0)^{\beta}}{s-\sigma_0+1}\right)\bigg|_{s=\sigma_0},\\
c_{j,k}&\colonequals\frac{(-1)^{k-j}(\sigma_0-1)}{(\beta-k-1)!\sigma_0^{k-j+1}},
\end{align*}
and the implicit constant in the error term depends at most on $\beta,\sigma_0,\theta_0,\delta,\epsilon_0$. Notably, one gains an asymptotic for $S(x)$ with a power-saving error term uniformly in $B$, in contrast to what is provided by (\ref{Equ:S(x)}). Furthermore, suppose that there exists a constant $\lambda>0$ such that  $\alpha(p^\nu)=O\left(\left(\lambda p^{\sigma_0-1}\right)^{\nu}\right)$ for all prime powers $p^{\nu}$. Let 
\[F(s,a)\colonequals\prod_{p\mid a}\left(1-\left(\sum_{\nu=0}^{\infty}\alpha(p^{\nu})p^{-\nu s}\right)^{-1}\right)\]
for $s\in\C$ with $\sigma\ge\theta_0$ and squarefree $a\in\N$. When $s=\sigma_0$, this definition coincides with the one introduced in Lemma \ref{lemMV2}. As in the proof of Lemma \ref{lemMV2}, it is not hard to show that 
\begin{equation}\label{Equ:F(s,p)}
F(s,p)=\frac{\alpha(p)}{p^s}+O\left(\frac{\alpha(p)^2}{p^{2\sigma}}+\frac{1}{p^{2(\sigma-\sigma_0+1)}}\right)
\end{equation}
for all $s\in\C$ with $\sigma\geq\theta_0$ and all sufficiently large $p$. In addition, we observe that 
\begin{align*}
\sum_{\substack{n=1\\a\mid n}}^{\infty}\frac{\alpha(n)}{n^s}&=\left(\prod_{p\nmid a}\sum_{\nu=0}^{\infty}\alpha(p^{\nu})p^{-\nu s}\right)\left(\prod_{p\mid a}\sum_{\nu=1}^{\infty}\alpha(p^{\nu})p^{-\nu s}\right)\\
&=F(s)\left(\prod_{p\mid a}\sum_{\nu=0}^{\infty}\alpha(p^{\nu})p^{-\nu s}\right)^{-1}\left(\prod_{p\mid a}\sum_{\nu=1}^{\infty}\alpha(p^{\nu})p^{-\nu s}\right)=F(s)F(s,a)
\end{align*}
for $s\in\C$ with $\sigma>\sigma_0$ and squarefree $a\in\N$. Applying (\ref{Equ:S(x) via S-D}) to the above Dirichlet series expansion of $F(s)F(s,a)$ and using (\ref{Equ:F(s,p)}) to obtain upper bounds for $H_{\beta}(s)F(s,a)$ uniformly in $\sigma\ge\theta_0$, we see that there exist constants $\epsilon\in(0,1)$, $Q_0\ge2$, and $d_{j,k}\in\R$, where $0\leq j<k<\beta$, such that 
\begin{align}\label{Equ:sum_{n<=x,a|n}alpha(n) via S-D}
\sum_{\substack{n\leq x\\a\mid n}}\alpha(n)&=\frac{\mu_0(\beta)F(\sigma_0,a)}{(\beta-1)!\sigma_0}x^{\sigma_0}(\log x)^{\beta-1}+x^{\sigma_0}(\log x)^{\beta-1}\sum_{k=1}^{\beta-1}\sum_{j=0}^{k}d_{j,k}\frac{F^{(j)}(\sigma_0,a)}{(\log x)^k}\nonumber\\
&\hspace*{3mm}+O\left(B2^{O(\omega(a))}a^{\sigma_0-1}(x/a)^{\theta}\right)
\end{align}
uniformly for all $x\geq3$, $\theta\in(\sigma_0-\epsilon,\sigma_0)$ and square-free $a\in\N$ with $P^-(a)>Q_0$, where $F^{(j)}(\sigma_0,a)$ is the $j$th order derivative of $F(s,a)$ with respect to $s$ evaluated at $s=\sigma_0$. Again, one may compare this result with Lemma \ref{lemMV2}. 
\par Now, if $f\colon\N\to\R$ is a strongly additive function with $|f(p)|\le M$ for all primes $p$, where $M>0$ is constant, and if $0<h_0<(3/2)^{2/3}$ is fixed but arbitrary, then we obtain, by using (\ref{Equ:sum_{n<=x,a|n}alpha(n) via S-D}) as a substitute for Lemma \ref{lemMV2} and arguing as before with the adoption of the technique used in \cite[Section 4.2]{KMS}, that
\[M(x;m)=C_mB(x)^{\frac{m}{2}}\left(\chi_m+O\left(\frac{m^{\frac{3}{2}}}{\sqrt{B(x)}}\right)\right)\]
uniformly for all sufficiently large $x$ and all $1\le m\le h_0(B(x)/M^2)^{1/3}$, provided that $B(x)\to\infty$ as $x\to\infty$. Analogously, let $f\colon\N\to\R$ is strongly additive such that $f(p)=O(\sqrt{B(p)})$ for all primes $p$, $B(x)\to\infty$ as $x\to\infty$, and 
\[\sum_{\substack{p\le x\\|f(p)|>\epsilon \sqrt{B(x)}}}\alpha(p)\frac{f(p)^2}{p}=o(B(x))\]
for any given $\epsilon>0$. Then $M(x;m)=(\mu_m+o(1))B(x)^{\frac{m}{2}}$ for every fixed $m\in\N$. These results supplement Theorems \ref{thmMT2} and \ref{thmWeightedDH}. It may be worth pointing out that in the proofs of these results one can simply take $z=x^{1/v}$ with $v$ being a suitable constant multiple of $m$. We invite the reader to fill in the details.
\par One of the key ingredients in the proof of Theorem \ref{thmMT2} is an asymptotic formula for 
\[\sum_{\substack{n\le x\\d\mid n}}\alpha(n),\]
which is provided by Lemma \ref{lemMV2}. More generally, let $\mathcal{A}(x)=\{a_n\}_{n\le x}$ be a non-decreasing sequence of positive integers, and suppose that 
\begin{equation}\label{Equ:A_d}
\mathcal{A}_{d,\alpha}(x)\colonequals\sum_{\substack{n\le x\\d\mid a_n}}\alpha(n)=\rho(d)S(x)+r_d(x)
\end{equation}
for square-free integers $d\in\N$, where $\rho\colon\N\to[0,1]$ is a multiplicative function, and $r_d(x)$ is a remainder term which is expected to be small for all $d$ or small on average over $d$. Here, $\rho(d)$ can be viewed as the density of the set $\{n\in\N\colon d\mid a_n\}$ with respect to the probability measure induced by $\alpha$. In this sieve-theoretic setting one can derive, without much difficulty, an analogue of \cite[Proposition 4]{GS}. It may be of interest to determine if such an analogue can be used to obtain general weighted Erd\H{o}s--Kac theorems for various interesting sequences $\{a_n\}$ studied relatively recently, including $g(p_n)$, $\varphi(n)$, the Carmichael function $\lambda(n)$, and the aliquot sum $s(n)\colonequals\sigma(n)-n$, where $g\in\Z[x]$ is an irreducible polynomial, $p_n$ is the $n$th prime, and $\lambda(n)$ denotes the exponent of $(\Z/n\Z)^{\times}$
(see \cite{Hal3}, \cite{EP,EGPS} and \cite{PT}). Besides, the same approach may also be adapted to prove results of weighted Erd\H{o}s--Kac type for short intervals as well as in the function field setting. We hope to return to these and other related problems in the future.

\medskip
{\noindent\bf Acknowledgment.} The author thanks his advisor Carl Pomerance for stimulating discussions and helpful comments. He is also grateful to Prof. Paul Pollack for his cheerful encouragement and valuable feedback. Finally, he would like to express his gratitude to the anonymous referee for a careful reading of the manuscript and for detailed comments and suggestions which helped improve the paper considerably.

\medskip

\bibliography{bibliography}
\end{document}

%% file: Macros.tex
\usepackage{amssymb}
\usepackage{amsthm}
\usepackage{amsmath}
\usepackage{amsxtra}
\usepackage{breqn}
\usepackage{mathrsfs}
\usepackage[all]{xy}
\usepackage{colonequals}
\usepackage{enumitem}
\usepackage{fullpage}
\usepackage{stmaryrd}

\numberwithin{equation}{section}
\theoremstyle{plain}
\newtheorem{thm}{Theorem}[section]
\newtheorem{lem}[thm]{Lemma} 

\newtheorem{cor}[thm]{Corollary}

\theoremstyle{definition}
\newtheorem*{ack*}{Acknowledgment}
\theoremstyle{remark}

\newtheorem{rmk}{Remark}[section]

\DeclareMathOperator{\Li}{Li}

\DeclareMathOperator{\Res}{Res}

\newcommand{\C}{\mathbb C}

\newcommand{\N}{\mathbb N}
\newcommand{\PP}{\mathbb P}
\newcommand{\Q}{\mathbb Q}
\newcommand{\R}{\mathbb R}
\newcommand{\Z}{\mathbb Z}

\newcommand{\rad}{\text{rad}}

\input cyracc.def


%% file: WEK.bbl
\begin{bibdiv}
\begin{biblist}

\bib{Alla}{incollection}{
      author={Alladi, K.},
       title={A study of the moments of additive functions using {L}aplace
  transforms and sieve method},
   booktitle={in \textit{{N}umber {T}heory ({O}otacamund, 1984)}, {Lecture
  Notes in Math., vol. 1122, Springer, Berlin, 1985}},
       pages={1\ndash 37},
}

\bib{Billingsley}{article}{
      author={Billingsley, P.},
       title={{On the central limit theorem for the prime divisor functions}},
        date={1969},
     journal={Amer. Math. Monthly},
      volume={76},
       pages={132\ndash 139},
}

\bib{BillingsleyText}{book}{
      author={Billingsley, P.},
       title={Probability and measure},
      series={3rd. ed., Wiley Ser. Probab. Stat.},
   publisher={Wiley-Intersci. Publ., John Wiley \& Sons, Inc.},
     address={New York},
        date={1995},
}

\bib{BFI}{article}{
      author={Bombieri, E.},
      author={Friedlander, J.~B.},
      author={Iwaniec, H.},
       title={{Primes in arithmetic progressions to large moduli}},
        date={1986},
     journal={Acta Arith.},
      volume={156},
      number={3--4},
       pages={203\ndash 251},
}

\bib{BT}{article}{
      author={de~la Bret\`{e}che, R.},
      author={Tenenbaum, G.},
       title={{Remarks on the Selberg--Delange method}},
        date={2021},
     journal={Acta Arith.},
      volume={200},
      number={4},
       pages={349\ndash 369},
}

\bib{Delange1}{article}{
      author={Delange, H.},
       title={{Sur le nombre des diviseurs premiers de $n$}},
        date={1953},
     journal={C. R. Acad. Sci. (Paris)},
      volume={237},
       pages={542\ndash 544},
}

\bib{Delange2}{article}{
      author={Delange, H.},
       title={{Sur un th\'{e}or\`{e}me d'Erd\H{o}s et Kac}},
        date={1956},
     journal={Acad. Roy. Belg. Bull. Cl. Sci. (5)},
      volume={42},
       pages={130\ndash 144},
}

\bib{DH}{article}{
      author={Delange, H.},
      author={Halberstam, H.},
       title={{A note on additive functions}},
        date={1957},
     journal={Pacific J. Math.},
      volume={7},
       pages={1551\ndash 1556},
}

\bib{Deligne}{article}{
      author={Deligne, P.},
       title={{La conjecture de Weil. I.}},
        date={1974},
     journal={Inst. Hautes \'{E}tudes Sci. Publ. Math.},
      number={43},
       pages={273\ndash 307},
}

\bib{Dur}{book}{
      author={Durrett, R.},
       title={Probability: theory and examples},
      series={Camb. Ser. Stat. Probab. Math., vol. 49},
   publisher={Cambridge University Press},
     address={Cambridge},
        date={2019},
}

\bib{EG}{article}{
      author={Elboim, D.},
      author={Gorodetsky, O.},
       title={{Multiplicative arithmetic functions and the generalized Ewens
  measure}},
     journal={preprint (2022), to appear in Israel J. Math.},
  pages={\href{https://arxiv.org/abs/1909.00601}{\fontfamily{cmtt}\selectfont
  arXiv:1909.00601}},
}

\bib{Elliott1}{article}{
      author={Elliott, P. D. T.~A.},
       title={{A central limit theorem for Ramanujan's tau function}},
        date={2012},
     journal={Ramanujan J.},
      volume={29},
      number={1--3},
       pages={145\ndash 161},
}

\bib{Elliott2}{article}{
      author={Elliott, P. D. T.~A.},
       title={{Corrigendum: central limit theorems for classical cusp forms}},
        date={2015},
     journal={Ramanujan J.},
      volume={36},
      number={1--2},
       pages={99\ndash 102},
}

\bib{EGPS}{incollection}{
      author={Erd\H{o}s, P.},
      author={Granville, A.},
      author={Pomerance, C.},
      author={Spiro, C.},
       title={On the normal behavior of the iterates of some arithmetic
  functions},
   booktitle={in \textit{Analytic {N}umber {T}heory, {P}roc. {C}onf. in honor
  of {P}aul {T}. {B}ateman}, {B.\,C. Berndt et al. (eds.), Progr. Math., vol.
  85, Birkh\"auser, Boston, 1990}},
       pages={165\ndash 204},
}

\bib{EK}{article}{
      author={Erd\H{o}s, P.},
      author={Kac, M.},
       title={{The Gaussian law of errors in the theory of additive number
  theoretic functions}},
        date={1940},
     journal={Amer. J. Math.},
      volume={62},
       pages={738\ndash 742},
}

\bib{EP}{article}{
      author={Erd\H{o}s, P.},
      author={Pomerance, C.},
       title={{On the normal number of prime factors of $\varphi(n)$}},
        date={1985},
     journal={Rocky Mountain J. Math.},
      volume={15},
      number={2},
       pages={343\ndash 352},
}

\bib{GS}{incollection}{
      author={Granville, A.},
      author={Soundararajan, K.},
       title={Sieving and the {E}rd{\H o}s--{K}ac theorem},
   booktitle={in \textit{Equidistribution in {N}umber {T}heory, an
  introduction}, {A. Granville and Z. Rudnick (eds.), NATO Sci. Ser. II Math.
  Phys. Chem., vol. 237, Springer, Dordrecht, 2007}},
       pages={15\ndash 27},
}

\bib{Hal1}{article}{
      author={Halberstam, H.},
       title={{On the distribution of additive number-theoretic functions}},
        date={1955},
     journal={J. London Math. Soc.},
      volume={30},
       pages={43\ndash 53},
}

\bib{Hal2}{article}{
      author={Halberstam, H.},
       title={{On the distribution of additive number-theoretic functions
  (II)}},
        date={1956},
     journal={J. London Math. Soc.},
      volume={31},
       pages={1\ndash 14},
}

\bib{Hal3}{article}{
      author={Halberstam, H.},
       title={{On the distribution of additive number-theoretic functions
  (III)}},
        date={1956},
     journal={J. London Math. Soc.},
      volume={31},
       pages={14\ndash 27},
}

\bib{HR}{article}{
      author={Hardy, G.~H.},
      author={Ramanujan, S.},
       title={{The normal number of prime factors of a number $n$}},
        date={1917},
     journal={Quart. J. Pure. Appl. Math.},
      volume={48},
       pages={76\ndash 92},
}

\bib{HW}{book}{
      author={Hardy, G.~H.},
      author={Wright, E.~M.},
       title={An introduction to the theory of numbers},
      series={6th. ed.},
   publisher={Oxford Univ. Press},
     address={Oxford},
        date={2008},
        note={Revised by D. R. Heath-Brown and J. H. Silverman. With a foreword
  by Andrew Wiles},
}

\bib{Kac}{article}{
      author={Kac, M.},
       title={{Probability methods in some problems of analysis and number
  theory}},
        date={1949},
     journal={Bull. Amer. Math. Soc.},
      volume={55},
       pages={641\ndash 665},
}

\bib{KMS}{article}{
      author={Khan, R.},
      author={Milinovich, M.~B.},
      author={Subedi, U.},
       title={{A weighted version of the Erd\H{o}s--Kac theorem}},
        date={2022},
     journal={J. Number Theory},
      volume={239},
       pages={1\ndash 20},
}

\bib{LeVeque}{article}{
      author={LeVeque, Wm.~J.},
       title={{On the size of certain number-theoretic functions}},
        date={1949},
     journal={Trans. Amer. Math. Soc.},
      volume={66},
       pages={440\ndash 463},
}

\bib{Ngu}{article}{
      author={Nguyen, D.~T.},
       title={{Generalized divisor functions in arithmetic progressions: I}},
        date={2021},
     journal={J. Number Theory},
      volume={227},
       pages={30\ndash 93},
}

\bib{PT}{article}{
      author={Pollack, P.},
      author={Troupe, L.},
       title={{Sums of proper divisors follow the Erd\H{o}s--Kac law}},
        date={2023},
     journal={Proc. Amer. Math. Soc.},
      volume={151},
      number={3},
       pages={977\ndash 988},
}

\bib{RT}{article}{
      author={R\'{e}nyi, A.},
      author={Tur\'{a}n, P.},
       title={{On a theorem of Erd\"{o}s--Kac}},
        date={1958},
     journal={Acta Arith.},
      volume={4},
       pages={71\ndash 84},
}

\bib{Sha}{article}{
      author={Shapiro, H.~N.},
       title={{Distribution functions of additive arithmetic functions}},
        date={1956},
     journal={Proc. Nat. Acad. Sci. U.S.A.},
      volume={42},
       pages={426\ndash 430},
}

\bib{Song}{article}{
      author={Song, J.~M.},
       title={{Sums of nonnegative multiplicative functions over integers
  without large prime factors I}},
        date={2001},
     journal={Acta Arith.},
      volume={97},
      number={4},
       pages={329\ndash 351},
}

\bib{TenenText}{book}{
      author={Tenenbaum, G.},
       title={Introduction to analytic and probabilistic number theory},
      series={3rd. ed., Grad. Stud. Math.},
   publisher={Amer. Math. Soc.},
     address={Providence, RI.},
        date={2015},
      volume={163},
        note={Translated from the 2008 French edition by Patrick D. F. Ion},
}

\bib{Tenen1}{article}{
      author={Tenenbaum, G.},
       title={{Moyennes effectives de fonctions multiplicatives complexes}},
        date={2017},
     journal={Ramanujan J.},
      volume={44},
      number={3},
       pages={641\ndash 701},
}

\bib{Tenen2}{article}{
      author={Tenenbaum, G.},
       title={{Correction to: moyennes effectives de fonctions multiplicatives
  complexes}},
        date={2020},
     journal={Ramanujan J.},
      volume={53},
      number={1},
       pages={243\ndash 244},
}

\bib{Turan}{article}{
      author={Tur\'an, P.},
       title={{On a theorem of Hardy and Ramanujan}},
        date={1934},
     journal={Ramanujan J.},
      volume={9},
      number={4},
       pages={274\ndash 276},
}

\end{biblist}
\end{bibdiv}
